\documentclass[final,3p]{elsarticle}
\makeatletter 
\def\ps@pprintTitle{%
 \def\@oddfoot{\footnotesize\itshape
       Published in \ifx\@journal\@empty Elsevier
       \else\@journal\fi\hfill February 1, 2020}}
\makeatother

\graphicspath{{./graphics/}}

\usepackage{amsmath,amsfonts,amssymb} 
\usepackage{mathtools} 
\usepackage{upgreek} 
\usepackage{easybmat} 
\usepackage{fixmath} 
\usepackage[font=footnotesize,list=true]{subcaption} 
\usepackage{booktabs} 
\usepackage[exponent-product=\cdot]{siunitx} 
\renewcommand{\ang}[1]{{{#1}^\circ}}
\usepackage{bm} 
\usepackage{xcolor}
\usepackage{microtype}
\usepackage{lmodern}

\newcommand{\graphicsFolder}{graphics}

\renewcommand{\vec}[1]{\mathbold #1} 
\newcommand{\besselj}{\mathrm{j}} 	
\newcommand{\legendre}{\mathrm{P}} 	
\newcommand{\hankel}{\mathrm{h}}   	
\newcommand{\Diff}{\upDelta} 

\newcommand{\idiff}{\, \mathrm{d}} 

\newcommand{\R}{\mathbb{R}}

\newcommand{\Z}{\mathbb{Z}}
\newcommand{\N}{\mathbb{N}}
\newcommand{\PI}{\uppi}
\newcommand{\euler}{\mathrm{e}}
\newcommand{\imag}{\mathrm{i}}

\providecommand*{\pderiv}[3][]{\frac{\partial^{#1}#2}{\partial #3^{#1}}}

\renewcommand{\geq}{\geqslant}
\renewcommand{\leq}{\leqslant}

\DeclareMathSymbol{\GAMMA}{\mathalpha}{operators}{0}

\DeclareMathOperator{\fourier}{\mathcal{F}}
\DeclareMathOperator{\TS}{TS}

\renewcommand\Re{\operatorname{Re}}

\renewcommand{\Xi}{{\vec{t}_1}}

\newcommand{\energyNorm}[2]{%
  {\left\vert\kern-0.25ex\left\vert\kern-0.25ex\left\vert #1 
    \right\vert\kern-0.25ex\right\vert\kern-0.25ex\right\vert}_{#2}
}

\let\originalleft\left
\let\originalright\right
\renewcommand{\left}{\mathopen{}\mathclose\bgroup\originalleft}
\renewcommand{\right}{\aftergroup\egroup\originalright}
\newcommand{\round}[1]{\ensuremath{\left\lfloor{#1}\right\rceil}}

\newcommand{\ceil}[1]{\ensuremath{\left\lceil{#1}\right\rceil}}
\newcommand{\st}{\,:\,}

\usepackage{xspace}
\newcommand{\MATLAB}{\textsc{Matlab}\xspace}

\def\ccpp{{C/C\nolinebreak[4]\hspace{-.05em}\raisebox{.4ex}{\tiny\bf ++}}\xspace}
\newcommand{\COMSOL}{\textsc{Comsol} Multiphysics\textregistered\xspace}

\usepackage{hyperref} 
\hypersetup{hidelinks,colorlinks=true,citecolor=blue,linkcolor=blue,urlcolor=blue} 
\usepackage[noabbrev]{cleveref} 
\usepackage{longtable}
\usepackage{enumitem}

\newtheorem{theorem}{Theorem}

\journal{Computer Methods in Applied Mechanics and Engineering}
\begin{document}

\begin{frontmatter}

\title{Isogeometric boundary element method for acoustic scattering by a submarine}

\author[venas]{Jon Vegard Ven{\aa}s\texorpdfstring{\corref{cor}}{}}
\ead{Jon.Venas@ntnu.no}

\author[venas]{Trond Kvamsdal}
\ead{Trond.Kvamsdal@ntnu.no}


\address[venas]{Department of Mathematical Sciences, Norwegian University of Science and Technology, 7034 Trondheim, Norway}
\cortext[cor]{Corresponding author.}

\begin{abstract}
Isogeometric analysis with the boundary element method (IGABEM) has recently gained interest. In this paper, the approximability of IGABEM on 3D acoustic scattering problems will be investigated and a new improved BeTSSi submarine will be presented as a benchmark example. Both Galerkin and collocation are considered in combination with several boundary integral equations (BIE). In addition to the conventional BIE, regularized versions of this BIE will be considered. Moreover, the hyper-singular BIE and the Burton--Miller formulation are also considered. A new adaptive integration routine is presented, and the numerical examples show the importance of the integration procedure in the boundary element method. The numerical examples also include comparison between standard BEM and IGABEM, which again verifies the higher accuracy obtained from the increased inter-element continuity of the spline basis functions. One of the main objectives in this paper is benchmarking acoustic scattering problems, and the method of manufactured solution will be used frequently in this regard.
\end{abstract}


\end{frontmatter}
\section{Introduction}
\begin{figure}
	\centering
	\begin{subfigure}[b]{0.3\textwidth}
		\centering
		\includegraphics[width=0.9\textwidth]{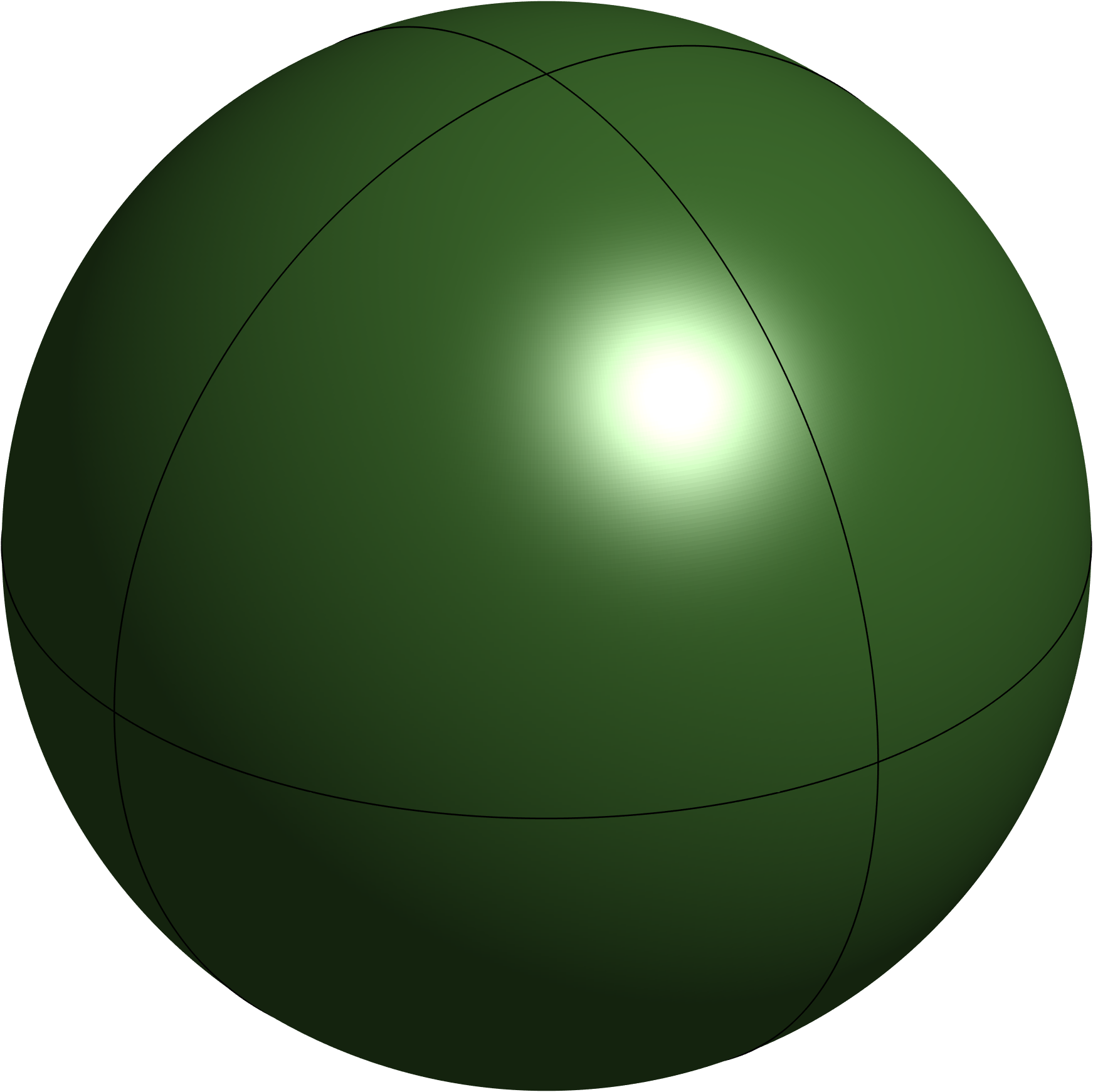}
		\caption{A sphere with 8 elements.}
		\label{Fig3:sphere}
	\end{subfigure}
	~    
	\begin{subfigure}[b]{0.35\textwidth}
		\centering
		\includegraphics[width=0.9\textwidth]{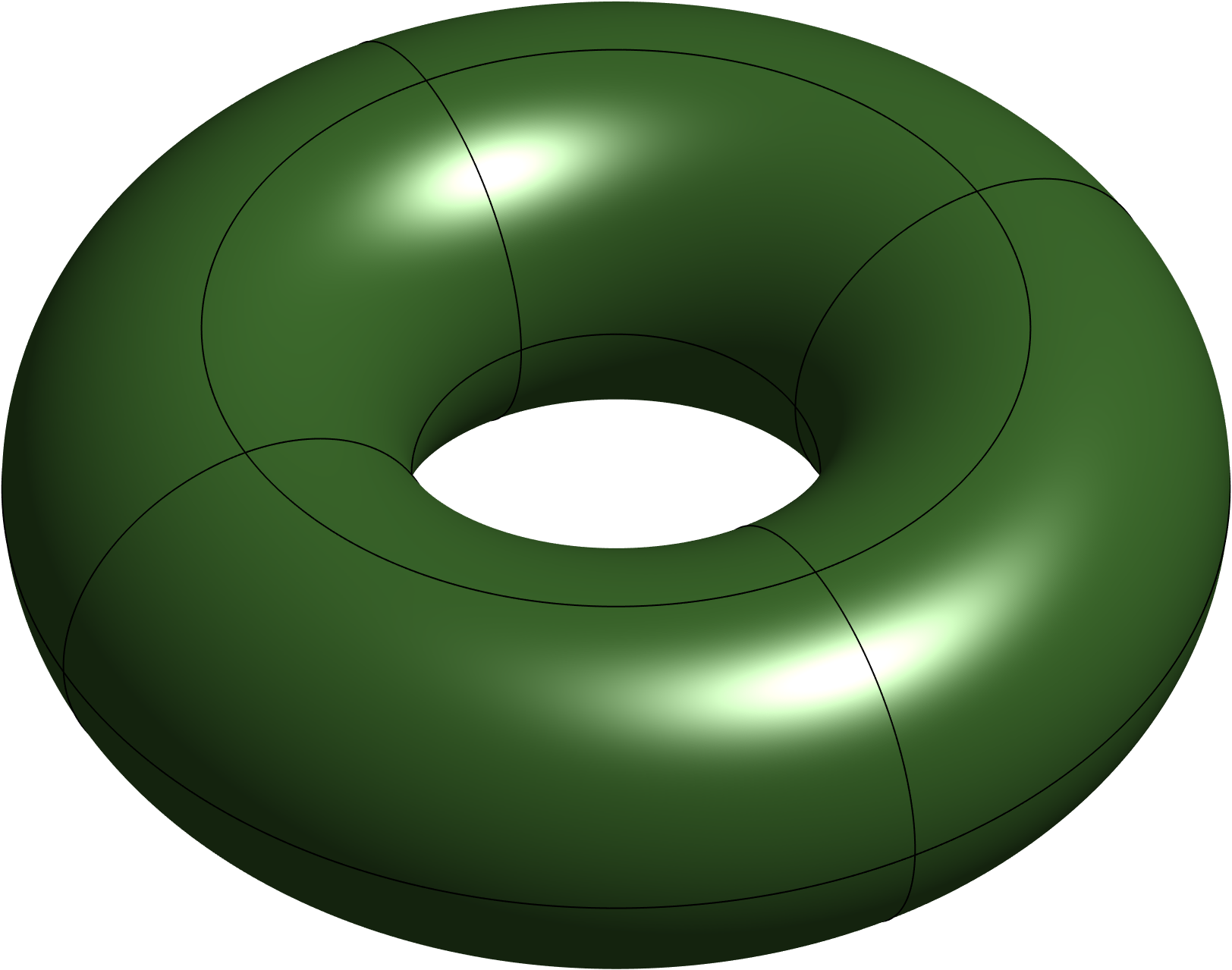}
		\caption{A torus with 16 elements.}
		\label{Fig3:Torus}
	\end{subfigure}
	\caption{Examples of exact NURBS geometries of second degree.}
	\label{Fig3:NURBSexamples}
\end{figure}

Isogeometric analysis (IGA) was introduced in 2005 by Hughes et al.~\cite{Hughes2005iac}, followed by the book~\cite{Cottrell2009iat} in 2009. Since then, IGA has received a great deal of attention in the effort of bridging the gap between finite element analysis (FEA) and computer aided design (CAD) tools. The initial problem that sparked the IGA movement was the cumbersome mesh generating process when converting the design models from CAD into the FEA programs, and the analysis could often imply a rerun of this tedious process. The problem being that the geometry was represented differently in CAD and FEA. An example is the geometries illustrated in \Cref{Fig3:NURBSexamples} which can be represented exactly using NURBS but is outside the space of standard (Lagrangian) FEM geometries.
Using the same geometry representation as in CAD, IGA features exact geometry, which remains true in all mesh refinement procedures. Moreover, it turns out that using the non-uniform rational B-splines (NURBS) as basis functions not only for representing the geometry, but also the solution space, greatly enhances the numerical accuracy, see~\cite{BeiraodaVeiga2011sef} and~\cite{BeiraodaVeiga2014mao}. This motivates the use of IGA even further, as IGA enables control of the continuity of the basis function up to $C^{\check{p}-1}$ where $\check{p}$ is the polynomial degree (in contrast with the $C^0$-continuity restriction in classical FEA). For exterior problems, one can introduce an artificial boundary to obtain a bounded domain introducing the difficulty of surface-to-volume parametrization. The boundary element method (BEM) avoids this issue entirely as it only relies on a computational domain on the surface of the scatterer. Moreover, solid domains are usually represented by surfaces in CAD-systems, such that if modeling of an elastic scatterer is required, the BEM solves this problem as well without the need of surface-to-volume parametrization. This then represents an even further improvement of the quality of the design-analysis bridging development.

\begin{figure}
	\centering
	\includegraphics[scale=1]{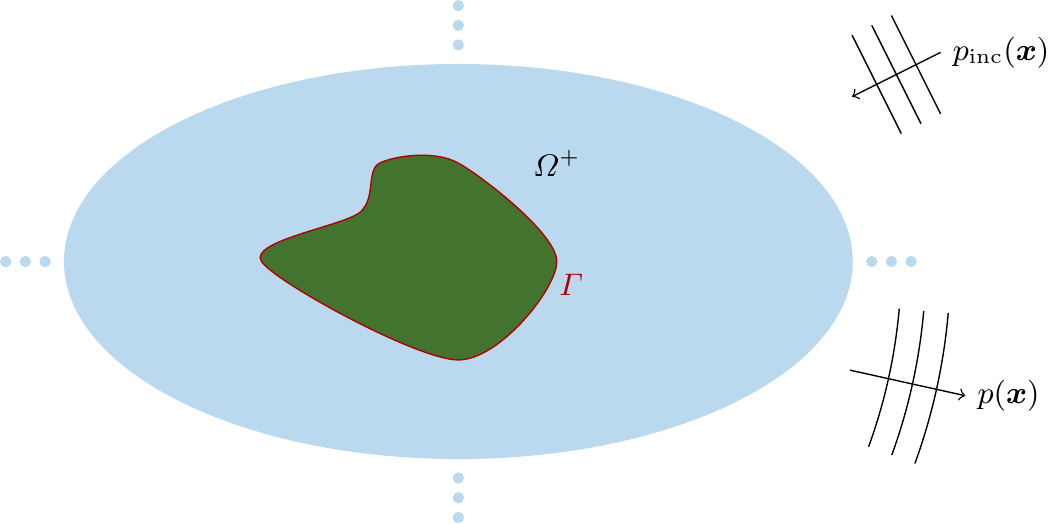}
	\caption{Illustration of the physical problem. A plane incident wave, $p_{\mathrm{inc}}(\vec{x})$, is scattered by the scatterer, represented by the closed boundary $\Gamma$, in an unbounded domain, $\Omega^+\subset\R^d$, resulting in the scattered pressure, $p(\vec{x})$.}
	\label{Fig3:physicalProblem}
\end{figure}
This work is only concerned with 3D acoustic scattering (with $d=3$). The main objective is scattering by plane waves, $p_{\rm inc}$, as illustrated in \Cref{Fig3:physicalProblem}. In scattering problems, it is often of interest to compute the target strength, $\TS$, of the scatterer in the far field. As an application of this work, the target strength is the quantity of interest for the acoustical aspects of constructing a submarine and is for this reason investigated in this work.

Assuming harmonic time dependency, all time dependent functions may be written as $\breve{F}=\breve{F}(\vec{x},t) = F(\vec{x})\euler^{-\imag\omega t}$ where $\omega$ is the angular frequency and $\imag = \sqrt{-1}$ the imaginary unit. This enables us to model the pressure $p$ in the fluid with the Helmholtz equation given by
\begin{equation}\label{Eq3:HelmholtzEquationIntro}
	\nabla^2 p + k^2 p = 0
\end{equation}
with the wave number $k=\frac{\omega}{c_{\mathrm{f}}}$ (where $c_{\mathrm{f}}$ is the wave speed in the fluid\footnote{Throughout this work we shall use $c_{\mathrm{f}}=\SI{1500}{m/s}$.}). Other important quantities include the frequency $f=\frac{\omega}{2\PI}$ and the wavelength $\lambda = \frac{2\PI}{k}$.

Some literature already exists for solving acoustic problems using IGABEM including \cite{Simpson2014aib,Keuchel2017eoh,Peake2013eib,Peake2014eai,Peake2015eib,Coox2017aii, Dolz2016aib,Dolz2018afi,Sun2019dib,Wu2020iib}. Arguably there is a lack of work in the approximability for IGABEM simulations for more complex geometries, and one of the aims of this work is to contribute to fill this gap.

The exterior Helmholtz problem is presented in \Cref{Sec3:exteriorHelmholtz}, and the corresponding boundary integral equations are given in \Cref{Sec3:BIE}. Discretization of these integral equations either with the use of collocation or a Galerkin approach yields the boundary element method which is presented in \Cref{Sec3:BEM}. The weakly singular boundary integral equation requires care when using numerical quadrature and is discussed in \Cref{Sec3:numericalQuad}. In \Cref{Sec3:resultsDisc} the results for several benchmark problems are presented. Not only are these benchmark problems important in bug testing for code development, but it is also important to establish reliable results for several geometries ranging in complexity. Finally, conclusions and suggested future work can be found in \Cref{Sec3:conclusions}.

\section{Helmholtz problems}
\label{Sec3:exteriorHelmholtz}
The Helmholtz problem is given by
\begin{alignat}{3}
	\nabla^2 p + k^2 p &= 0 	\quad &&\text{in}\quad \Omega,\label{Eq3:HelmholtzEqn}\\
	\partial_n p &= g						&&\text{on}\quad \Gamma,\label{Eq3:HelmholtzEqnNeumannCond}
\end{alignat}
where $\partial_n$ denotes the partial derivative in the normal direction, $\vec{n}$, on the surface $\Gamma$. Throughout this work, $\vec{n}$ is always pointing ``into'' $\Omega^+$. If $\Omega=\Omega^-$ is inside a closed boundary $\Gamma$, the problem is referred to as an interior problem. If, on the other hand, $\Omega=\Omega^+$ is the unbounded domain outside $\Gamma$ (as illustrated in \Cref{Fig3:physicalProblem}), the problem is referred to as an exterior problem where we must impose the Sommerfeld condition~\cite{Sommerfeld1949pde} 
\begin{equation}\label{Eq3:sommerfeldCond}
	\pderiv{p}{r}-\imag k p = o\left(r^{-1}\right)\quad \text{with}\quad r=|\vec{x}|\
\end{equation}
in order to restrict the field in the limit $r\to\infty$ uniformly in $\hat{\vec{x}}=\frac{\vec{x}}{r}$, such that no waves originate from infinity (to obtain uniqueness of the solution $p$).

A common approach for solving unbounded scattering problems with the FEM is to introduce an artificial boundary that encloses the scatterer. On the artificial boundary some sort of absorbing boundary condition (ABC) is prescribed. The problem is then reduced to a finite domain problem, and the bounded domain between the scatterer and the artificial boundary can be discretized with finite elements. Several methods exist for handling the exterior Helmholtz problem (on unbounded domain), including
\begin{itemize}
	\item the perfectly matched layer (PML) method after B{\'e}renger~\cite{Berenger1994apm,Berenger1996pml}
	\item the boundary element method~\cite{Sauter2011bem,Schanz2007bea,Marburg2008cao,Chandler_Wilde2012nab}
	\item Dirichlet to Neumann-operators (DtN-operators)~\cite{Givoli2013nmf}
	\item local differential ABC operators~\cite{Shirron1995soe,Bayliss1982bcf,Hagstrom1998afo,Tezaur2001tdf}
	\item the infinite element method.~\cite{Bettess1977ie,Bettess1977dar}
\end{itemize}
Due to the complexity of the BeTSSi geometry considered in this work we conveniently consider the boundary element method to solve the Helmholtz problem in order to avoid the surface-to-volume parametrization discussed in the introduction.

The Neumann condition (in \Cref{Eq3:HelmholtzEqnNeumannCond}), given by the function $g$, will in the case of rigid scattering be given in terms of the incident wave $p_{\mathrm{inc}}$. Zero displacement of the fluid normally on the scatterer (rigid scattering) implies that $\partial_n(p+p_{\mathrm{inc}}) = 0$, and hence
\begin{equation}
	g = -\pderiv{p_{\mathrm{inc}}}{n}.
\end{equation}
Plane incident waves (with amplitude $P_{\mathrm{inc}}$) traveling in the direction $\vec{d}_{\mathrm{s}}$ can be written as
\begin{equation}\label{Eq3:p_inc}
	p_{\mathrm{inc}} = P_{\mathrm{inc}}\euler^{\imag k\vec{d}_{\mathrm{s}}\cdot\vec{x}}.
\end{equation}
The normal derivative on the surface of any smooth geometry may then be computed by
\begin{align}
	\pderiv{p_{\mathrm{inc}}}{n} &= \vec{n}\cdot\nabla p_{\mathrm{inc}} = \imag k\vec{d}_{\mathrm{s}}\cdot\vec{n} p_{\mathrm{inc}}.
\end{align}

\subsection{Far field pattern}
If the field at the scatterer is known, one can compute the solution in the exterior domain, $\Omega^+$, using the following integral solution (cf.~\cite[Theorem 2.21]{Chandler_Wilde2012nab})
\begin{equation}\label{Eq3:KirchhoffIntegral}
	p(\vec{x}) = \int_{\Gamma}\left[ p(\vec{y})\pderiv{\Phi_k(\vec{x},\vec{y})}{n(\vec{y})} - \Phi_k(\vec{x},\vec{y})\pderiv{p(\vec{y})}{n(\vec{y})}\right]\idiff \Gamma(\vec{y}),\quad\vec{x}\in\Omega^+
\end{equation}
where $\vec{y}$ is a point on the surface $\Gamma$, $\vec{n}$ lies on $\Gamma$ pointing ``into'' $\Omega^+$ at $\vec{y}$, and $\Phi_k$ is the free space Green's function for the Helmholtz equation in \Cref{Eq3:HelmholtzEqn} given (in 3D) by
\begin{equation}\label{Eq3:FreeSpaceGrensFunction}
	\Phi_k(\vec{x},\vec{y}) = \frac{\euler^{\imag kR}}{4\PI R},\quad\text{where}\quad R = |\vec{x} - \vec{y}|.
\end{equation} 
For later convenience, we note that
\begin{align*}
	\pderiv{\Phi_k(\vec{x},\vec{y})}{n(\vec{y})} &= \frac{\Phi_k(\vec{x},\vec{y})}{R}(\imag kR-1)\pderiv{R}{n(\vec{y})}\\
	\pderiv{\Phi_k(\vec{x},\vec{y})}{n(\vec{x})} &= \frac{\Phi_k(\vec{x},\vec{y})}{R}(\imag kR-1)\pderiv{R}{n(\vec{x})}\\
	\frac{\partial^2\Phi_k(\vec{x},\vec{y})}{\partial\vec{n}(\vec{y})\partial\vec{n}(\vec{x})} &= -\frac{\Phi_k(\vec{x},\vec{y})}{R^2}\left[\vec{n}(\vec{x})\cdot\vec{n}(\vec{y})(\imag kR-1) + \left(k^2R^2+3(\imag kR - 1)\right)\pderiv{R}{n(\vec{x})}\pderiv{R}{n(\vec{y})}\right]
\end{align*}
where
\begin{equation*}
	\pderiv{R}{n(\vec{x})} = \frac{(\vec{x}-\vec{y})\cdot\vec{n}(\vec{x})}{R}\quad\text{and}\quad\pderiv{R}{n(\vec{y})} = -\frac{(\vec{x}-\vec{y})\cdot\vec{n}(\vec{y})}{R}.
\end{equation*}
The \textit{far field pattern} for the scattered pressure $p$, is defined by
\begin{equation}\label{Eq3:farfield}
	p_0(\hat{\vec{x}}) =  \lim_{r\to\infty} r \euler^{-\imag k r}p(r\hat{\vec{x}}),
\end{equation}
with $r = |\vec{x}|$ and $\hat{\vec{x}} = \vec{x}/|\vec{x}|$. Using the limits
\begin{equation}\label{Eq3:Phi_k_limits}
	\lim_{r\to\infty} r\euler^{-\imag k r}\Phi_k(r\hat{\vec{x}},\vec{y}) = \frac{1}{4\PI}\euler^{-\imag k \hat{\vec{x}}\cdot\vec{y}}\quad\text{and}\quad 
	\lim_{r\to\infty} r\euler^{-\imag k r}\pderiv{\Phi_k(r\hat{\vec{x}},\vec{y})}{n(\vec{y})} = -\frac{\imag k}{4\PI}\euler^{-\imag k \hat{\vec{x}}\cdot\vec{y}}\hat{\vec{x}}\cdot\vec{n}(\vec{y})
\end{equation}
the formula in \Cref{Eq3:KirchhoffIntegral} simplifies in the far field to (cf.~\cite[p. 32]{Ihlenburg1998fea})
\begin{equation}\label{Eq3:KirchhoffIntegralFarField}
	p_0(\hat{\vec{x}}) = -\frac{1}{4\PI}\int_{\Gamma}\left[ \imag k p(\vec{y})\hat{\vec{x}}\cdot\vec{n}(\vec{y}) + \pderiv{p(\vec{y})}{n(\vec{y})}\right]\euler^{-\imag k \hat{\vec{x}}\cdot\vec{y}}\idiff \Gamma(\vec{y}).
\end{equation}
From the far field pattern, the \textit{target strength}, $\TS$, can be computed. It is defined by
\begin{equation}\label{Eq3:TS}
	\TS = 20\log_{10}\left(\frac{|p_0(\hat{\vec{x}})|}{|P_{\mathrm{inc}}|}\right)
\end{equation}
where $P_{\mathrm{inc}}$ is the amplitude of the incident wave at the geometric center of the scatterer (i.e. the origin). Note that the $\TS$ is independent of $P_{\mathrm{inc}}$, which is a result of the linear dependency of the amplitude of the incident wave in scattering problems (i.e. doubling the amplitude of the incident wave will double the amplitude of the scattered wave).
\section{Boundary integral equations}
\label{Sec3:BIE}
We adopt the following notation from \cite{Chandler_Wilde2012nab}. The single- and double layer potential operator are given by
\begin{equation*}
	\mathcal{S}_k \phi(\vec{x}) = \int_\Gamma \Phi_k(\vec{x},\vec{y}) \phi(\vec{y}) \idiff\Gamma(\vec{y})\qquad \vec{x}\in\R^d\setminus\Gamma,
\end{equation*}
and
\begin{equation*}
	\mathcal{D}_k \phi(\vec{x}) = \int_\Gamma \pderiv{\Phi_k(\vec{x},\vec{y})}{n(\vec{y})} \phi(\vec{y}) \idiff\Gamma(\vec{y})\qquad \vec{x}\in\R^d\setminus\Gamma,
\end{equation*}
respectively. Here, the normal vector $\vec{n}$ at the surface $\Gamma$ always points from the interior domain $\Omega^-$ into the exterior domain $\Omega^+$.

For $D\subset\R^d$, define the spaces (for details see~\cite{Chandler_Wilde2012nab})
\begin{align*}
	L_{\mathrm{loc}}^2(D) &= \left\{u|_G\in L^2(G)\st \forall G\subset D,\, G \text{ bounded and measurable}\right\}\\
	H_{\mathrm{loc}}^1(D) &= \left\{u\in L_{\mathrm{loc}}^2(D)\st vu\in H^1(D),\, v\in C_{\mathrm{comp}}^\infty(\overline{D})\right\}\\
	H_{\mathrm{loc}}^1(D;\nabla) &= \left\{u\in L_{\mathrm{loc}}^2(D)\st \nabla u\in \left[L_{\mathrm{loc}}^2(D)\right]^d,\, \nabla^2 u\in L_{\mathrm{loc}}^2(D)\right\}\\
	H^s(\R^d) &= \left\{u\in L^2(\R^d)\st \fourier^{-1}\left[\left(1+|\vec{\xi}|^2\right)^s\fourier u\right]\in L^2(\R^d)\right\}\\
	H^s(D) &= \left\{u|_D\st u\in H^s(\R^d)\right\}\\
	H^s(\Gamma) &= \left\{\phi\in L^2(\Gamma)\st \phi_f\in H^s(\R^{d-1})\right\}
\end{align*}
with the Fourier transform
\begin{equation*}
	(\fourier u)(\vec{\xi}) = (2\PI)^{-d/2}\int_{\R^d}\euler^{-\imag \vec{x}\cdot \vec{\xi}}u(\vec{x})\idiff\Omega(\vec{x}),\quad\vec{\xi}\in\R^d.
\end{equation*}
By defining $\gamma^\pm$ to be the trace operator from $H^s(\Omega^\pm)\to H^{s-1/2}(\Gamma)$ for $\frac{1}{2}<s<\frac{3}{2}$ and $\partial_n^\pm$ to be the normal derivative from $H^1(\Omega^\pm; \nabla)\to H^{1/2}(\Gamma)$, we restate two important theorems for BEM analysis from \cite{Chandler_Wilde2012nab}, namely theorem 2.20 and 2.21:
\begin{theorem}\label{Thm:InteriorProblem}
	If $p\in H^1(\Omega^-)\cup C^2(\Omega^-)$ and, for some $k\geq 0$, $\nabla^2 p + k^2p = 0$ in $\Omega^-$, then
	\begin{equation*}
		\mathcal{S}_k\partial_n^- p(\vec{x}) - \mathcal{D}_k\gamma^- p(\vec{x}) = \begin{cases} p(\vec{x}), & x\in\Omega^-,\\
		0 & x\in\Omega^+.\end{cases}
\end{equation*}	 
\end{theorem}
\begin{theorem}\label{Thm:ExteriorProblem}
	If $p\in H_{\mathrm{loc}}^1(\Omega^+)\cup C^2(\Omega^+)$ and, for some $k > 0$, $\nabla^2 p + k^2 p = 0$ in $\Omega^+$ and $p$ satisfies the Sommerfeld radiation condition in $\Omega^+$, that is,
\begin{equation*}
	\pderiv{p(\vec{x})}{r}-\imag k p(\vec{x}) = o\left(r^{-\frac{d-1}{2}}\right)\quad r=|\vec{x}|
\end{equation*}	
	as $r\to\infty$ uniformly in $\hat{\vec{x}}=\frac{\vec{x}}{r}$, then
\begin{equation*}
		-\mathcal{S}_k\partial_n^+ p(\vec{x}) + \mathcal{D}_k\gamma^+ p(\vec{x}) = \begin{cases} p(\vec{x}), & x\in\Omega^+,\\
		0 & x\in\Omega^-.\end{cases}
\end{equation*}	 
\end{theorem}
The \textit{acoustic} single- and double layer potential operator are respectively given by
\begin{equation*}
	S_k \phi(\vec{x}) = \int_\Gamma \Phi_k(\vec{x},\vec{y}) \phi(\vec{y}) \idiff\Gamma(\vec{y})\qquad \vec{x}\in\Gamma
\end{equation*}
and
\begin{equation*}
	D_k \phi(\vec{x}) = \int_\Gamma \pderiv{\Phi_k(\vec{x},\vec{y})}{n(\vec{y})} \phi(\vec{y}) \idiff\Gamma(\vec{y})\qquad \vec{x}\in\Gamma
\end{equation*}
and the  acoustic adjoint double-layer operator and the hypersingular operator are respectively given by
\begin{equation*}
	D'_k \phi(\vec{x}) = \int_\Gamma \pderiv{\Phi_k(\vec{x},\vec{y})}{n(\vec{x})} \phi(\vec{y}) \idiff\Gamma(\vec{y})\qquad \vec{x}\in\Gamma
\end{equation*}
and
\begin{equation*}
	H_k \phi(\vec{x}) = \int_\Gamma \frac{\partial^2\Phi_k(\vec{x},\vec{y})}{\partial\vec{n}(\vec{y})\partial\vec{n}(\vec{x})} \phi(\vec{y}) \idiff\Gamma(\vec{y})\qquad \vec{x}\in\Gamma.
\end{equation*}
By following the notation in \cite[p. 117]{Chandler_Wilde2012nab} we let
\begin{equation*}
	M_k = \begin{bmatrix}
		D_k & -S_k\\
		H_k & -D'_k
	\end{bmatrix}\quad\text{and}\quad c^\pm p = \begin{bmatrix}
		\gamma^\pm p\\
		\partial_n^\pm p
	\end{bmatrix}
\end{equation*}
such that the boundary integral equations (BIE) for the exterior- and interior problem are respectively given by
\begin{equation*}
	\mp\frac{1}{2}c^\pm p = M_k c^\pm p.
\end{equation*}
We can write this more explicitly as
\begin{align*}
	&\mp \frac{1}{2}p(\vec{x}) + \int_\Gamma \pderiv{\Phi_k(\vec{x},\vec{y})}{n(\vec{y})} p(\vec{y}) \idiff\Gamma(\vec{y}) = \int_\Gamma \Phi_k(\vec{x},\vec{y}) \pderiv{p(\vec{y})}{n(\vec{y})} \idiff\Gamma(\vec{y})\\
	&\mp \frac{1}{2}\pderiv{p(\vec{x})}{n(\vec{x})} +  \int_\Gamma \frac{\partial^2\Phi_k(\vec{x},\vec{y})}{\partial\vec{n}(\vec{y})\partial\vec{n}(\vec{x})} p(\vec{y}) \idiff\Gamma(\vec{y}) = \int_\Gamma \pderiv{\Phi_k(\vec{x},\vec{y})}{n(\vec{x})} \pderiv{p(\vec{y})}{n(\vec{y})} \idiff\Gamma(\vec{y})
\end{align*}
for almost all $\vec{x}\in\Gamma$. These integral equations need a modification if $\Gamma$ is not smooth at $\vec{x}$. With the jump term defined as (cf. \cite{Hwang1997hbi})
\begin{equation}\label{Eq3:jumpTerm}
	C^\pm(\vec{x}) = \begin{cases}-\frac{1}{2}(1\pm 1) & \vec{x}\in\Omega^+\\-\frac{1}{2}(1\pm 1) - \int_\Gamma \pderiv{\Phi_0(\vec{x},\vec{y})}{n(\vec{y})}\idiff\Gamma(\vec{y}) &\vec{x}\in\Gamma\\
	\frac{1}{2}(1\mp 1)& \vec{x}\in\Omega^-\end{cases}
\end{equation}
the conventional BIE (CBIE) and hypersingular BIE (HBIE) are respectively given by 
\begin{align}\label{Eq3:CBIE}
	& C^\pm(\vec{x})p(\vec{x}) + \int_\Gamma \pderiv{\Phi_k(\vec{x},\vec{y})}{n(\vec{y})} p(\vec{y}) \idiff\Gamma(\vec{y}) = \int_\Gamma \Phi_k(\vec{x},\vec{y}) \pderiv{p(\vec{y})}{n(\vec{y})} \idiff\Gamma(\vec{y})\\\label{Eq3:HBIE}
	& C^\pm(\vec{x})\pderiv{p(\vec{x})}{n(\vec{x})} +  \int_\Gamma \frac{\partial^2\Phi_k(\vec{x},\vec{y})}{\partial\vec{n}(\vec{y})\partial\vec{n}(\vec{x})} p(\vec{y}) \idiff\Gamma(\vec{y}) = \int_\Gamma \pderiv{\Phi_k(\vec{x},\vec{y})}{n(\vec{x})} \pderiv{p(\vec{y})}{n(\vec{y})} \idiff\Gamma(\vec{y}).
\end{align}
Note that using the divergence theorem it is possible to show the following (cf. \cite[p. 126]{Sauter2011bem})
\begin{equation*}
	\int_\Gamma \pderiv{\Phi_0(\vec{x},\vec{y})}{n(\vec{y})}\idiff\Gamma(\vec{y})=\begin{cases}
		0 & \vec{x}\in\Omega^+\\
		-\frac{1}{2} & \vec{x}\in\Gamma,\text{ if } \Gamma \text{ is smooth at } \vec{x}\\
		-1 & \vec{x}\in\Omega^-.
	\end{cases}
\end{equation*}
This result may be generalized for the case that $\Gamma$ is not smooth at $\vec{x}$, namely in terms of the solid angle~\cite{Sun2015bri}
\begin{equation}\label{Eq3:Phi0SolidAngle}
	\int_\Gamma \pderiv{\Phi_0(\vec{x},\vec{y})}{n(\vec{y})}\idiff\Gamma(\vec{y}) = -\frac{c_0}{4\PI}
\end{equation}
where the solid angle $c_0$ can be computed by
\begin{equation*}
	c_0 = \lim_{\varepsilon\to 0^+}\frac{|\partial B_\varepsilon(\vec{x})\cap\Omega^-|}{\varepsilon^2}
\end{equation*}
where $B_\varepsilon(\vec{x})$ is a ball of radius $\varepsilon$ centered at $\vec{x}$. In other words, the integral in \Cref{Eq3:Phi0SolidAngle} is given by the negative relative size of the surface of a infinitesimal small sphere centered at $\vec{x}$ that is inside $\Omega^-$. This enables simple exact calculation of this integral for most standard geometries. For example, if $\Omega^-$ is a cube, the integral in \Cref{Eq3:Phi0SolidAngle} takes the value $-\frac14$ and $-\frac{1}{8}$ if $\vec{x}$ is at an edge or at a vertex, respectively. This can be used to test the numerical integration involved in solving BIEs.

Combining the CBIE in~\Cref{Eq3:CBIE} and the HBIE in~\Cref{Eq3:HBIE} yields the Burton--Miller (BM) formulation which can conceptually be written as
\begin{equation*}
	\mathrm{CBIE} + \alpha\cdot\mathrm{HBIE} = 0
\end{equation*}
with the usual choice of the coupling parameter $\alpha=\frac{\imag}{k}$~\cite{Zheng2015itb}. More precisely, the BM formulation is given by
\begin{equation}\label{Eq3:BM}
\begin{aligned}
	&C^\pm(\vec{x})p(\vec{x}) + \int_\Gamma \pderiv{\Phi_k(\vec{x},\vec{y})}{n(\vec{y})} p(\vec{y}) \idiff\Gamma(\vec{y}) + \alpha\int_\Gamma \frac{\partial^2\Phi_k(\vec{x},\vec{y})}{\partial\vec{n}(\vec{y})\partial\vec{n}(\vec{x})} p(\vec{y}) \idiff\Gamma(\vec{y})\\
	\quad&=\int_\Gamma \Phi_k(\vec{x},\vec{y}) \pderiv{p(\vec{y})}{n(\vec{y})} \idiff\Gamma(\vec{y})+\alpha\int_\Gamma \pderiv{\Phi_k(\vec{x},\vec{y})}{n(\vec{x})} \pderiv{p(\vec{y})}{n(\vec{y})} \idiff\Gamma(\vec{y}) - \alpha C^\pm(\vec{x})\pderiv{p(\vec{x})}{n(\vec{x})}.
\end{aligned}
\end{equation}
As in \cite{Simpson2014aib}, we restrict our analysis to direct IGABEM formulations (indirect IGABEM formulations are considered in~\cite{Coox2017aii,Dolz2018afi,Wu2020iib}).

\subsection{Regularization techniques}
Using~\Cref{Eq3:jumpTerm} the CBIE can be regularized as follows
\begin{equation}\label{Eq3:regulCBIE}
	-\frac{1}{2}p(\vec{x})(1\pm 1) + \int_\Gamma \pderiv{\Phi_k(\vec{x},\vec{y})}{n(\vec{y})}p(\vec{y}) - \pderiv{\Phi_0(\vec{x},\vec{y})}{n(\vec{y})} p(\vec{x}) \idiff\Gamma(\vec{y}) = \int_\Gamma \Phi_k(\vec{x},\vec{y}) \pderiv{p(\vec{y})}{n(\vec{y})} \idiff\Gamma(\vec{y}).
\end{equation}
With the identities~\cite{Liu1999anf,Taus2016iao}
\begin{align}
	\int_\Gamma\frac{\partial^2\Phi_0(\vec{x},\vec{y})}{\partial\vec{n}(\vec{y})\partial\vec{n}(\vec{x})} (\vec{y}-\vec{x})\idiff\Gamma(\vec{y}) &= \int_\Gamma \pderiv{\Phi_0(\vec{x},\vec{y})}{n(\vec{x})} \vec{n}(\vec{y})
	+ \pderiv{\Phi_0(\vec{x},\vec{y})}{n(\vec{y})}\vec{n}(\vec{x}) \idiff\Gamma(\vec{y})\label{Eq3:secondIdentity}\\
	\int_\Gamma\frac{\partial^2\Phi_0(\vec{x},\vec{y})}{\partial\vec{n}(\vec{y})\partial\vec{n}(\vec{x})}\idiff\Gamma(\vec{y}) &= 0\nonumber
\end{align}
the regularization of the HBIE is given by
\begin{align}\label{Eq3:regulHBIE}
\begin{split}
	&\int_\Gamma\left[\frac{\partial^2\Phi_k(\vec{x},\vec{y})}{\partial\vec{n}(\vec{y})\partial\vec{n}(\vec{x})} - \frac{\partial^2\Phi_0(\vec{x},\vec{y})}{\partial\vec{n}(\vec{y})\partial\vec{n}(\vec{x})}\right] p(\vec{y}) \idiff\Gamma(\vec{y}) \\
	&+\int_\Gamma \frac{\partial^2\Phi_0(\vec{x},\vec{y})}{\partial\vec{n}(\vec{y})\partial\vec{n}(\vec{x})}\left[ p(\vec{y}) - p(\vec{x}) - \pderiv{p(\vec{x})}{v_j} \vec{v}_j\cdot(\vec{y}-\vec{x})\right] \idiff\Gamma(\vec{y})\\
	&+\pderiv{p(\vec{x})}{v_j}\int_\Gamma \pderiv{\Phi_0(\vec{x},\vec{y})}{n(\vec{x})}\vec{v}_j\cdot \vec{n}(\vec{y})
	+ \pderiv{\Phi_0(\vec{x},\vec{y})}{n(\vec{y})}\vec{v}_j\cdot \vec{n}(\vec{x}) \idiff\Gamma(\vec{y}) \\
	&= \int_\Gamma \left[\pderiv{\Phi_k(\vec{x},\vec{y})}{n(\vec{x})}+\pderiv{\Phi_0(\vec{x},\vec{y})}{n(\vec{y})}\right] \pderiv{p(\vec{y})}{n(\vec{y})} \idiff\Gamma(\vec{y}) \\
	&{\hskip1em\relax}-\int_\Gamma \pderiv{\Phi_0(\vec{x},\vec{y})}{n(\vec{y})}\left[ \pderiv{p(\vec{y})}{n(\vec{y})} - \pderiv{p(\vec{x})}{n(\vec{x})}\right]\idiff\Gamma(\vec{y}) + \frac{1}{2}\pderiv{p(\vec{x})}{n(\vec{x})}(1\pm 1)\\
	&{\hskip1em\relax}-\pderiv{p(\vec{x})}{n(\vec{x})}\left[\int_\Gamma \pderiv{\Phi_0(\vec{x},\vec{y})}{n(\vec{x})}\vec{n}(\vec{x})\cdot \vec{n}(\vec{y}) + \pderiv{\Phi_0(\vec{x},\vec{y})}{n(\vec{y})} \idiff\Gamma(\vec{y}) - \int_\Gamma \frac{\partial^2\Phi_0(\vec{x},\vec{y})}{\partial\vec{n}(\vec{y})\partial\vec{n}(\vec{x})}\vec{n}(\vec{x})\cdot(\vec{y}-\vec{x}) \idiff\Gamma(\vec{y})\right]
	\end{split}
\end{align}
where the summation over the indices $j=1,2$ is implied, and $\vec{v}_j$ (cf. \cite[Fig. 2]{Simpson2014aib}) is an orthonormal set of (unit) vectors at $\vec{x}$ such that $\vec{v}
_3 = \vec{n}$, $\vec{v}_1=\vec{e}_\upxi$ and $\vec{v}_2 = \vec{v}_3\times\vec{v}_1$ with the following notation
\begin{equation*}
	\vec{e}_\upxi = \frac{1}{h_\upxi}\pderiv{\vec{x}}{\xi},\quad \vec{e}_\upeta = \frac{1}{h_\upeta}\pderiv{\vec{x}}{\eta},\quad  h_\upxi = \left|\pderiv{\vec{x}}{\xi}\right|,\quad h_\upeta = \left|\pderiv{\vec{x}}{\eta}\right|.
\end{equation*}
Here, $\xi$ and $\eta$ are the parameters for the surface parametrization. Note that~\cite[p. 219]{Scott2013ibe}
\begin{align*}
\pderiv{p(\vec{x})}{v_1} &= \frac{1}{h_\upxi}\pderiv{p(\vec{x})}{\xi}\\
\pderiv{p(\vec{x})}{v_2} &= -\frac{1}{h_\upxi}\frac{\cos\theta}{\sin\theta}\pderiv{p(\vec{x})}{\xi} + \frac{1}{h_\upeta}\frac{1}{\sin\theta}\pderiv{p(\vec{x})}{\eta}\\
\end{align*}
where $\theta$ is the angle between $\vec{e}_\upxi$ and $\vec{e}_\upeta$. The integrals in \Cref{Eq3:regulCBIE,Eq3:regulHBIE} are at most weakly singular.

In practice~\cite{Scott2013ibe}, the integrals in the BIEs are discretized individually using the same quadrature points making several terms cancel.

Another approach for regularizing the CBIE in \Cref{Eq3:CBIE} is presented in \cite{Sun2015bri}. Consider the function
\begin{equation*}
	\Psi(\vec{y}) = p(\vec{x})\Psi_1(\vec{y}) + \pderiv{p}{n}\Big\vert_{\vec{y}=\vec{x}}\Psi_2(\vec{y})
\end{equation*}
where $\Psi_1(\vec{y})$ and $\Psi_2(\vec{y})$ solve
\begin{equation*}
	\nabla^2\Psi_1(\vec{y}) + k^2\Psi_1(\vec{y}) = 0,\quad \Psi_1(\vec{x})= 1\quad \nabla\Psi_1(\vec{x})\cdot \vec{n}(\vec{x}) = 0
\end{equation*}
and
\begin{equation*}
	\nabla^2\Psi_2(\vec{y}) + k^2\Psi_2(\vec{y}) = 0,\quad \Psi_2(\vec{x})= 0\quad \nabla\Psi_2(\vec{x})\cdot \vec{n}(\vec{x}) = 1.
\end{equation*}
The idea is that $\Psi(\vec{y})$ also solves BIEs such that a subtraction of two such BIEs yields regularization of the integrand. There exist a lot of freedom in choosing functions $\Psi_1$ and $\Psi_2$ that satisfy these constraints. The original ones suggested by \cite{Sun2015bri} are given by
\begin{equation}\label{Eq3:psi_12_1}
	\Psi_1^{(1)}(\vec{y}) = \frac{C_1\cos[k(R_1-C_1)]}{R_1} + \frac{\sin[k(R_1-C_1)]}{kR_1}\quad\text{and}\quad \Psi_2^{(1)}(\vec{y}) = \frac{C_1^2\sin[k(R_1-C_1)]}{C_2kR_1}
\end{equation}
where
\begin{equation*}
	R_1(\vec{y}) = |\vec{y}-\vec{x}_1|, \quad C_1 = |\vec{x}-\vec{x}_1|,\quad C_2 = (\vec{x}-\vec{x}_1)\cdot\vec{n}(\vec{x}).
\end{equation*}
The point $\vec{x}_1$ must lie outside the solution domain and chosen such that $C_2\neq 0$ (for the sphere and the torus geometry in this work, we use $\vec{x}_1=\vec{x} - \vec{n}(\vec{x})$). However, these functions do not satisfy an exterior problem (as they do not satisfy the Sommerfeld radiation condition). This problem is resolved by adding a non-vanishing integral at infinity as described in \cite{Sun2015bri}.

One can easily create functions that also satisfy the Sommerfeld radiation condition, simply by basing the functions on the fundamental solutions in \Cref{Eq3:FreeSpaceGrensFunction}
\begin{equation}\label{Eq3:psi_12_2}
	\Psi_1^{(2)}(\vec{y}) = \frac{1}{C_1}\frac{\Phi_k(\vec{x}_1,\vec{y})}{\Phi_k(\vec{x}_1,\vec{x})}+\left(1-\frac{1}{C_1}\right)\frac{\Phi_k(\vec{x}_2,\vec{y})}{\Phi_k(\vec{x}_2,\vec{x})},\quad\text{and}\quad \Psi_2^{(2)}(\vec{y}) = \frac{1}{C_2}\left[\frac{\Phi_k(\vec{x}_1,\vec{y})}{\Phi_k(\vec{x}_1,\vec{x})}-\frac{\Phi_k(\vec{x}_2,\vec{y})}{\Phi_k(\vec{x}_2,\vec{x})}\right]
\end{equation}
where
\begin{equation*}
	C_1 = 1- \frac{r_2^2(\imag k r_1-1)\left(\vec{x}_1-\vec{x}\right)\cdot \vec{n}(\vec{x})}{r_1^2(\imag k r_2-1)\left(\vec{x}_2-\vec{x}\right)\cdot \vec{n}(\vec{x})},\quad C_2 = \frac{C_1}{r_2^2} (\imag k r_2-1)\left(\vec{x}_2-\vec{x}\right)\cdot \vec{n}(\vec{x}),\quad r_1 = |\vec{x}_1-\vec{x}|,\quad r_2 = |\vec{x}_2-\vec{x}|
\end{equation*}
The points $\vec{x}_1$ and $\vec{x}_2$ must lie outside the solution domain and chosen such that $C_1\neq 0$ and $C_2\neq 0$ (for the sphere and the torus geometry in this work, we use $\vec{x}_1=\vec{x} - \frac{1}{2}\vec{n}(\vec{x})$ and $\vec{x}_2=\vec{x} - \vec{n}(\vec{x})$, respectively).

Alternatively, for the interior problem one could choose
\begin{equation*}
	\Psi_1^{(3)}(\vec{y}) = \frac{\vec{k}_2\cdot\vec{n}(\vec{x})\euler^{\imag\vec{k}_1\cdot(\vec{y}-\vec{x})}-\vec{k}_1\cdot\vec{n}(\vec{x})\euler^{\imag\vec{k}_2\cdot(\vec{y}-\vec{x})}}{(\vec{k}_2-\vec{k}_1)\cdot\vec{n}(\vec{x})}\quad\text{and}\quad \Psi_2^{(3)}(\vec{y}) = \frac{\euler^{\imag\vec{k}_2\cdot(\vec{y}-\vec{x})}-\euler^{\imag\vec{k}_1\cdot(\vec{y}-\vec{x})}}{\imag(\vec{k}_2-\vec{k}_1)\cdot\vec{n}(\vec{x})}
\end{equation*}
where $\vec{k}_2 = k\vec{d}_2$ and $\vec{k}_1 = k\vec{d}_1$ are the wave vectors for the plane wave in the direction of the unit vectors $\vec{d}_1$ and $\vec{d}_2$, respectively. Choosing $\vec{d}_2=\vec{d}_1+\vec{n}(\vec{x})$ we get (with $|\vec{n}(\vec{x})|=1$)
\begin{equation*}
	\Psi_1^{(3)}(\vec{y}) = \left(\vec{d}_1\cdot\vec{n}(\vec{x})+1\right)\euler^{\imag k\vec{d}_1(\vec{y}-\vec{x})} - \vec{d}_1\cdot\vec{n}(\vec{x})\euler^{\imag k(\vec{d}_1 + \vec{n}(\vec{x}))(\vec{y}-\vec{x})}
\end{equation*}
and
\begin{equation*}
	\Psi_2^{(3)}(\vec{y}) = \frac{\imag}{k}\left(\euler^{\imag k\vec{d}_1\cdot(\vec{y}-\vec{x})}-\euler^{\imag k(\vec{d}_1+\vec{n}(\vec{x}))\cdot(\vec{y}-\vec{x})}\right)
\end{equation*}
where
\begin{equation}\label{Eq3:generald1}
	\vec{d}_1 = \begin{cases} \frac{\sqrt{3}}{2\sqrt{1-n_1(\vec{x})^2}}\begin{bmatrix}
		(1-n_1(\vec{x})^2)\cos\theta_1\\
		-n_1(\vec{x}) n_2(\vec{x}) \cos\theta_1 + n_3(\vec{x}) \sin\theta_1\\
		-n_1(\vec{x}) n_3(\vec{x}) \cos\theta_1 - n_2(\vec{x}) \sin\theta_1
	\end{bmatrix}-\frac{1}{2}\vec{n}(\vec{x}) & |n_1(\vec{x})| < \frac{1}{\sqrt{2}}\\
	\frac{\sqrt{3}}{2\sqrt{1-n_2(\vec{x})^2}}\begin{bmatrix}
		-n_1(\vec{x}) n_2(\vec{x}) \sin\theta_2 - n_3(\vec{x}) \cos\theta_2\\
		(1-n_2(\vec{x})^2)\sin\theta_2\\
		-n_2(\vec{x}) n_3(\vec{x}) \sin\theta_2 + n_1(\vec{x}) \cos\theta_2
	\end{bmatrix}-\frac{1}{2}\vec{n}(\vec{x}) & \text{otherwise,}
	\end{cases}
\end{equation}
for some free parameters $\theta_1$ and $\theta_2$. Choosing $\theta_1=-\PI/2$ and $\theta_2=-\PI$ yields
\begin{equation}\label{Eq3:d1d2}
	\vec{d}_1 = \begin{cases}\frac{\sqrt{3}}{2\sqrt{1-n_1(\vec{x})^2}}\vec{e}_1\times\vec{n}(\vec{x})-\frac{1}{2}\vec{n}(\vec{x}) & |n_1(\vec{x})| < \frac{1}{\sqrt{2}}\\
	\frac{\sqrt{3}}{2\sqrt{1-n_2(\vec{x})^2}}\vec{e}_2\times\vec{n}(\vec{x})-\frac{1}{2}\vec{n}(\vec{x}) & \text{otherwise.}
	\end{cases}
\end{equation}
Then, $\vec{d}_1\cdot\vec{n}(\vec{x}) = -\frac{1}{2}$ and
\begin{equation}
	\Psi_1^{(3)}(\vec{y}) = \frac12\left(\euler^{\imag k\vec{d}_1(\vec{y}-\vec{x})} +\euler^{\imag k\vec{d}_2(\vec{y}-\vec{x})}\right)\quad\text{and}\quad \Psi_2^{(3)}(\vec{y}) = \frac{\imag}{k}\left(\euler^{\imag k\vec{d}_1\cdot(\vec{y}-\vec{x})}-\euler^{\imag k\vec{d}_2\cdot(\vec{y}-\vec{x})}\right).
\end{equation}
The advantage of this choice over the former two choices is that it does not require finding points ($\vec{x}_1$ and $\vec{x}_2$) outside the solution domain that satisfy a given criterion.

If
\begin{equation*}
	\Psi(\vec{y}) = p(\vec{x})\Psi_1^{(1)}(\vec{y}) + \pderiv{p}{n}\Big\vert_{\vec{y}=\vec{x}}\Psi_2^{(1)}(\vec{y})
\end{equation*}
then\footnote{Recall that the upper plus sign in $\pm$ (and negative sign for $\mp$) is chosen for the exterior problem while the negative sign in $\pm$ (and positive sign for $\mp$) is chosen for the interior problem.} (cf.~\cite{Sun2015bri}) 
\begin{equation}\label{Eq3:RCBIE1}
\begin{aligned}
	&\frac{1}{2}p(\vec{x})\left[1\mp 1 - \left(1+\frac{\imag}{kC_1}\right)\left(1-\euler^{2\imag kC_1}\right)\right] + \int_\Gamma \left(p(\vec{y}) - p(\vec{x})\Psi_1(\vec{y}) - \pderiv{p}{n}\Big\vert_{\vec{y}=\vec{x}}\Psi_2(\vec{y})\right)\pderiv{\Phi_k(\vec{x},\vec{y})}{n(\vec{y})}  \idiff\Gamma(\vec{y}) \\
	&{\hskip5em\relax}= \frac{\imag C_1}{2kC_2} \left(1-\euler^{2\imag kC_1}\right)\pderiv{p}{n}\Big\vert_{\vec{y}=\vec{x}} +  \int_\Gamma \left(\pderiv{p(\vec{y})}{n(\vec{y})} - p(\vec{x})\pderiv{\Psi_1(\vec{y})}{n(\vec{y})} - \pderiv{p}{n}\Big\vert_{\vec{y}=\vec{x}}\pderiv{\Psi_2(\vec{y})}{n(\vec{y})}\right)\Phi_k(\vec{x},\vec{y}) \idiff\Gamma(\vec{y}).
\end{aligned}
\end{equation}
We refer to this integral equation as the first regularized CBIE (RCBIE1). If
\begin{equation*}
	\Psi(\vec{y}) = p(\vec{x})\Psi_1^{(2)}(\vec{y}) + \pderiv{p}{n}\Big\vert_{\vec{y}=\vec{x}}\Psi_2^{(2)}(\vec{y})
\end{equation*}
then $\Psi(\vec{y})$ solves the exterior problem of \Cref{Eq3:CBIE} such that
\begin{equation}\label{Eq3:RCBIE2}
\begin{aligned}
	&\frac{1}{2}p(\vec{x})(1\mp 1) + \int_\Gamma \left(p(\vec{y}) - p(\vec{x})\Psi_1(\vec{y}) - \pderiv{p}{n}\Big\vert_{\vec{y}=\vec{x}}\Psi_2(\vec{y})\right)\pderiv{\Phi_k(\vec{x},\vec{y})}{n(\vec{y})}  \idiff\Gamma(\vec{y}) \\
	&= \int_\Gamma \left(\pderiv{p(\vec{y})}{n(\vec{y})} - p(\vec{x})\pderiv{\Psi_1(\vec{y})}{n(\vec{y})} - \pderiv{p}{n}\Big\vert_{\vec{y}=\vec{x}}\pderiv{\Psi_2(\vec{y})}{n(\vec{y})}\right)\Phi_k(\vec{x},\vec{y}) \idiff\Gamma(\vec{y}).
\end{aligned}
\end{equation}
We refer to this integral equation as the second regularized CBIE (RCBIE2). If
\begin{equation*}
	\Psi(\vec{y}) = p(\vec{x})\Psi_1^{(3)}(\vec{y}) + \pderiv{p}{n}\Big\vert_{\vec{y}=\vec{x}}\Psi_2^{(3)}(\vec{y})
\end{equation*}
then $\Psi(\vec{y})$ solves the interior problem of \Cref{Eq3:CBIE} such that
\begin{equation}\label{Eq3:RCBIE3}
\begin{aligned}
	&-\frac{1}{2}p(\vec{x})(1\pm 1) + \int_\Gamma \left(p(\vec{y}) - p(\vec{x})\Psi_1(\vec{y}) - \pderiv{p}{n}\Big\vert_{\vec{y}=\vec{x}}\Psi_2(\vec{y})\right)\pderiv{\Phi_k(\vec{x},\vec{y})}{n(\vec{y})}  \idiff\Gamma(\vec{y}) \\
	&{\hskip4em\relax}= \int_\Gamma \left(\pderiv{p(\vec{y})}{n(\vec{y})} - p(\vec{x})\pderiv{\Psi_1(\vec{y})}{n(\vec{y})} - \pderiv{p}{n}\Big\vert_{\vec{y}=\vec{x}}\pderiv{\Psi_2(\vec{y})}{n(\vec{y})}\right)\Phi_k(\vec{x},\vec{y}) \idiff\Gamma(\vec{y}).
\end{aligned}
\end{equation}
We refer to this integral equation as the third regularized CBIE (RCBIE3). These integrals have bounded integrands~\cite{Klaseboer2012nsb} and are thus a further regularization of \Cref{Eq3:regulCBIE}.

\subsection{Rigid scattering problems}
For rigid (exterior) scattering problems the boundary integral equations are simplified somewhat. Consider an incident plane wave
\begin{equation*}
	p_{\mathrm{inc}}(\vec{x}) = P_{\mathrm{inc}}\euler^{\imag \vec{k}\cdot\vec{x}}
\end{equation*}
scattered by the boundary $\Gamma$. Here, $P_{\mathrm{inc}}$ is the amplitude, and $\vec{k}$ is the wave vector. Combining \Cref{Thm:InteriorProblem} and \Cref{Thm:ExteriorProblem} we can write
\begin{equation*}
	p_{\mathrm{tot}}(\vec{x}) = p_{\mathrm{inc}}(\vec{x}) +\mathcal{D}_k\gamma^+ p_{\mathrm{tot}}(\vec{x}) - \mathcal{S}_k\partial_n^+ p_{\mathrm{tot}}(\vec{x})
\end{equation*}
where $p_{\mathrm{tot}}=p+p_{\mathrm{inc}}$ is the total field and $p$ is the scattered field satisfying the assumptions of \Cref{Thm:ExteriorProblem}. 

For rigid scattering we have $\partial_n^+ p_{\mathrm{tot}}(\vec{x}) = 0$, such that the regularized CBIE in \Cref{Eq3:regulCBIE} and HBIE in \Cref{Eq3:regulHBIE} reduce to\footnote{Note that this CBIE formulation no longer contains weakly singular integrals (only integrals with bounded integrands).}
\begin{equation}\label{Eq3:regulCBIErigid}
	-p_{\mathrm{tot}}(\vec{x}) + \int_\Gamma \pderiv{\Phi_k(\vec{x},\vec{y})}{n(\vec{y})}p_{\mathrm{tot}}(\vec{y}) - \pderiv{\Phi_0(\vec{x},\vec{y})}{n(\vec{y})} p_{\mathrm{tot}}(\vec{x}) \idiff\Gamma(\vec{y}) = -p_{\mathrm{inc}}(\vec{x})
\end{equation}
and 
\begin{align*}
	&\int_\Gamma\left[\frac{\partial^2\Phi_k(\vec{x},\vec{y})}{\partial\vec{n}(\vec{y})\partial\vec{n}(\vec{x})} - \frac{\partial^2\Phi_0(\vec{x},\vec{y})}{\partial\vec{n}(\vec{y})\partial\vec{n}(\vec{x})}\right] p_{\mathrm{tot}}(\vec{y}) \idiff\Gamma(\vec{y}) \\
	&\qquad+\int_\Gamma \frac{\partial^2\Phi_0(\vec{x},\vec{y})}{\partial\vec{n}(\vec{y})\partial\vec{n}(\vec{x})}\left[ p_{\mathrm{tot}}(\vec{y}) - p_{\mathrm{tot}}(\vec{x}) - \pderiv{p_{\mathrm{tot}}(\vec{x})}{v_j} \vec{v}_j\cdot(\vec{y}-\vec{x})\right] \idiff\Gamma(\vec{y})\\
	&\qquad+\pderiv{p_{\mathrm{tot}}(\vec{x})}{v_j}\int_\Gamma \pderiv{\Phi_0(\vec{x},\vec{y})}{n(\vec{x})}\vec{v}_j\cdot \vec{n}(\vec{y})
	+ \pderiv{\Phi_0(\vec{x},\vec{y})}{n(\vec{y})}\vec{v}_j\cdot \vec{n}(\vec{x}) \idiff\Gamma(\vec{y}) = -\pderiv{p_{\mathrm{inc}}(\vec{x})}{n(\vec{x})},
\end{align*}
respectively.
%
%
In a similar fashion \Cref{Eq3:RCBIE1}, \Cref{Eq3:RCBIE2} and \Cref{Eq3:RCBIE3} can be reformulated as 
\begin{align*}
	&-\frac{1}{2}p_{\mathrm{tot}}(\vec{x})\left(1+\frac{\imag}{kC_1}\right)\left(1-\euler^{2\imag kC_1}\right) \\
	&\qquad+\int_\Gamma \left(p_{\mathrm{tot}}(\vec{y}) - p_{\mathrm{tot}}(\vec{x})\Psi_1(\vec{y})\right)\pderiv{\Phi_k(\vec{x},\vec{y})}{n(\vec{y})}+ p_{\mathrm{tot}}(\vec{x})\pderiv{\Psi_1(\vec{y})}{n(\vec{y})} \Phi_k(\vec{x},\vec{y})\idiff\Gamma(\vec{y}) =-p_{\mathrm{inc}}(\vec{x}),
\end{align*}
\begin{align*}
	&\int_\Gamma \left(p_{\mathrm{tot}}(\vec{y}) - p_{\mathrm{tot}}(\vec{x})\Psi_1(\vec{y})\right)\pderiv{\Phi_k(\vec{x},\vec{y})}{n(\vec{y})}+ p_{\mathrm{tot}}(\vec{x})\pderiv{\Psi_1(\vec{y})}{n(\vec{y})} \Phi_k(\vec{x},\vec{y})\idiff\Gamma(\vec{y}) =-p_{\mathrm{inc}}(\vec{x})
\end{align*}
and
\begin{align*}
	&-p_{\mathrm{tot}}(\vec{x})+\int_\Gamma \left(p_{\mathrm{tot}}(\vec{y}) - p_{\mathrm{tot}}(\vec{x})\Psi_1(\vec{y})\right)\pderiv{\Phi_k(\vec{x},\vec{y})}{n(\vec{y})} + p_{\mathrm{tot}}(\vec{x})\pderiv{\Psi_1(\vec{y})}{n(\vec{y})} \Phi_k(\vec{x},\vec{y})\idiff\Gamma(\vec{y}) =-p_{\mathrm{inc}}(\vec{x}),
\end{align*}
respectively.

\section{Collocation and Galerkin formulations}
\label{Sec3:BEM}
For the discretization procedure we consider a finite dimensional trial space $V_h\subset V=H^{1/2}(\Gamma)$ which is built up by the same NURBS basis functions used to represent the CAD geometry. In this work, the geometry is assumed to be constructed by tensorial NURBS patches such that the geometry for each patch can be written as
\begin{equation*}
	\vec{X}(\xi,\eta) = \sum_{i=1}^n\sum_{j=1}^m R_{i,j}^{\check{p},\check{q}}(\xi,\eta)
\end{equation*}
with notation taken from and explained in~\cite[p. 51]{Cottrell2006iao}. For convenience we simplify the notation $R_{i,j}^{\check{p},\check{q}}$ to $R_{\tilde{i}}$ where the index $\tilde{i}$ represents a map from local indices to global indices (over all patches).

For the collocation formulations, we evaluate the BIEs at $n_{\mathrm{dofs}}$ collocation points, $\vec{x}_i\in\Gamma$. This forms an algebraic system of equations which can be solved to obtain the numerical solution. Throughout this work, the collocation points are chosen to be the Greville abscissae as described in \cite{Scott2013ibe}. 

The Galerkin formulations are obtained by multiplying the BIEs with a test function $q_{\mathrm{tot}}(\vec{x})$ and integrating over $\Gamma$. For brevity we only here consider rigid scattering problems with the CBIE formulation in \Cref{Eq3:regulCBIErigid}
\begin{align*}
	&-\int_\Gamma p_{\mathrm{tot}}(\vec{x})q_{\mathrm{tot}}(\vec{x})\idiff\Gamma(\vec{x}) + \int_\Gamma q_{\mathrm{tot}}(\vec{x})\int_\Gamma \pderiv{\Phi_k(\vec{x},\vec{y})}{n(\vec{y})}p_{\mathrm{tot}}(\vec{y}) - \pderiv{\Phi_0(\vec{x},\vec{y})}{n(\vec{y})} p_{\mathrm{tot}}(\vec{x}) \idiff\Gamma(\vec{y})\idiff\Gamma(\vec{x}) \\
	&{\hskip30em\relax}= -\int_\Gamma p_{\mathrm{inc}}(\vec{x})q_{\mathrm{tot}}(\vec{x})\idiff\Gamma(\vec{x}).
\end{align*}
Letting
\begin{equation*}
	p_{\mathrm{tot}}(\vec{x}) = \sum_{j=1}^{n_{\mathrm{dofs}}} u_j R_j(\vec{x}),
\end{equation*}
we get (by choosing $q_{\mathrm{tot}}(\vec{x})=R_i(\vec{x})$)
\begin{align*}
	&-\sum_{j=1}^{n_{\mathrm{dofs}}}u_j \left[\int_\Gamma R_j(\vec{x})R_i(\vec{x})\idiff\Gamma(\vec{x})+ \int_\Gamma R_i(\vec{x})\int_\Gamma \pderiv{\Phi_k(\vec{x},\vec{y})}{n(\vec{y})}R_j(\vec{y}) \idiff\Gamma(\vec{y})\idiff\Gamma(\vec{x})\right.\\ 
	&\left.{\hskip5em\relax}-\int_\Gamma R_i(\vec{x}) R_j(\vec{x})\int_\Gamma \pderiv{\Phi_0(\vec{x},\vec{y})}{n(\vec{y})} \idiff\Gamma(\vec{y})\idiff\Gamma(\vec{x})\right]= -\int_\Gamma p_{\mathrm{inc}}(\vec{x})R_i(\vec{x})\idiff\Gamma(\vec{x}),\quad \forall i=1,\dots, n_{\mathrm{dofs}},
\end{align*}
which results in a linear system of $n_{\mathrm{dofs}}$ equations. Instead of looping through all basis functions $R_i(\vec{x})$, it is advantageous to loop through the elements as done in finite element methods~\cite{Simpson2014aib}.

For the collocation formulations we prepend a letter ``C'' (i.e. CCBIE, CBM, CRCBIE1, etc.) and for the Galerkin formulations we prepend a letter ``G'' (i.e. GCBIE, GBM, GRCBIE1, etc.).

\section{Numerical evaluation of the boundary integral equations}
\label{Sec3:numericalQuad}
In~\cite[p. 286]{Simpson2014aib} an adaptive integration technique is used around the collocation points in order to resolve the singular behavior of the integrand. Every element not containing the source point is divided into\footnote{The function $\round{\cdot}$ is the rounding function, i.e. $\round{x} = \left\lfloor x+\frac{1}{2}\right\rfloor$, where $\lfloor x\rfloor = \max\{n\in\Z\st n\leq x\}$.}
\begin{equation}\label{Eq3:numberOfSubElements}
	n_{\mathrm{div}} = \left(1+\round{\frac{s_1 h}{l}}\right)^{d-1}
\end{equation}
sub elements at which standard quadrature is applied. Here, $l$ is the distance from the center of the element to the source point\footnote{Arguably, a better choice for $l$ would be the minimal distance between the source point and any point in the element as outlined in~\cite{Taus2015iaf}. It is not clear to the authors if this is an optimization as it requires additional computational effort.}, $h$ is the element size (largest diagonal of the element) and $s_1$ is a user defined parameter controlling the adaptivity in terms of quadrature point density. For the element containing the source point, the element is divided into 2 to 4 (triangular) sub elements (depending on the locations of the source point; at a corner, on an edge, or within an element) as described in~\cite{Scott2013ibe}. A polar integration is then applied to each triangle such that the weakly singular integrands are regularized.

We use $\check{p}_\upxi+1+n_{\mathrm{eqp},1}$ quadrature points within each sub-element in the $\xi$-direction, and $\check{p}_\upeta+1+n_{\mathrm{eqp},1}$ in the $\eta$-direction. In the polar integration we use $\check{p}_{\mathrm{max}}+1+n_{\mathrm{eqp},2}$ in each parameter direction where $\check{p}_{\mathrm{max}} = \max\{\check{p}_\upxi,\check{p}_\upeta\}$ for the Simpson method.

\begin{figure}
	\centering
	\includegraphics{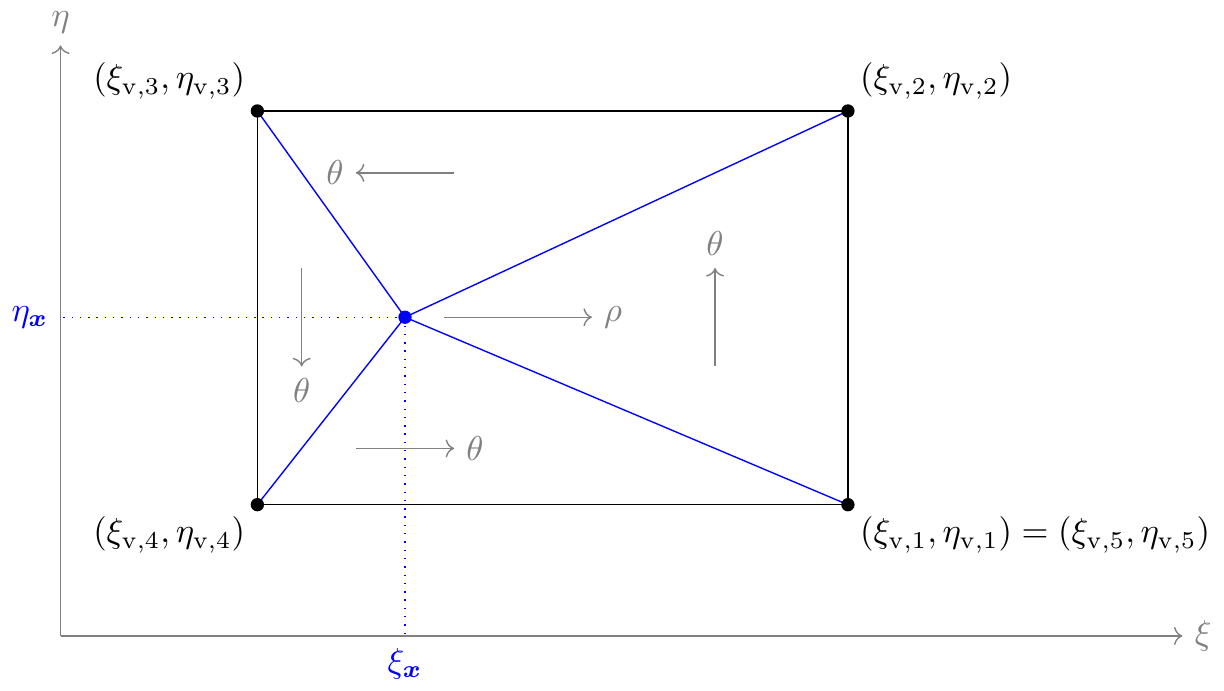}%
	\caption{\textbf{Numerical evaluation of the boundary integral equations}: The element containing the source point $\vec{x}$ is divided in (upto) 4 triangles in the parameter domain.}
	\label{Fig3:Duffy}
\end{figure}

In this work we present a modification to this routine inspired by Taus et al.~\cite{Taus2016iao,Taus2015iaf}. For each element not containing the source point, each (sub) element is divided into 4 until $s_1 h/l < 1$ where $h$ is the size of the (sub) element and $l$ is the distance from the (sub) element center to the source point. Whenever a (sub) element fulfills this requirement, standard quadrature is used with $\round{(\check{p}_\upxi+1)(s_1 h/l+1)}$ quadrature points in the $\xi$-direction and $\round{(\check{p}_\upeta+1)(s_1 h/l+1)}$ quadrature points in the $\eta$-direction. An alternative approach to the polar integration is here used. It is based on the transformation in~\cite{Duffy1982qoa} (for details see~\cite{Sauter1997qfh}), which avoids the problem of awkward integration limits opposite to the triangle vertex containing the singularity. Each triangular sub element is bilinearly transformed into the unit square. Consider the $i^{\mathrm{th}}$ triangular sub element with vertices $\{(\xi_{\vec{x}},\eta_{\vec{x}}),(\xi_{\mathrm{v},i},\eta_{\mathrm{v},i}),(\xi_{\mathrm{v},i+1},\eta_{\mathrm{v},i+1})\}$ in the parameter domain where $(\xi_{\vec{x}},\eta_{\vec{x}})$ is the parametric coordinate of $\vec{x}$ and $(\xi_{\mathrm{v},1},\eta_{\mathrm{v},1})=(\xi_{\mathrm{v},5},\eta_{\mathrm{v},5})$, $(\xi_{\mathrm{v},2},\eta_{\mathrm{v},2})$, $(\xi_{\mathrm{v},3},\eta_{\mathrm{v},3})$ and $(\xi_{\mathrm{v},4},\eta_{\mathrm{v},4})$ are the parametric coordinates for the four vertices of the element (see~\Cref{Fig3:Duffy}). The transformation is then given by ($\rho,\theta\in[0,1]$)
\begin{equation}
\begin{aligned}
	\xi &= \xi_{\vec{x}} + \rho(\xi_{\mathrm{v},i}-\xi_{\vec{x}}+(\xi_{\mathrm{v},i+1}-\xi_{\mathrm{v},i})\theta)\\
	\eta &= \eta_{\vec{x}} + \rho(\eta_{\mathrm{v},i}-\eta_{\vec{x}}+(\eta_{\mathrm{v},i+1}-\eta_{\mathrm{v},i})\theta)
\end{aligned}
\end{equation}
with Jacobian determinant given by
\begin{equation*}
	J_2 = \rho\left[(\xi_{\mathrm{v},i}-\xi_{\vec{x}} + (\xi_{\mathrm{v},i+1}-\xi_{\mathrm{v},i})\theta)(\eta_{\mathrm{v},i+1}-\eta_{\mathrm{v},i}) - (\eta_{\mathrm{v},i}-\eta_{\vec{x}} + (\eta_{\mathrm{v},i+1}-\eta_{\mathrm{v},i})\theta)(\xi_{\mathrm{v},i+1}-\xi_{\mathrm{v},i})\right].
\end{equation*}
The factor $\rho$ in the Jacobian determinant is responsible for regularizing the weakly singular integral. Note that $J_2=0$ for the collapsed triangle(s) when $\vec{x}$ lies on the edge (vertex) of the element. 
Each triangular sub element is divided into $n_{\mathrm{div},\uptheta}^{(i)}$ sub elements (in the $i^{\mathrm{th}}$ triangle) in the $\theta$-direction and $n_{\mathrm{div},\mathrm{r}}$ sub elements in the radial direction, where
\begin{equation*}
	n_{\mathrm{div},\uptheta}^{(i)} = \ceil{s_2\frac{\theta_{\mathrm{dir}}^{(i)}}{\ang{90}}},\quad n_{\mathrm{div},\mathrm{r}} = \ceil{s_2},\quad s_2 = \frac{\check{p}_{\mathrm{max}}+1+n_{\mathrm{eqp},2}}{2(\check{p}_{\mathrm{max}}+1)}.
\end{equation*}
Here, $\theta_{\mathrm{dir}}^{(i)}$ is the interior angle (in the parent domain) neighboring the source point of the initial sub triangle $i$. The reason for the sub division of the triangles (as opposed to use high order quadrature) is that a high number of quadrature points is here needed (which will later be illustrated). This sub division maps each sub element (in the $(\rho,\theta)$-domain) to the reference domain $[-1,1]^2$ by the linear transformation
\begin{equation}
\begin{aligned}
	\rho &= \rho_j + \frac12(\rho_{j+1}-\rho_j)(\tilde{\rho}+1),\qquad\rho_j = \frac{j}{n_{\mathrm{div},\mathrm{r}}},\qquad j=0,\dots,n_{\mathrm{div},\mathrm{r}}-1\\
	\theta &= \theta_l + \frac12(\theta_{l+1}-\theta_l)(\tilde{\theta}+1),\qquad\theta_l=\frac{l}{n_{\mathrm{div},\uptheta}^{(i)}},\qquad l=0,\dots,n_{\mathrm{div},\uptheta}^{(i)}-1
\end{aligned}
\end{equation}
with Jacobian determinant $J_3=1/(4n_{\mathrm{div},\uptheta}^{(i)}n_{\mathrm{div},\mathrm{r}})$. Each of these sub elements are now evaluated using $2(\check{p}_{\mathrm{max}}+1)$ quadrature points in both parametric directions.

For the Galerkin formulations the integral integrating the BIEs uses $(\check{p}_\upxi+1+n_{\mathrm{eqp},1})\times(\check{p}_\upeta+1+n_{\mathrm{eqp},1})$ quadrature points over each element. If not otherwise stated, we shall use $n_{\mathrm{eqp},1}=0$ throughout this work.

\begin{figure}
	\begin{subfigure}[t]{\textwidth}
		\centering
		\includegraphics[width=0.49\textwidth]{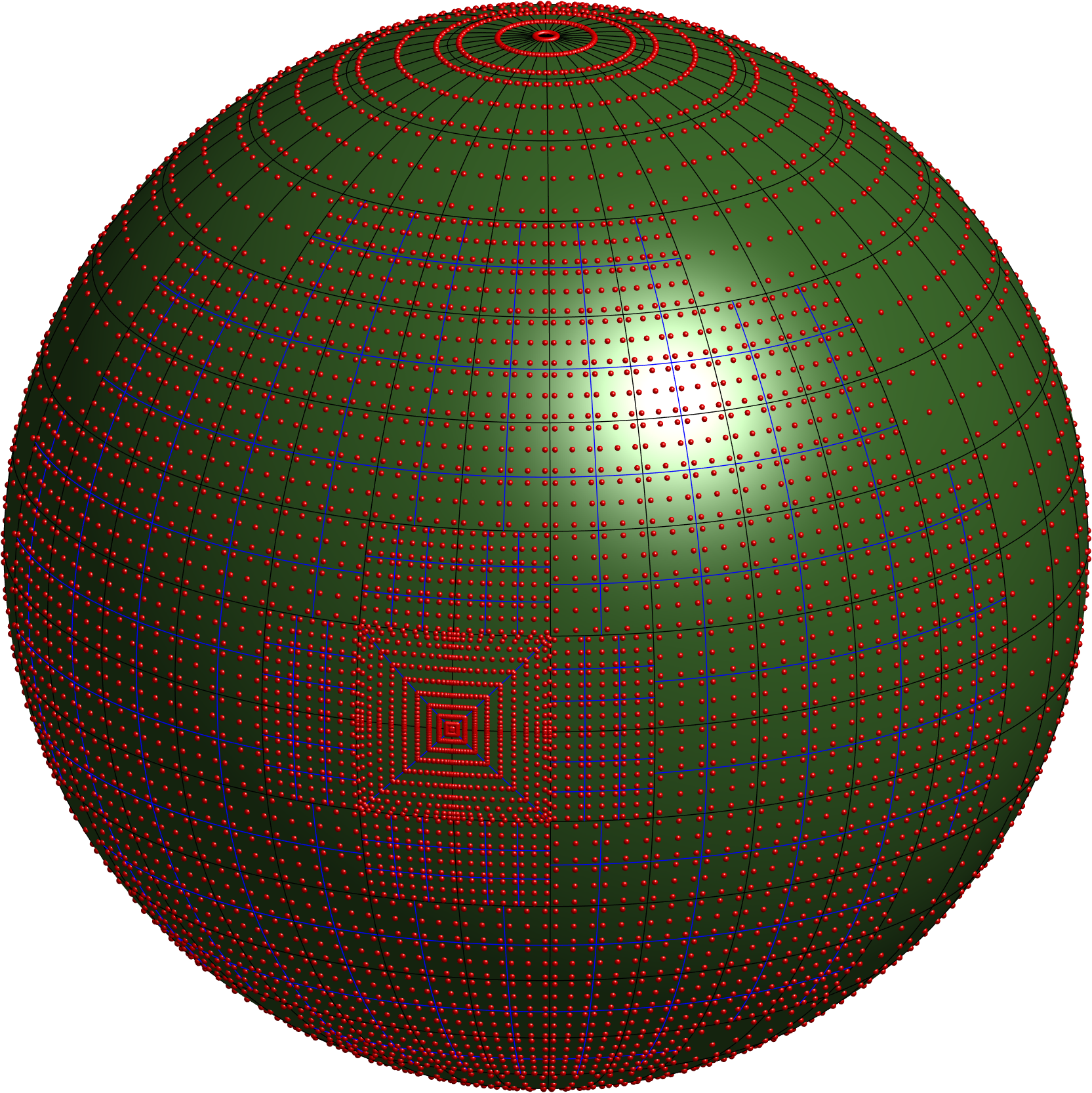}%
		\hspace*{0.02\textwidth}%
		\includegraphics[width=0.49\textwidth]{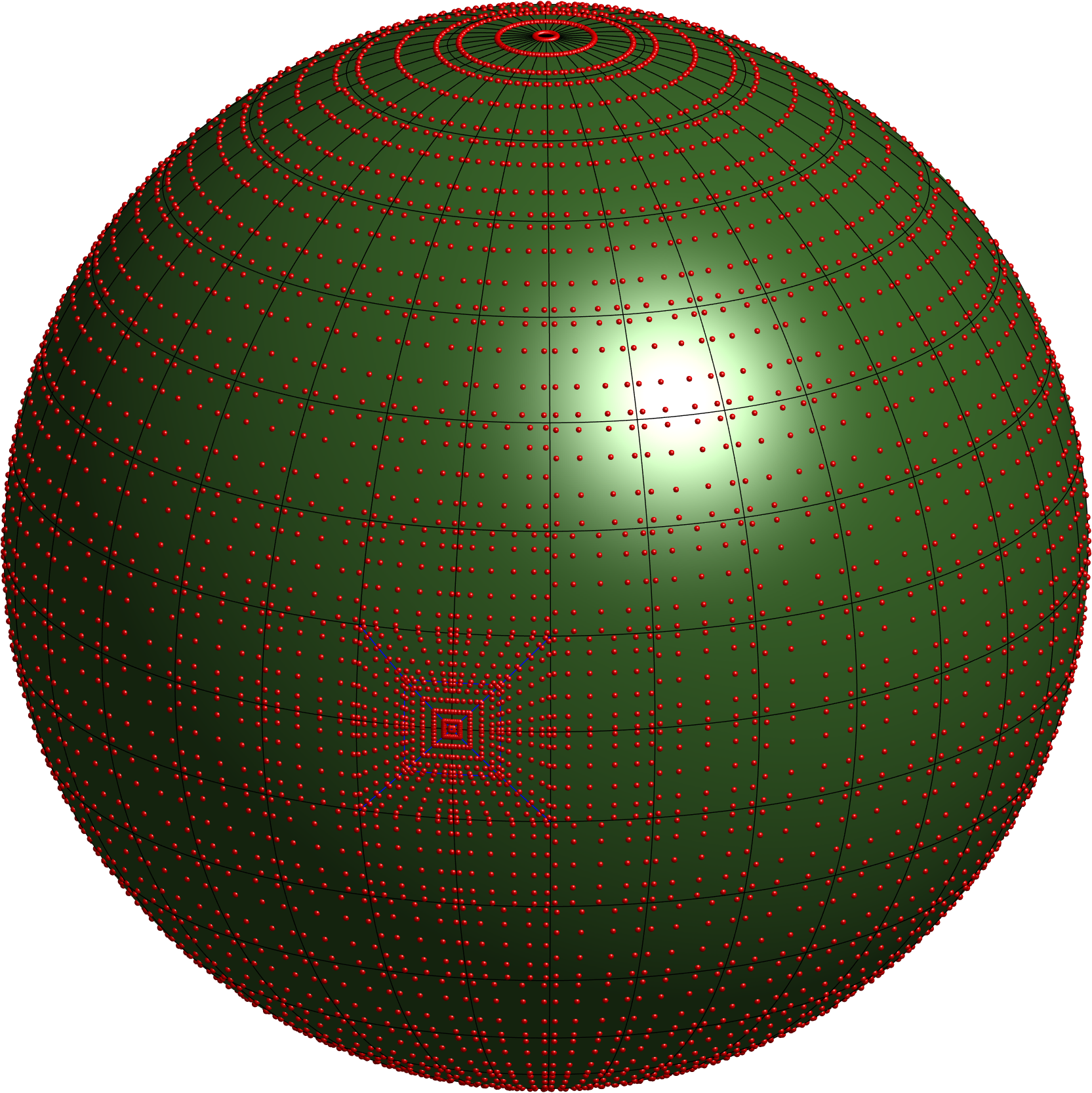}
		\caption{The $317^{\mathrm{th}}$ collocation point is at a vertex shared by four elements which are divided into two (triangular) sub elements. The total number of quadrature points for the left and right figures are 8492 and 6880, respectively.}
	\end{subfigure}
	\par\bigskip
	\par\bigskip
	\begin{subfigure}[t]{\textwidth}
		\centering
		\includegraphics[width=0.49\textwidth]{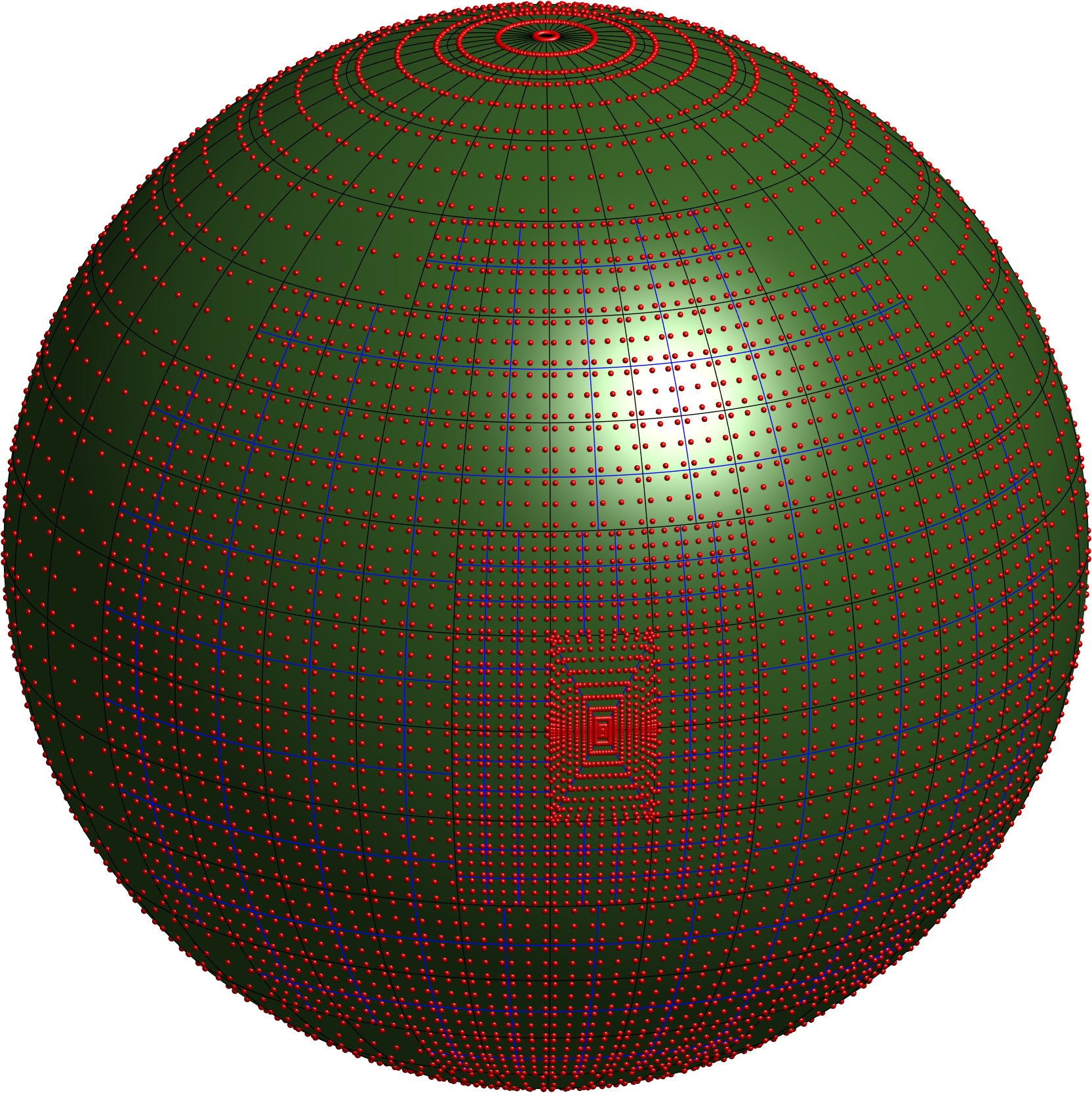}%
		\hspace*{0.02\textwidth}%
		\includegraphics[width=0.49\textwidth]{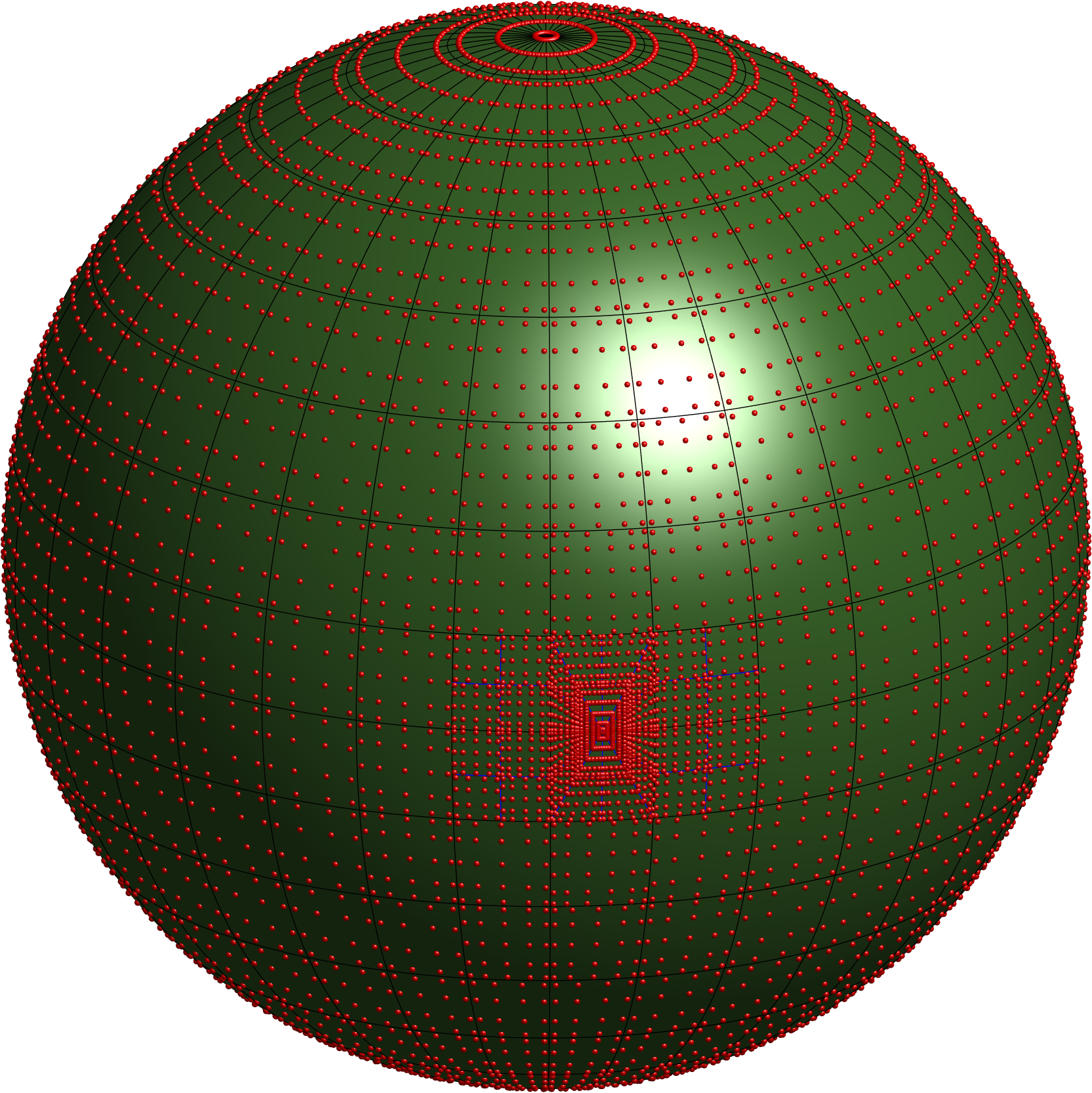}
		\caption{The $319^{\mathrm{th}}$ collocation point is on an element edge shared by two elements which are divided into three (triangular) sub elements. The total number of quadrature points for the left and right figures are 8250 and 7456, respectively.}
	\end{subfigure}
	\caption{\textbf{Numerical evaluation of the boundary integral equations}: The figures to the left are the integration procedure in~\cite{Simpson2014aib} (with $s_1=2$). The sub-element divisions are here shown by blue lines (the black lines are the element edges). The red points are the quadrature points. Here, $n_{\mathrm{eqp},1}=0$ and $n_{\mathrm{eqp},2}=8$, and we thus get $(\check{p}_{\mathrm{max}}+1+n_{\mathrm{eqp},2})\times(\check{p}_{\mathrm{max}}+1+n_{\mathrm{eqp},2})=11\times 11$ quadrature in each sub-element around the source point, and $(\check{p}_\upxi+1+n_{\mathrm{eqp},1})\times(\check{p}_\upeta+1+n_{\mathrm{eqp},1})=3\times 3$ in the remaining elements. The figures to the right are the new integration routine presented in this work with $s_1=1$.}
	\label{Fig3:Quadrature1}
\end{figure}

\begin{figure}
	\begin{subfigure}[t]{\textwidth}
		\centering
		\includegraphics[width=0.49\textwidth]{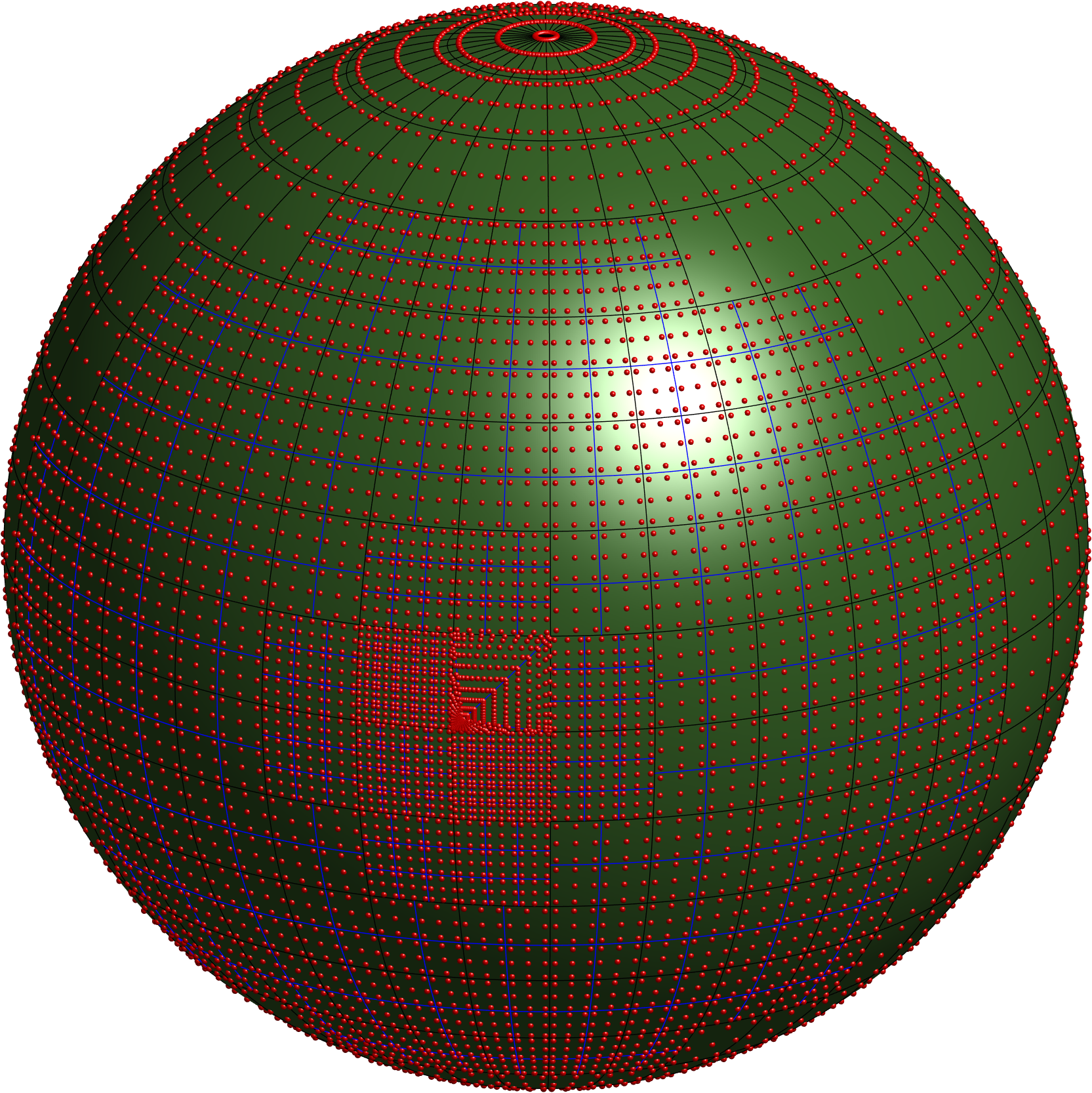}%
		\hspace*{0.02\textwidth}%
		\includegraphics[width=0.49\textwidth]{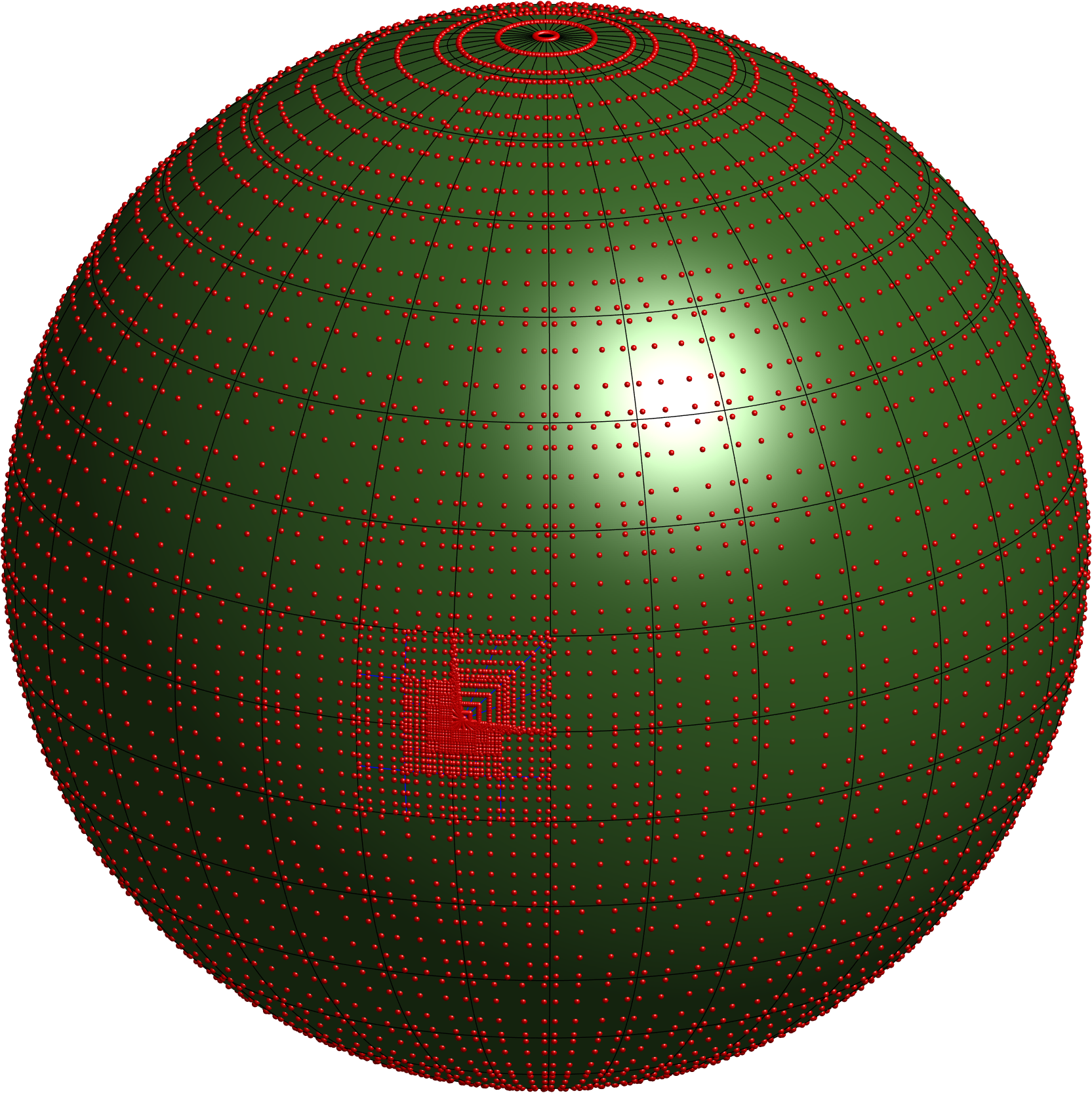}
		\caption{The $392^{\mathrm{th}}$ collocation point is at the center of an element which is divided into four (triangular) sub elements. The total number of quadrature points for the left and right figures are 8314 and 7371, respectively.}
	\end{subfigure}
	\par\bigskip
	\par\bigskip
	\begin{subfigure}[t]{\textwidth}
		\centering
		\includegraphics[width=0.49\textwidth]{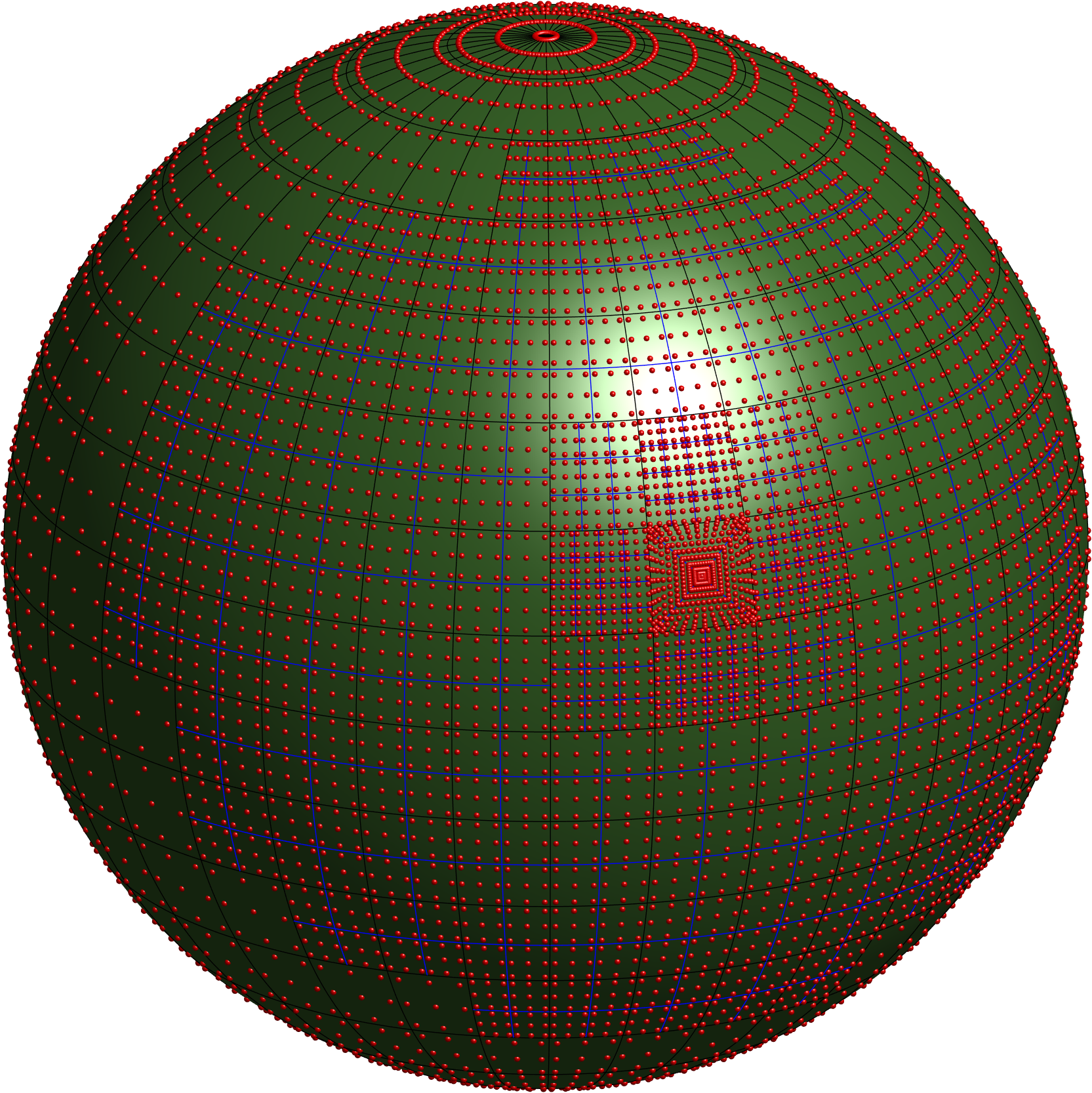}%
		\hspace*{0.02\textwidth}%
		\includegraphics[width=0.49\textwidth]{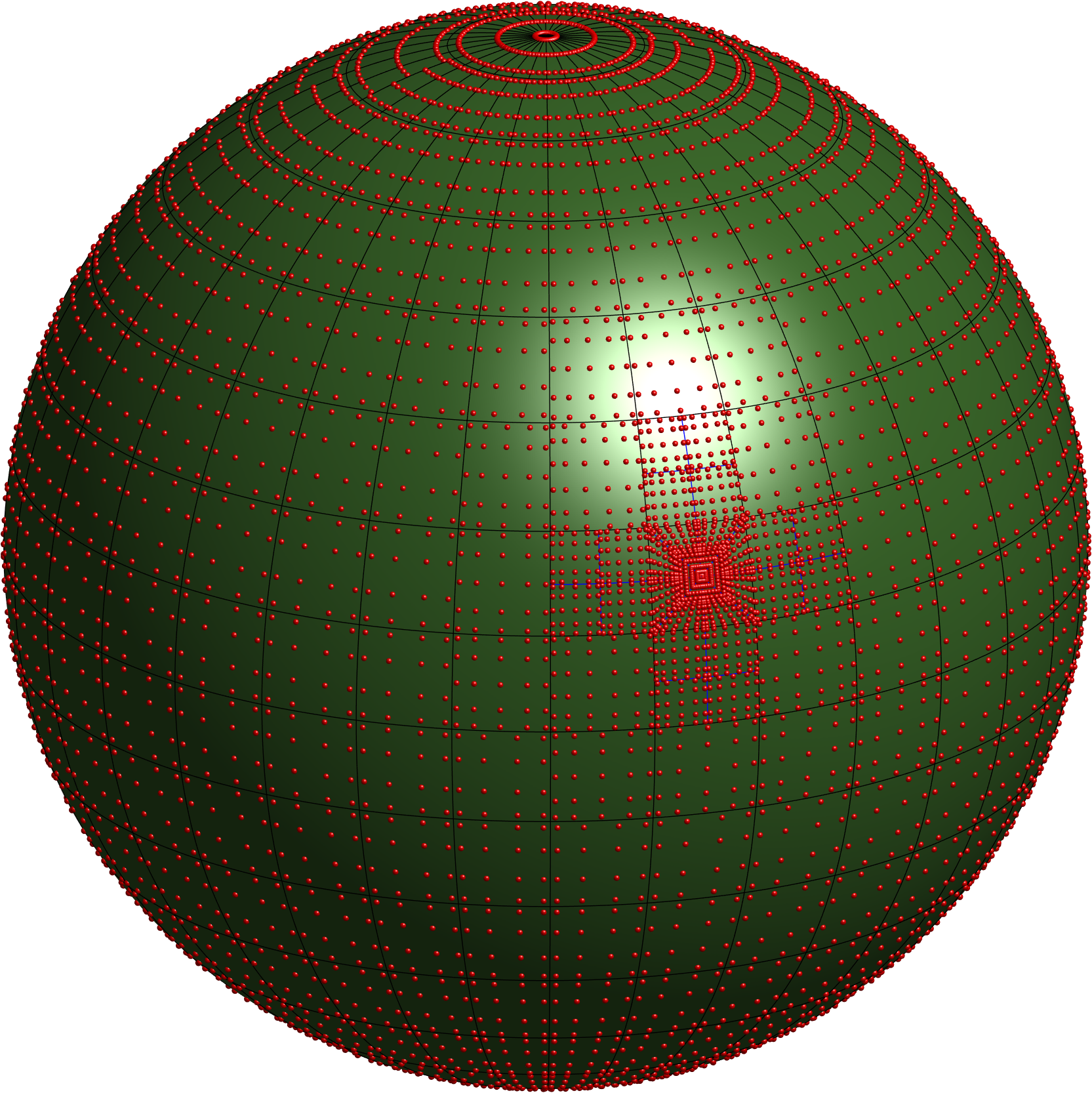}
		\caption{The source point is the corner quadrature point for Galerkin formulations inside an element which is divided into four (triangular) sub elements. The total number of quadrature points for the left and right figures are 8629 and 7956, respectively.}
	\end{subfigure}
	\caption{\textbf{Numerical evaluation of the boundary integral equations}: The figures to the left are the integration procedure in~\cite{Simpson2014aib} (with $s_1=2$). The sub-element divisions are here shown by blue lines (the black lines are the element edges). The red points are the quadrature points. Here, $n_{\mathrm{eqp},1}=0$ and $n_{\mathrm{eqp},2}=8$, and we thus get $(\check{p}_{\mathrm{max}}+1+n_{\mathrm{eqp},2})\times(\check{p}_{\mathrm{max}}+1+n_{\mathrm{eqp},2})=11\times 11$ quadrature in each sub-element around the source point, and $(\check{p}_\upxi+1+n_{\mathrm{eqp},1})\times(\check{p}_\upeta+1+n_{\mathrm{eqp},1})=3\times 3$ in the remaining elements. The figures to the right are the new integration routine presented in this work with $s_1=1$.}
	\label{Fig3:Quadrature2}
\end{figure}

In~\Cref{Fig3:Quadrature1,Fig3:Quadrature2} the locations of the quadrature points are illustrated on the third uniform mesh refinement of the coarse mesh in~\Cref{Fig3:sphere} (with $\check{p}=2$). 
\section{Numerical examples}
\label{Sec3:resultsDisc}
Acoustic scattering problems on a sphere are investigated in the following. These problems possess analytic solutions~\cite{Venas2019e3s} and are for this reason often used to verify numerical methods in acoustic scattering, e.g.~\cite{Gerdes1996so3,Ihlenburg1998fea,Simpson2014aib,Gerdes1998tcv,Gerdes1999otp,Coox2017aii}. In order to analyze convergence properties of IGABEM we also consider a torus, which can be represented by NURBS of polynomial order $\check{p}\geq 2$ with no poles in the parametrization. Also, a cube geometry will be investigated to check the behavior of the BIEs at $G^0$-geometries. We then continue be analyzing the BeTSSi\footnote{Benchmark Target Strength Simulation.} submarine. Before we consider the rigid scattering problem on this complex geometry, we present the method of manufactured solution. This method enables us to get some quality insurance of the underlying mesh to be used in the full scattering problem. Moreover, to some extent, the method can be used for quality insurance of the numerical solution of the scattering problem. Together with the benchmark problem on the sphere, these methods yield a solid basis for testing the correctness of the implemented code.

In this work, the test setting is chosen so that the present approach can be compared to other methods. In particular, the scattering on a rigid sphere example and the torus example is found in~\cite{Simpson2014aib}. Scattering on the BeTSSi submarine has been addressed at three workshops in the past 18 years~\cite{Nolte2014bib}. FWG\footnote{Forschungsanstalt f\"{u}r Wasserschall und Geophysik.} initiated the first workshop in 2001 (held in Kiel 2002) and delivered the generic BeTSSi submarine (for which the outer hull is described in~\Cref{Sec3:BeTSSi_description}). The second workshop took place in Kiel in 2014 and the third in the Hague in 2016. The best of these results will be used as reference solutions in this work. Additionally, we create our own reference simulations using \COMSOL~\cite{ComsolV54}. This benchmarking exercise is a crucial step to obtain reliable solutions for even more complex models.

The aim of these numerical examples is to investigate the approximability of IGABEM and its formulations. Moreover, we aim to establish highly accurate solutions for the BeTSSi submarine for benchmarking purposes and compare the accuracy and computational complexity of these results to existing simulations.

With the use of the Galerkin method the following quasi-optimal error estimate exists for the BEM~\cite[Theorem 2.49]{Chandler_Wilde2012nab} (with the Burton--Miller formulation)
\begin{equation}\label{Eq3:aprioriErrorEstimate}
	\|p-p_h\|_{L^2(\Gamma)} \leq C_1 \inf_{q_h\in V_h}\|p-q_h\|_{L^2(\Gamma)} \leq C_2 (hk)^{\check{p}+1}
\end{equation}
where $V_h$ is the finite dimensional subspace in which the solution is sought and the constants $C_1$ and $C_2$ may depend on the analytic solution $p$, the boundary $\Gamma$ and the wave number $k$. In this work we also aim to give numerical evidence for similar estimates for the other BEM formulations.

The simulations are based on the ASIGA\footnote{The ASIGA (Acoustic Scattering with IsoGeometric Analysis) library can be found at \href{https://github.com/Zetison/ASIGA}{this GiT-repository}.} library written in \MATLAB~\cite{MatlabR2019a}. The integration is here vectorized over the quadrature points, such that the effect of increasing the number of quadrature points is of less significance due to the efficiency of vectorization in \MATLAB. For this reason, we take the liberty of over integration the BIEs without suffering to much from computational cost. For optimization purposes, the library could be written in \ccpp which would require an accuracy-cost tradeoff study in this respect. Additionally, acceleration techniques exist for the boundary element method which have not been implemented in the ASIGA library. We refer to \cite{Dolz2016aib,Dolz2018afi,Beer2008tbe} for details. These optimizations are suggested as future work.

The BIE formulations listed in~\Cref{Tab3:BIEs} will be investigated both in terms of approximability and the presence of fictitious eigenfrequencies.
\begin{table}
	\centering
	\caption{Overview of the boundary integral equation (BIE) formulations considered in this work.}
	\label{Tab3:BIEs}
	\begin{tabular}{l l l}
		\toprule
		Abbreviation & Name & Definition\\
		\hline
		CBIE & Conventional BIE & \Cref{Eq3:CBIE}\\
		RCBIE1 & The first regularized CBIE & \Cref{Eq3:RCBIE1}\\
		RCBIE2 & The second regularized CBIE & \Cref{Eq3:RCBIE2}\\
		RCBIE3 & The third regularized CBIE & \Cref{Eq3:RCBIE3}\\
		HBIE & Hypersingular BIE & \Cref{Eq3:HBIE}\\
		BM & Burton--Miller & \Cref{Eq3:BM}\\
		\bottomrule
	\end{tabular}
\end{table}

The meshes will be generated from a coarse CAD model mesh (for example \Cref{Fig3:parm1} for the sphere) with mesh number $m=1$. We shall denote by ${\cal M}_{m,\check{p},\check{k}}^{\textsc{igabem}}$, mesh number $m$ with polynomial order $\check{p}$ and continuity $\check{k}$ across element boundaries\footnote{Except for (potentially) some $C^0$ lines in the initial CAD geometry.}. For the corresponding FEM meshes we denote by ${\cal M}_{m,\check{p},\mathrm{s}}^{\textsc{fembem}}$ and ${\cal M}_{m,\check{p},\mathrm{i}}^{\textsc{fembem}}$ the subparametric and isoparametric FEM meshes, respectively. These meshes are constructed by the procedure outlined in~\cite[p. 191]{Venas2018iao}.

\subsection{Pulsating sphere}
\label{Sec3:pulsatingSphere}
\begin{figure}
	\begin{subfigure}{0.49\textwidth}
		\centering
		\includegraphics[width=0.8\textwidth]{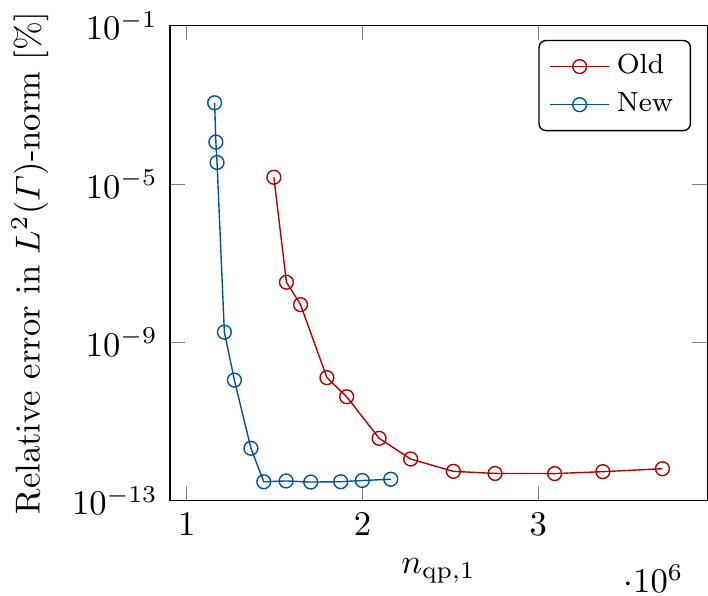}
		\caption{CCBIE}
	\end{subfigure}%
	\hspace*{0.02\textwidth}%
	\begin{subfigure}{0.49\textwidth}
		\centering
		\includegraphics[width=0.8\textwidth]{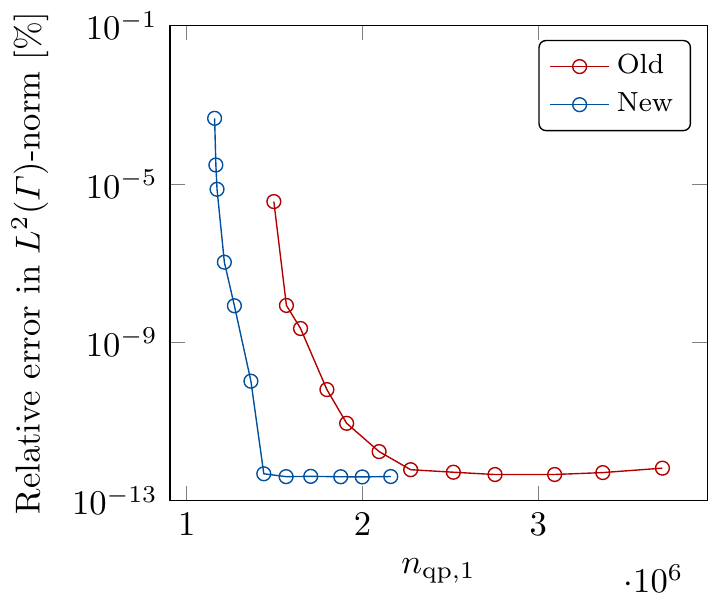}
		\caption{CRCBIE3}
	\end{subfigure}
	\par\bigskip
	\par\bigskip
	\begin{subfigure}{0.49\textwidth}
		\centering
		\includegraphics[width=0.8\textwidth]{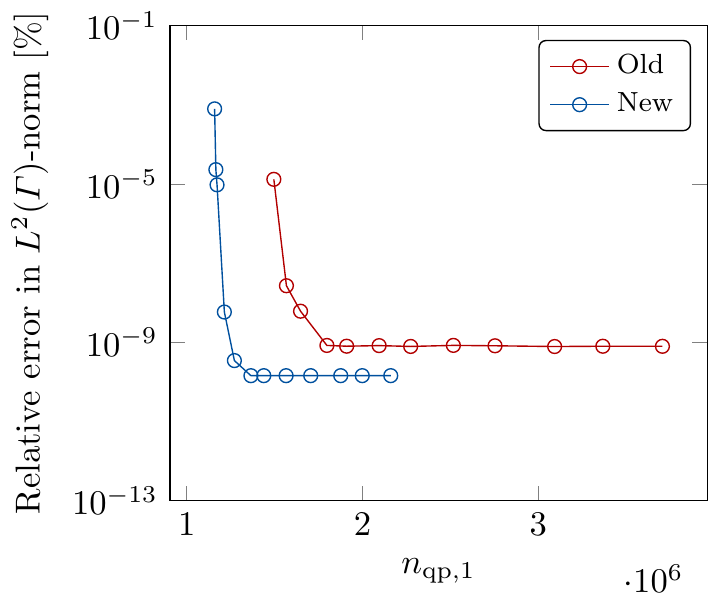}
		\caption{CHBIE}
	\end{subfigure}%
	\hspace*{0.02\textwidth}%
	\begin{subfigure}{0.49\textwidth}
		\centering
		\includegraphics[width=0.8\textwidth]{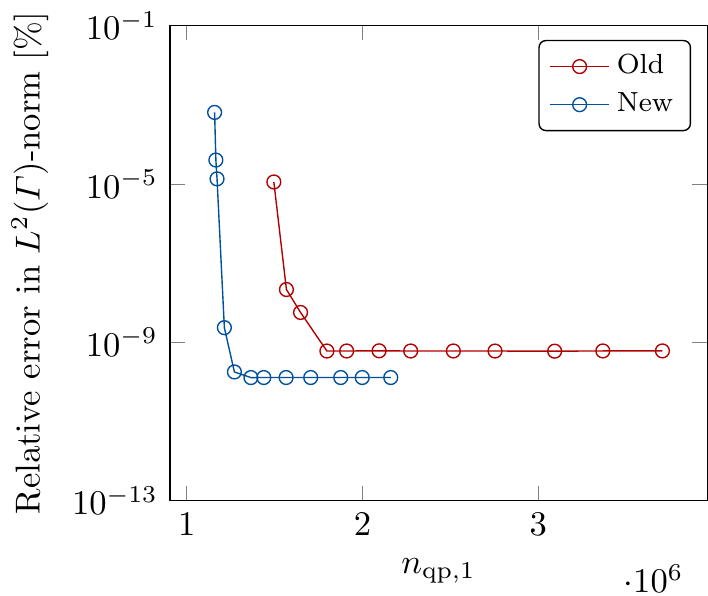}
		\caption{CBM}
	\end{subfigure}
	\par\bigskip
	\par\bigskip
	\begin{subfigure}{0.49\textwidth}
		\centering
		\includegraphics[width=0.8\textwidth]{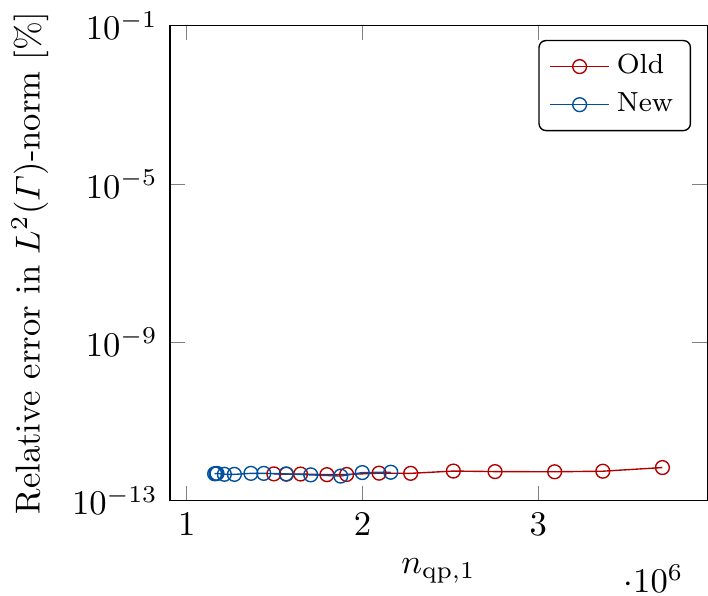}
		\caption{CRCBIE1}
	\end{subfigure}%
	\hspace*{0.02\textwidth}%
	\begin{subfigure}{0.49\textwidth}
		\centering
		\includegraphics[width=0.8\textwidth]{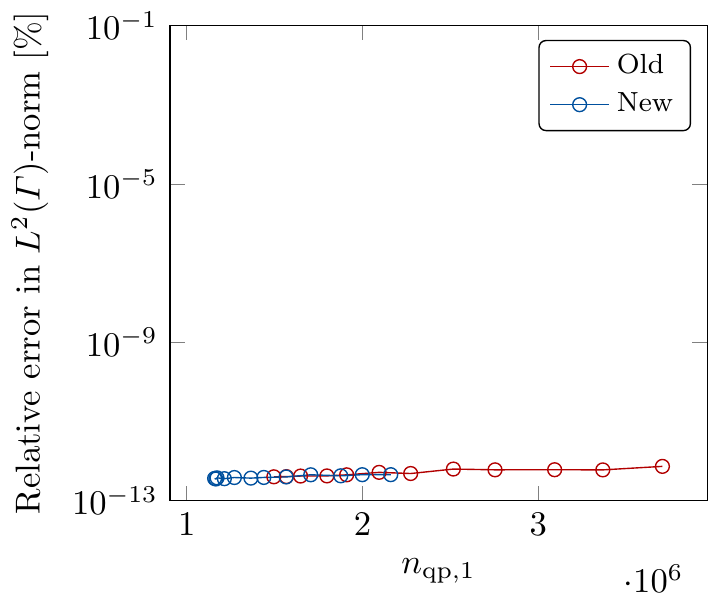}
		\caption{CRCBIE2}
	\end{subfigure}
	\caption{\textbf{Pulsating sphere}: Surface error as a function of the total number of quadrature points $n_{\mathrm{qp,1}}$ at $kR_0=1$.  The old adaptive quadrature scheme presented by Simpson in~\cite{Simpson2014aib} is compared to the new adaptive quadrature scheme presented in this work. The sample points correspond to $s_1\in\{1,2,\dots,12\}$ and $s_1\in\{1,2,\dots,12\}/5$ for the old and new method, respectively.}
	\label{Fig3:Simpson_PS_1}
\end{figure}

\begin{figure}
	\begin{subfigure}{0.49\textwidth}
		\centering
		\includegraphics[width=0.8\textwidth]{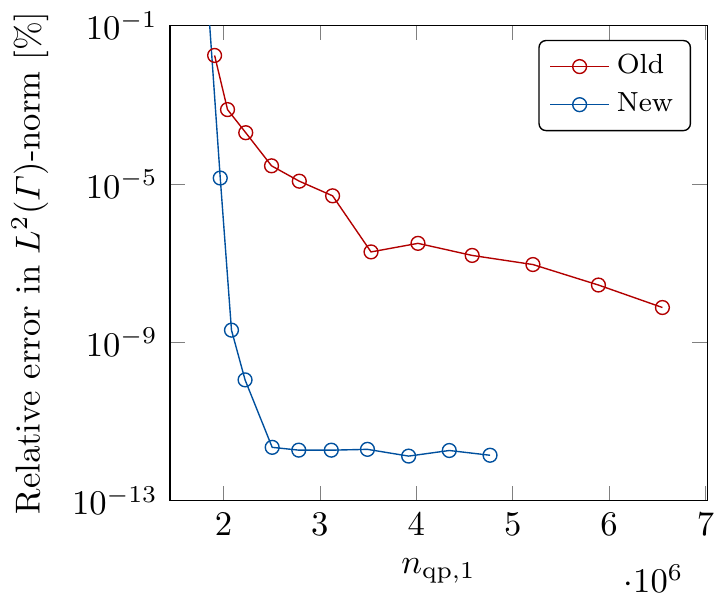}
		\caption{GCBIE}
	\end{subfigure}%
	\hspace*{0.02\textwidth}%
	\begin{subfigure}{0.49\textwidth}
		\centering
		\includegraphics[width=0.8\textwidth]{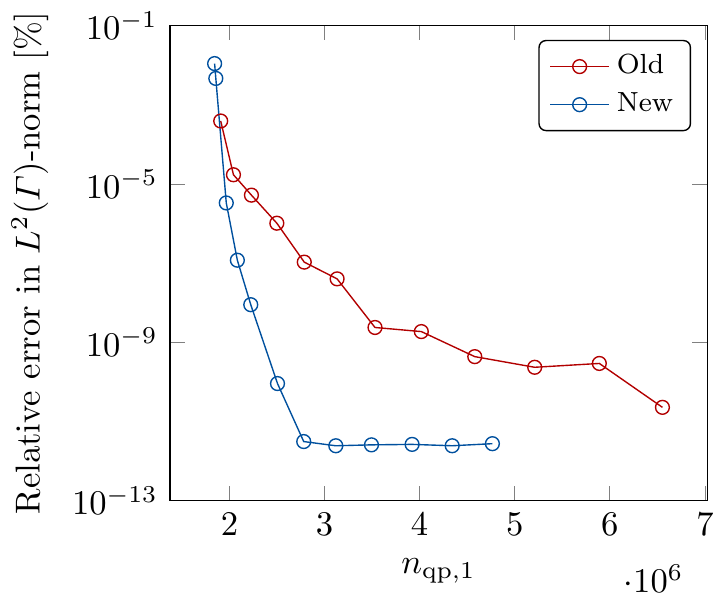}
		\caption{GRCBIE3}
	\end{subfigure}
	\par\bigskip
	\par\bigskip
	\begin{subfigure}{0.49\textwidth}
		\centering
		\includegraphics[width=0.8\textwidth]{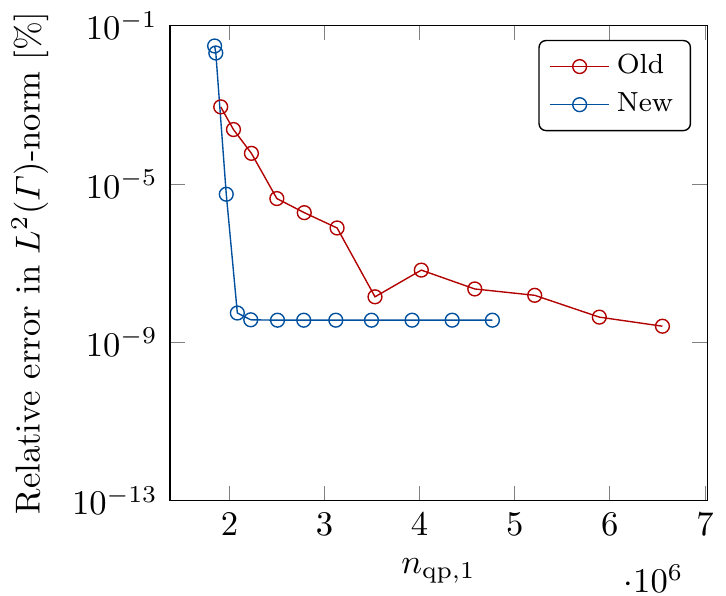}
		\caption{GHBIE}
	\end{subfigure}%
	\hspace*{0.02\textwidth}%
	\begin{subfigure}{0.49\textwidth}
		\centering
		\includegraphics[width=0.8\textwidth]{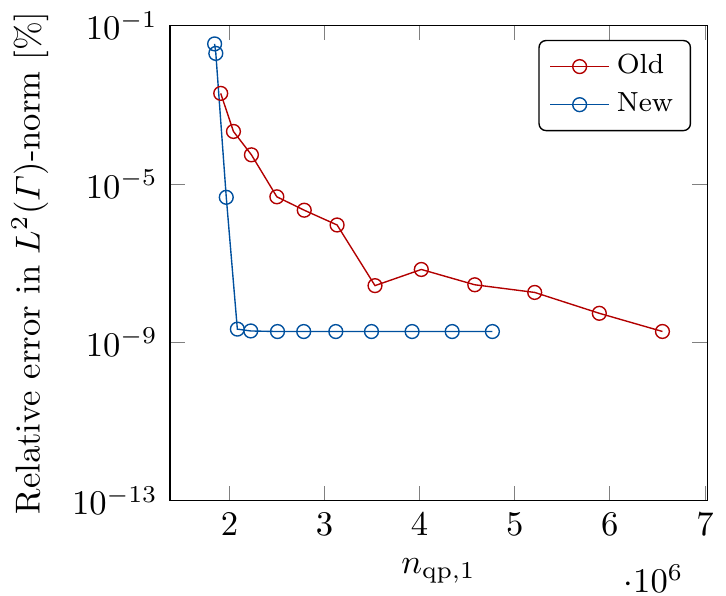}
		\caption{GBM}
	\end{subfigure}
	\par\bigskip
	\par\bigskip
	\begin{subfigure}{0.49\textwidth}
		\centering
		\includegraphics[width=0.8\textwidth]{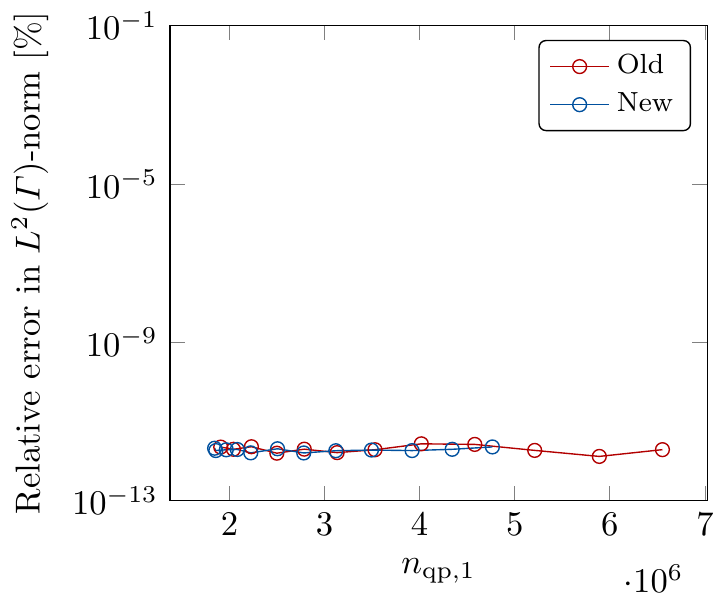}
		\caption{GRCBIE1}
	\end{subfigure}%
	\hspace*{0.02\textwidth}%
	\begin{subfigure}{0.49\textwidth}
		\centering
		\includegraphics[width=0.8\textwidth]{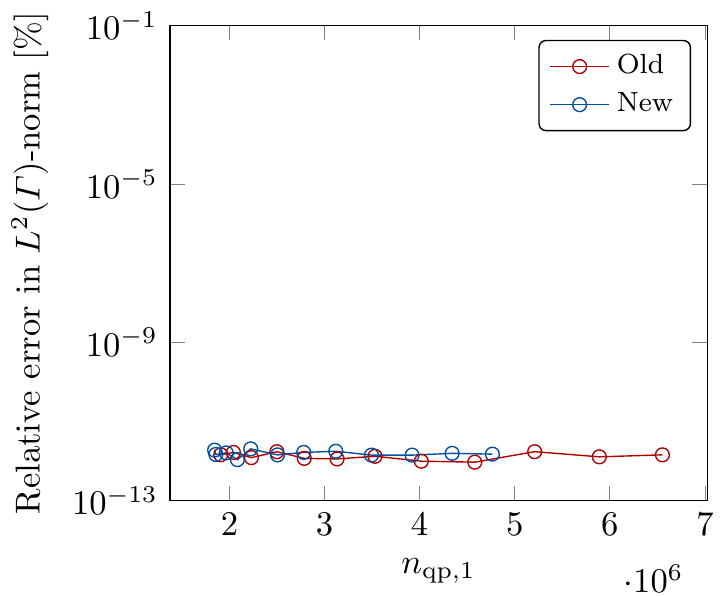}
		\caption{GRCBIE2}
	\end{subfigure}
	\caption{\textbf{Pulsating sphere}: Surface error as a function of the total number of quadrature points $n_{\mathrm{qp,1}}$ at $kR_0=1$.  The old adaptive quadrature scheme presented by Simpson in~\cite{Simpson2014aib} is compared to the new adaptive quadrature scheme presented in this work. The sample points correspond to $s_1\in\{1,2,\dots,12\}$ and $s_1\in\{1,2,\dots,12\}/5$ for the old and new method, respectively.}
	\label{Fig3:Simpson_PS_2}
\end{figure}

Consider a pulsating unit sphere centered at the origin (cf. \cite{Simpson2014aib,Zheng2015itb}) with analytic solution given by
\begin{equation}
	p(\vec{x}) = \frac{\euler^{\imag kR}}{4\PI R},\quad R=|\vec{x}|,\quad \vec{x}\in\Omega^+
\end{equation}
and with the (constant) Neumann condition
\begin{equation}
	g(\vec{x}) = \frac{\euler^{\imag k}}{4\PI}(\imag k -1),\quad \vec{x}\in\Gamma.
\end{equation}
This problem serves as a patch test for IGA as the analytic solution lies in the numerical solution space ($p(\vec{x})$ is constant at $\Gamma$). Contrary to FEM with affine mappings, (proper) Gaussian quadrature does not integrate the integrals in BEM exactly. Therefore, this example may be used to give some indication of the quality of the integration procedure. In~\Cref{Fig3:Simpson_PS_1,Fig3:Simpson_PS_2} we compare the two adaptive quadrature schemes (described in~\Cref{Sec3:numericalQuad}), where we set $n_{\mathrm{eqp},2}=100$ to avoid error originating from the integration over the element containing the source points. The $L^2$-error of the numerical solution is here plotted against $n_{\mathrm{qp},1}$; the total number of quadrature points, excluding quadrature points in elements containing the source point. The simulations are done on the coarsest mesh of the second NURBS parametrization in~\Cref{Fig3:parm2} (with $\check{p}=4$). The BM and HBIE formulations (for both collocation and Galerkin) have more round-off errors and are for this reason further away from machine epsilon precision results compared to the other formulations. In all cases, the new adaptive quadrature scheme obtains better results. Interestingly CBIE obtains slightly better results using the new adaptive quadrature scheme compared to RCBIE3, the latter being the regularized version of the former. This might be due to the reduction of symmetry in the RCBIE3 compared to CBIE for this problem.

Note that for this problem using RCBIE1 or RCBIE2 (\Cref{Eq3:RCBIE1,Eq3:RCBIE2}), results with machine epsilon precision are always obtained since the integrands are zero. This is due to the spherical symmetry of the problem and the functions involved.

Based on this study, a proper choice for the parameter $s_1$ is $s_1=1.4$ for the new adaptive method. If not otherwise stated, we shall use $s_1=1.4$ and $n_{\mathrm{eqp},2}=50$, which in most cases results in over integration. As was mentioned before, the cost of this is not significant due to the current implementation in \MATLAB.

\subsection{Rigid scattering on a sphere} 
Consider a plane wave, with the direction of incidence given by
\begin{equation}\label{Eq3:d_s}
	\vec{d}_{\mathrm{s}} = -\begin{bmatrix}
		\cos\beta_{\mathrm{s}}\cos\alpha_{\mathrm{s}}\\
		\cos\beta_{\mathrm{s}}\sin\alpha_{\mathrm{s}}\\
		\sin\beta_{\mathrm{s}}
	\end{bmatrix},
\end{equation}
with\footnote{The angles $\alpha$ and $\beta$ are the so-called aspect and elevation angle, respectively. Note that the aspect angle is equal to the spherical coordinate $\varphi$ (the azimuth angle).} $\alpha_{\mathrm{s}} = \ang{240}$ and $\beta_{\mathrm{s}} = \ang{30}$, scattered by a rigid sphere with radius $R_0=\SI{1}{m}$. 
\begin{figure}
	\centering
	\begin{subfigure}[t]{0.3\textwidth}
		\includegraphics[width=0.9\textwidth]{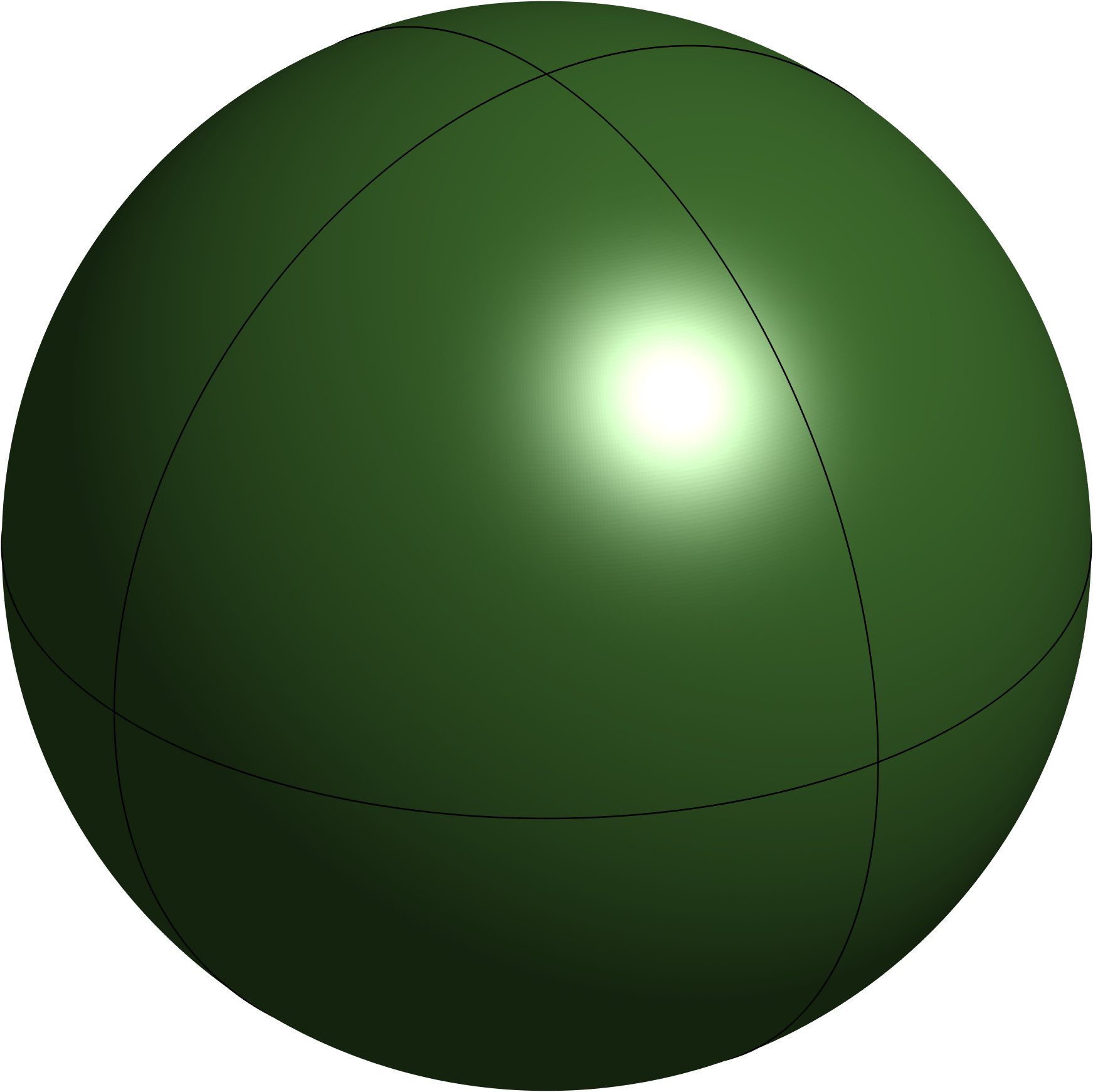}
		\caption{Parametrization 1.}
		\label{Fig3:parm1}
	\end{subfigure} 
	~
	\begin{subfigure}[t]{0.3\textwidth}
		\includegraphics[width=0.9\textwidth]{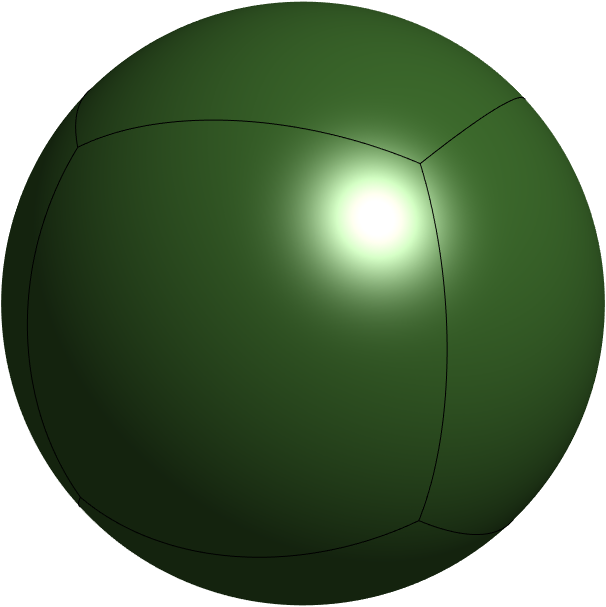}
		\caption{Parametrization 2.}
		\label{Fig3:parm2}
	\end{subfigure} 
	\caption{Two exact NURBS parametrizations of the sphere. Parametrization 1 uses a single patch with 8 elements of degree $\check{p}\geq 2$ while parametrization 2 uses 6 patches of degree $\check{p}\geq 4$. Parametrization 1 is described in~\Cref{Sec3:NURBSsphere1} and parametrization 2 is described in \Cref{Sec3:NURBSsphere2}.}
	\label{Fig3:SphericalShellParametrizations}
\end{figure}

For the rigid scattering problems considered in this work, the error is computed of $p_{\mathrm{tot}}$ and the best approximation (BA) is obtained by performing an $L^2$-projection of $p_{\mathrm{tot}}$ onto the discretized solution space.

\begin{figure}
	\centering    
	\begin{subfigure}[t]{0.49\textwidth}
		\centering
		\includegraphics{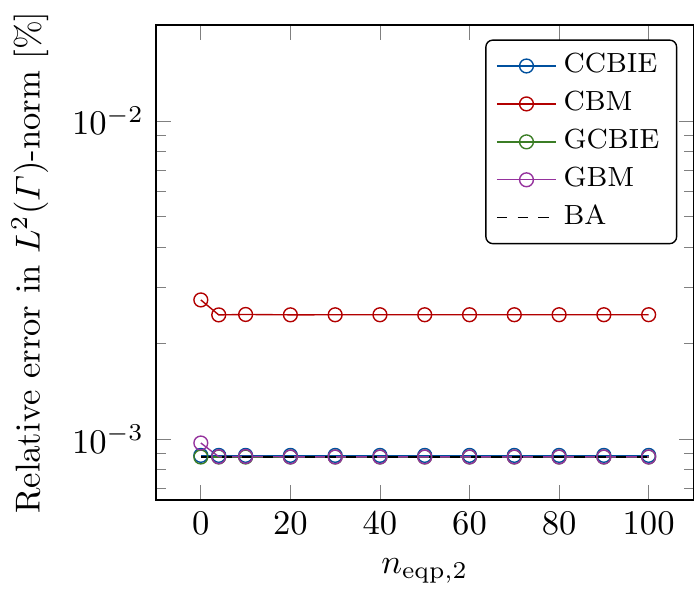}
	\end{subfigure}%
	\hspace*{0.02\textwidth}%
	\begin{subfigure}[t]{0.49\textwidth}
		\centering
		\includegraphics{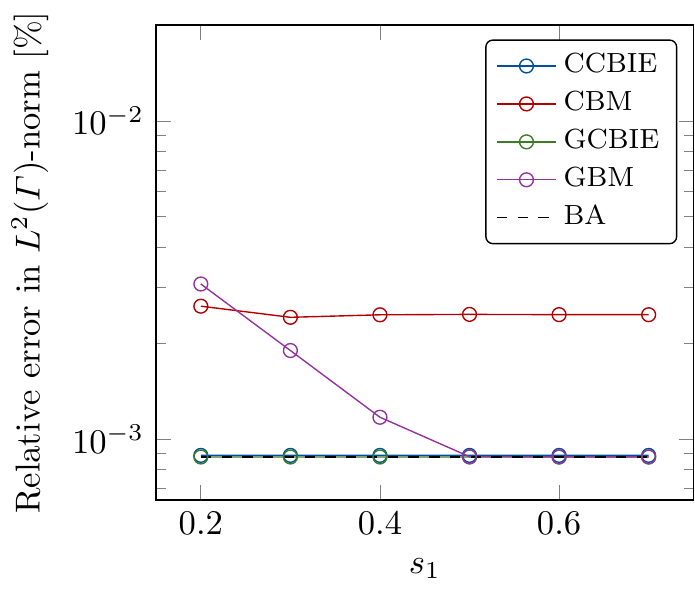}
	\end{subfigure}
	\caption{\textbf{Rigid scattering on a sphere}: Surface error as a function of the parameters $n_{\mathrm{eqp},2}$ and $s_1$ to the left and right, respectively, on the mesh ${\cal M}_{5,2,1}^{\textsc{igabem}}$.}
	\label{Fig3:S1_qp_M5p2}
\end{figure}

\begin{figure}
	\centering
	\begin{subfigure}[t]{0.49\textwidth}
		\centering
		\includegraphics{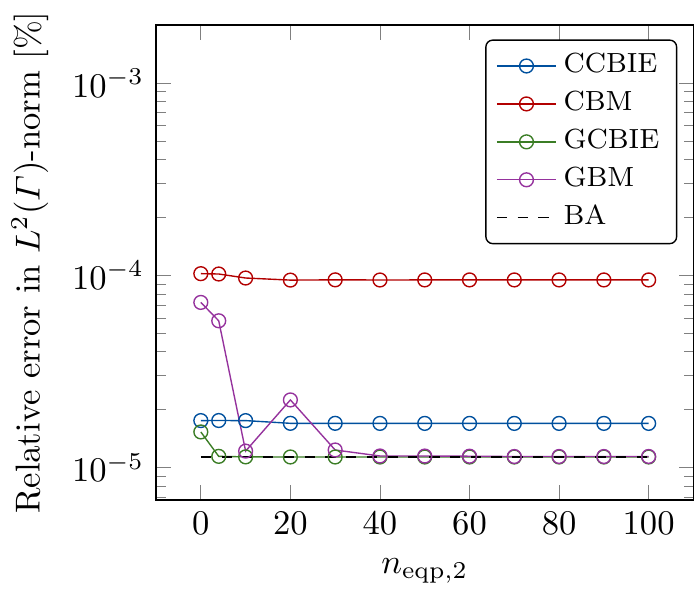}
	\end{subfigure}%
	\hspace*{0.02\textwidth}%
	\begin{subfigure}[t]{0.49\textwidth}
		\centering
		\includegraphics{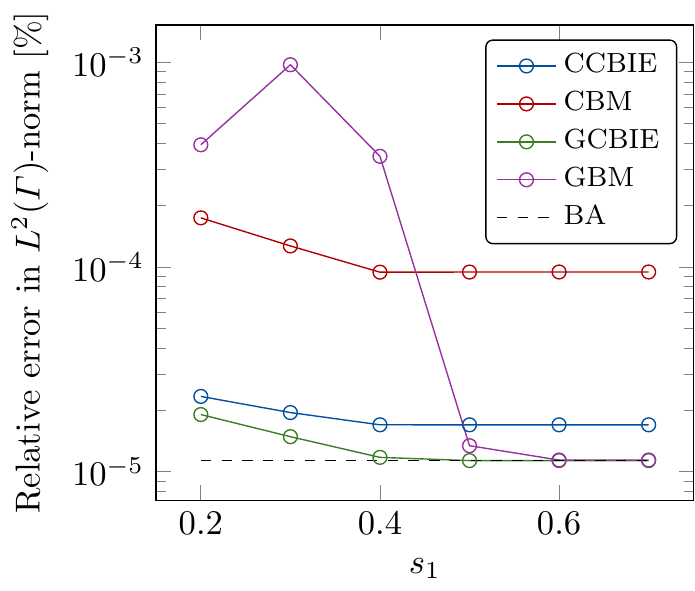}
	\end{subfigure}
	\caption{\textbf{Rigid scattering on a sphere}: Surface error as a function of the parameters $n_{\mathrm{eqp},2}$ and $s_1$ to the left and right, respectively, on the mesh ${\cal M}_{4,5,4}^{\textsc{igabem}}$.}
	\label{Fig3:S1_qp_M4p5}
\end{figure}

\begin{figure}
	\centering
	\begin{subfigure}[t]{0.49\textwidth}
		\centering
		\includegraphics{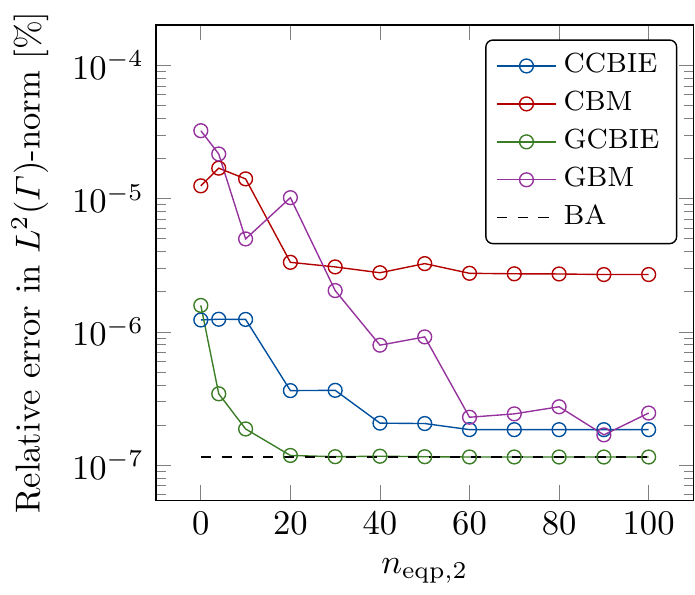}
	\end{subfigure}%
	\hspace*{0.02\textwidth}%
	\begin{subfigure}[t]{0.49\textwidth}
		\centering
		\includegraphics{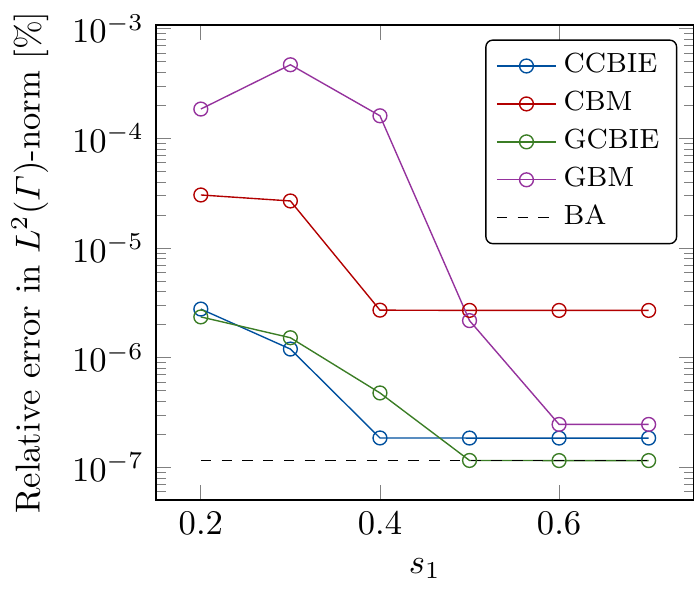}
	\end{subfigure}
	\caption{\textbf{Rigid scattering on a sphere}: Surface error as a function of the parameters $n_{\mathrm{eqp},2}$ and $s_1$ to the left and right, respectively, on the mesh ${\cal M}_{5,5,4}^{\textsc{igabem}}$.}
	\label{Fig3:S1_qp_M5p5}
\end{figure}

Continuing the study of numerical quadrature, we investigate the parameters $s_1$ and $n_{\mathrm{eqp},2}$ also for rigid scattering. The study for the parameter $s_1$ uses $n_{\mathrm{eqp},2}=100$ and the study for $n_{\mathrm{eqp},2}$ uses $s_1=0.7$. For FEM/IGA using $\check{p}+1$ quadrature points in each parametric direction in each element ensures accurate numerical integration regardless of the computational mesh. As can be observed from \Cref{Fig3:S1_qp_M5p2,Fig3:S1_qp_M4p5,Fig3:S1_qp_M5p5} this is not the case for BEM. Separate choices for the parameters $n_{\mathrm{eqp},2}$ and $s_1$ need to be made for each formulation. Contrary to FEM/IGA the optimal quadrature rule seems to be depending on $h$-refinement (not only $\check{p}$-refinement). Although the integrals in the CBIE formulation are regularized to contain no singular integrals, the parameter $s_1$ may still not be set to zero. This could be expected due to the gradients around the source points.


\begin{figure}
	\centering
	\begin{subfigure}[t]{\textwidth}
		\centering
		\includegraphics[width=\textwidth]{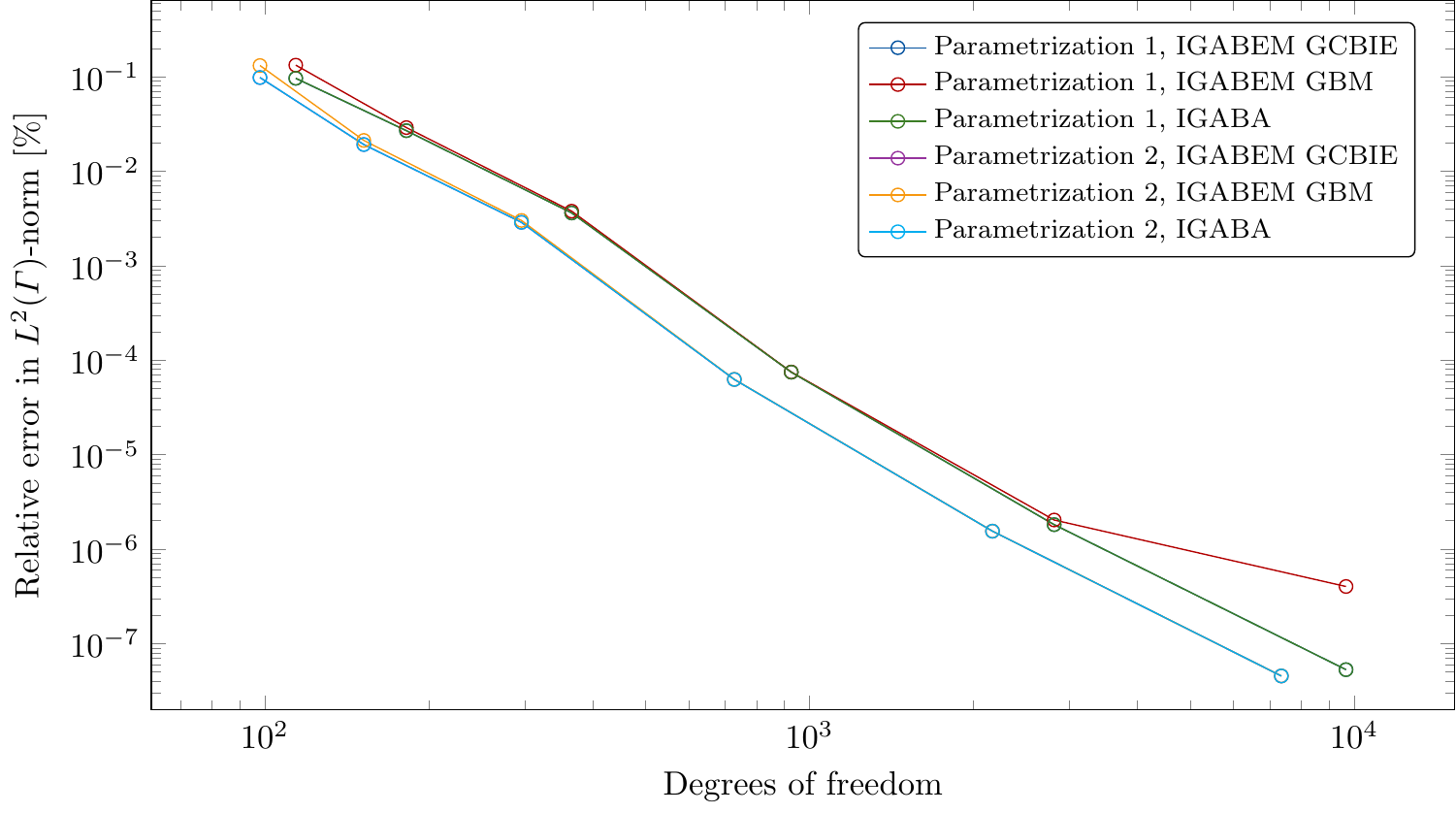}
		\caption{Galerkin formulations}
		\label{Fig3:S1parmCompG}
	\end{subfigure}
	\par\bigskip
	\par\bigskip
	\begin{subfigure}[t]{\textwidth}
		\centering
		\includegraphics[width=\textwidth]{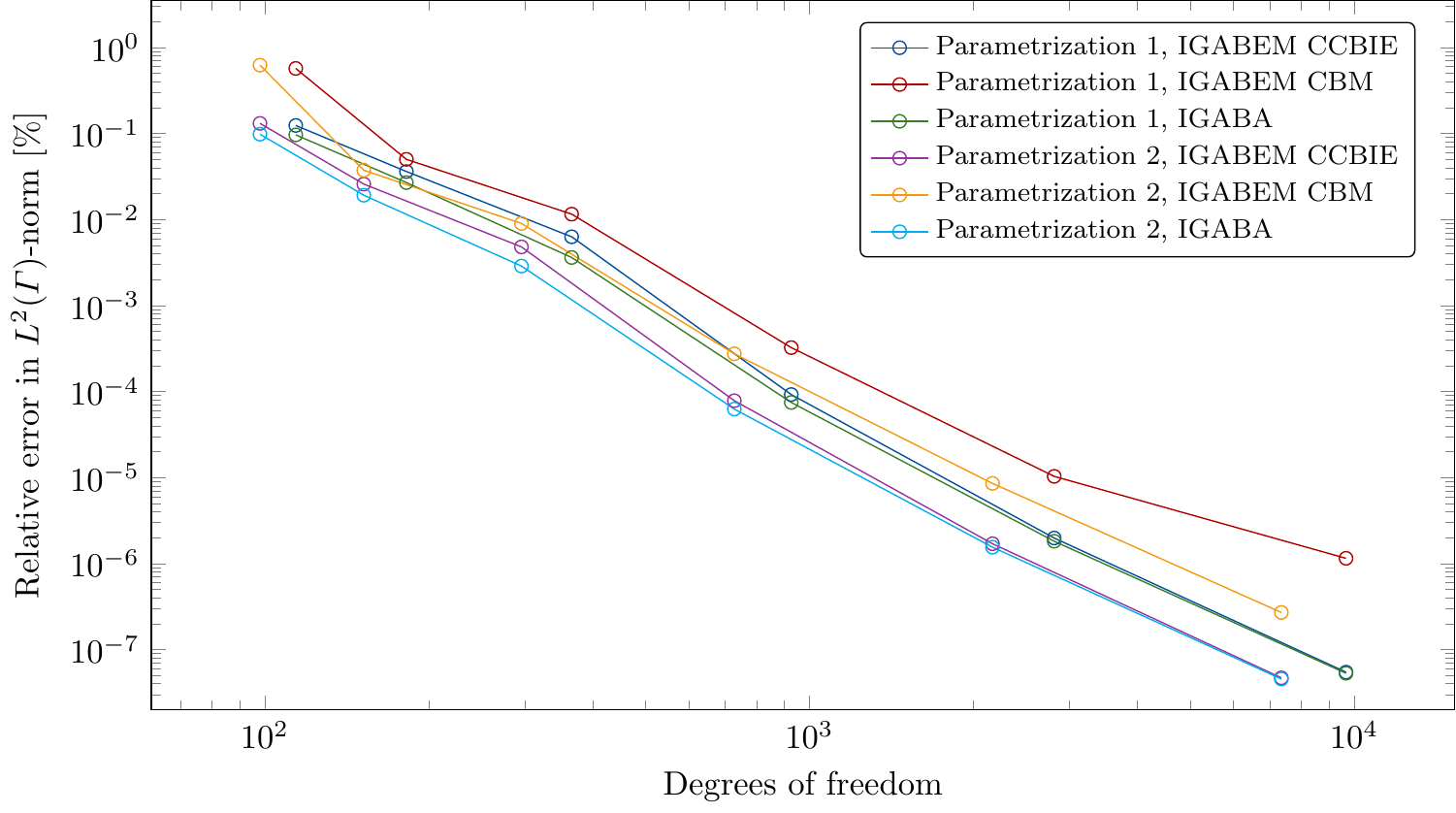}
		\caption{Collocation formulations}
		\label{Fig3:S1parmCompC}
	\end{subfigure}
	\caption{\textbf{Rigid scattering on a sphere}: Convergence analysis with $\check{p}=4$ and $kR_0=1$.}
	\label{Fig3:S1parmComp}
\end{figure}

\begin{figure}
	\centering
	\includegraphics{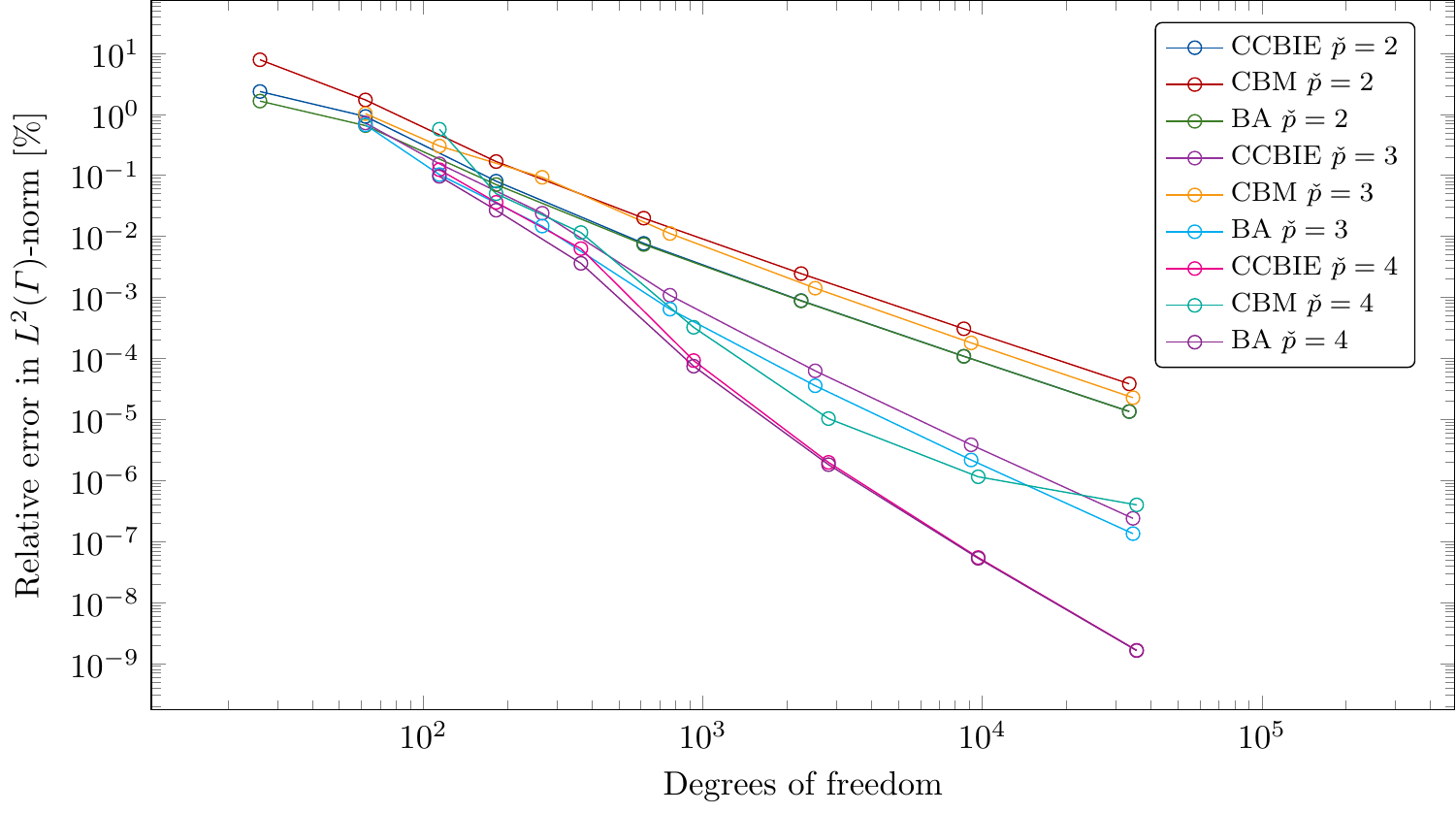}
	\caption{\textbf{Rigid scattering on a sphere}: Convergence analysis with $kR_0=1$.}
	\label{Fig3:S1_cgComp}
\end{figure}

For convenience we perturb the collocation points at the north and the south pole of the parametrization in \Cref{Fig3:parm1} in the HBIE and BM formulation for the ease of implementation. The perturbation is taken to be a distance $\frac12|\Diff\eta_e|/\check{p}_\upeta$ in the $\eta$-direction (in the parametric space), where $|\Diff\eta_e|$ is the element interval in the parametric domain in the $\eta$-direction. A similar strategy will be employed for the corresponding problematic areas on the BeTSSi submarine. This may be a sub optimal placement of collocation points, and as we can see from \Cref{Fig3:S1parmCompC}, the CBM formulation does not obtain the accuracy of the Galerkin formulation (\Cref{Fig3:S1parmCompG}). But this is also true for parametrization 2 (which contains no poles), and so this calls for an investigation of better placement of collocation points in general for the CHBIE and CBM than that of the Greville abscissae. The CBM formulation for parametrization 1 is visibly polluted by round-off errors similar to those seen in \Cref{Sec3:pulsatingSphere}.
In~\Cref{Fig3:S1_cgComp} we can observe that CBM loses one order of convergence for the odd degree $\check{p}=3$, which is similar to the effect discussed in~\cite{Gomez2016tvc}. However, this effect does not come into play in the same way for the CCBIE formulation, although it is still a significant difference between this simulation and the best approximation. This is in stark contrast to the CCBIE simulations of even degree which approaches the best approximation solution.

The plots in \Cref{Fig3:S1parmComp} also show the impact a sub optimal parametrization may have. Parametrization 1 has roughly 8\% higher errors compared to parametrization 2 in terms of degrees of freedom.

\begin{figure}
	\centering
	\includegraphics[width=\textwidth]{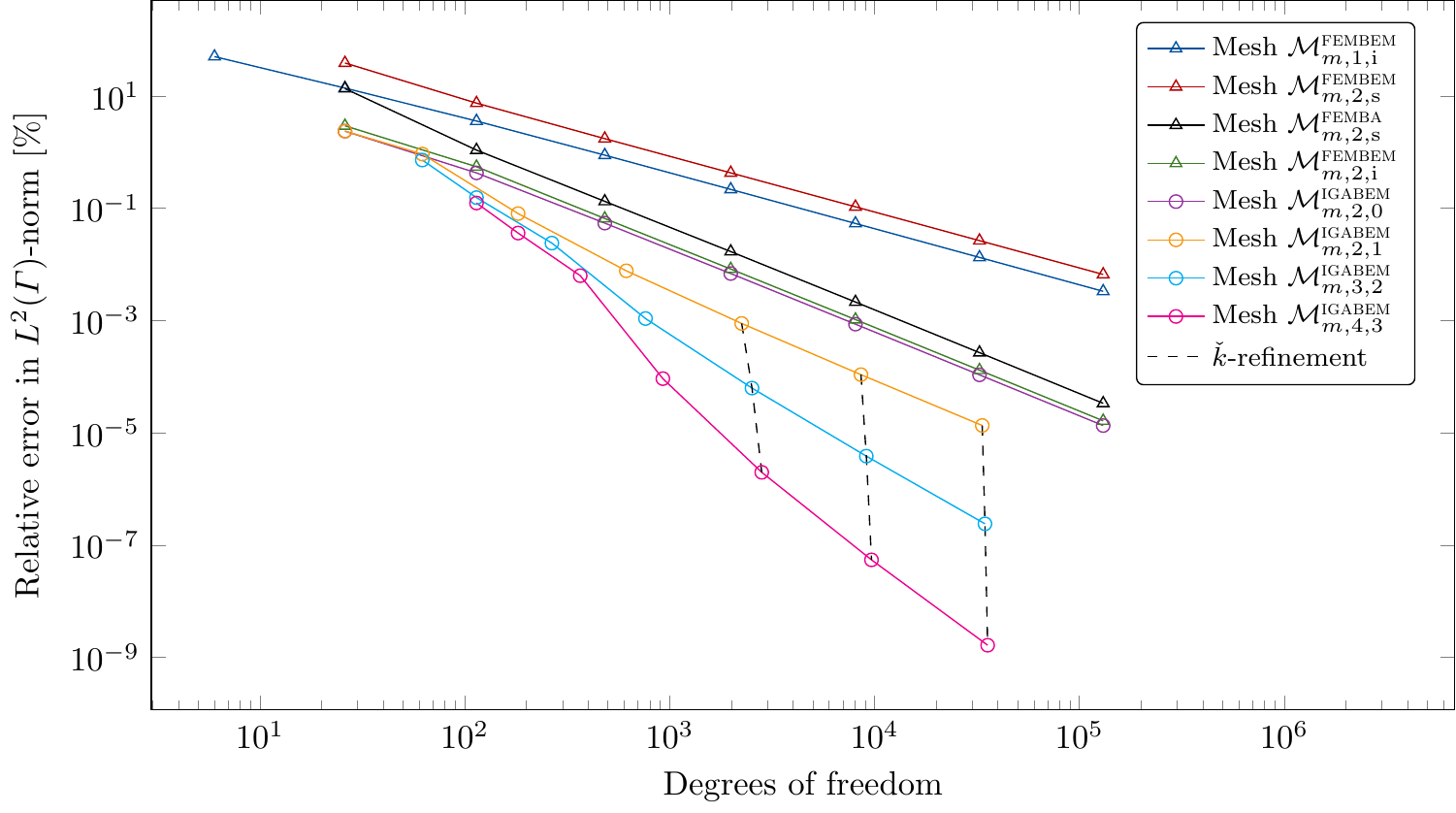}
	\caption{\textbf{Rigid scattering on a sphere}: Convergence analysis with the CCBIE formulation on parametrization 1 for $kR_0=1$.}
	\label{Fig3:convergenceAnalysis}
\end{figure}

\begin{figure}
	\centering
	\includegraphics[width=\textwidth]{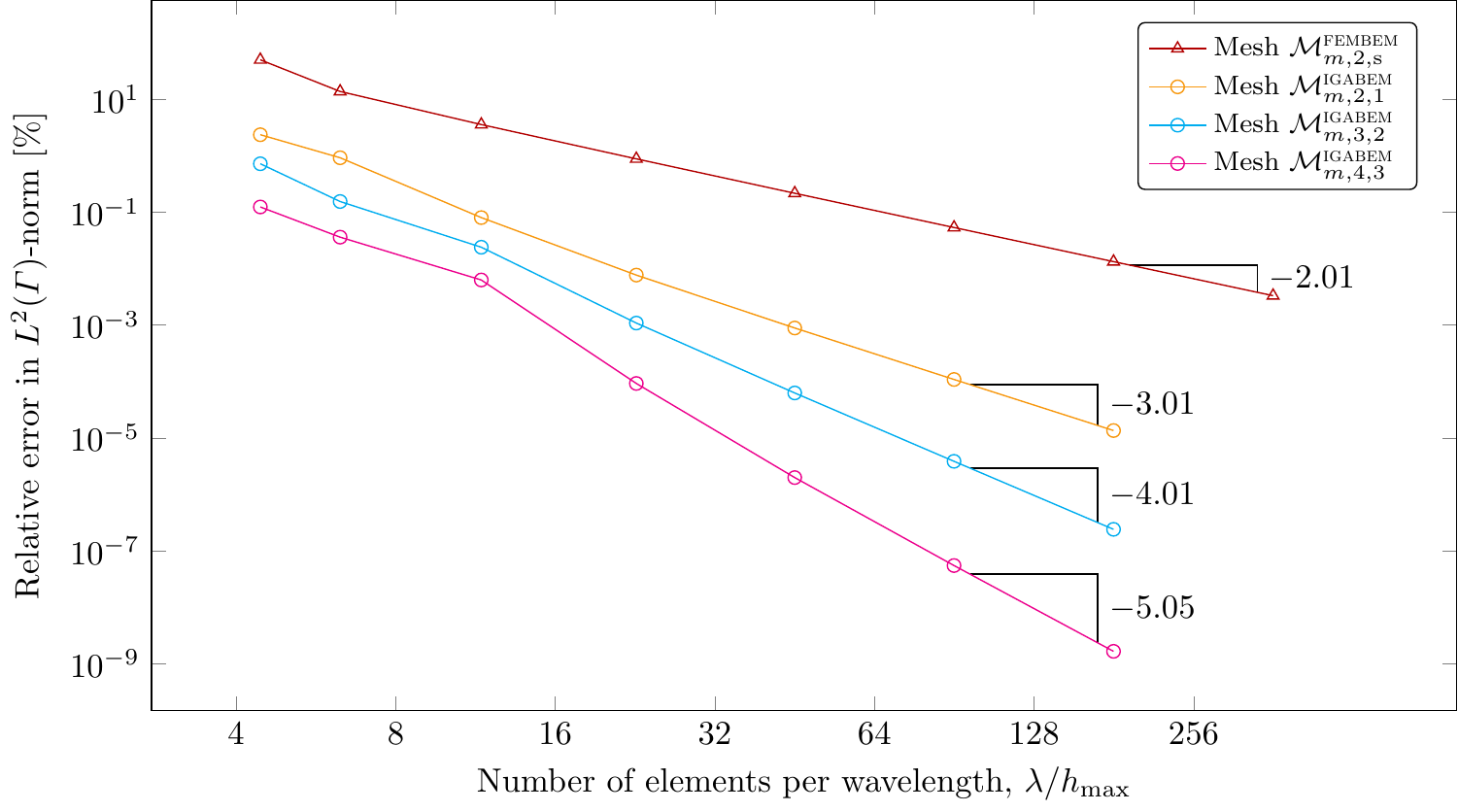}
	\caption{\textbf{Rigid scattering on a sphere}: Convergence analysis with the CCBIE formulation on parametrization 1 for $kR_0=1$.}
	\label{Fig3:convergenceAnalysis_nepw}
\end{figure}

In \Cref{Fig3:convergenceAnalysis} we compare the classical boundary element method (FEMBEM) with IGA. For the subparametric second order FEMBEM mesh a full convergence order (see \Cref{Fig3:convergenceAnalysis_nepw}) is lost in comparison with the best approximation for the same mesh (FEMBA). In fact, little is to be gained by increasing the polynomial order when using a linear approximation of the geometry. The exactness of the geometry is of less importance for isoparametric FEMBEM, which can be observed by comparing the results for mesh ${\cal M}_{m,2,\mathrm{i}}^{\textsc{fembem}}$ and mesh ${\cal M}_{m,2,0}^{\textsc{igabem}}$. Increasing the continuity ($\check{k}$-refinement) of the basis functions, however, improves the accuracy significantly as obtained for infinite isogeometric finite elements~\cite{Venas2018iao}.

\begin{figure}
	\centering
	\begin{subfigure}[t]{\textwidth}
		\includegraphics[width=\textwidth]{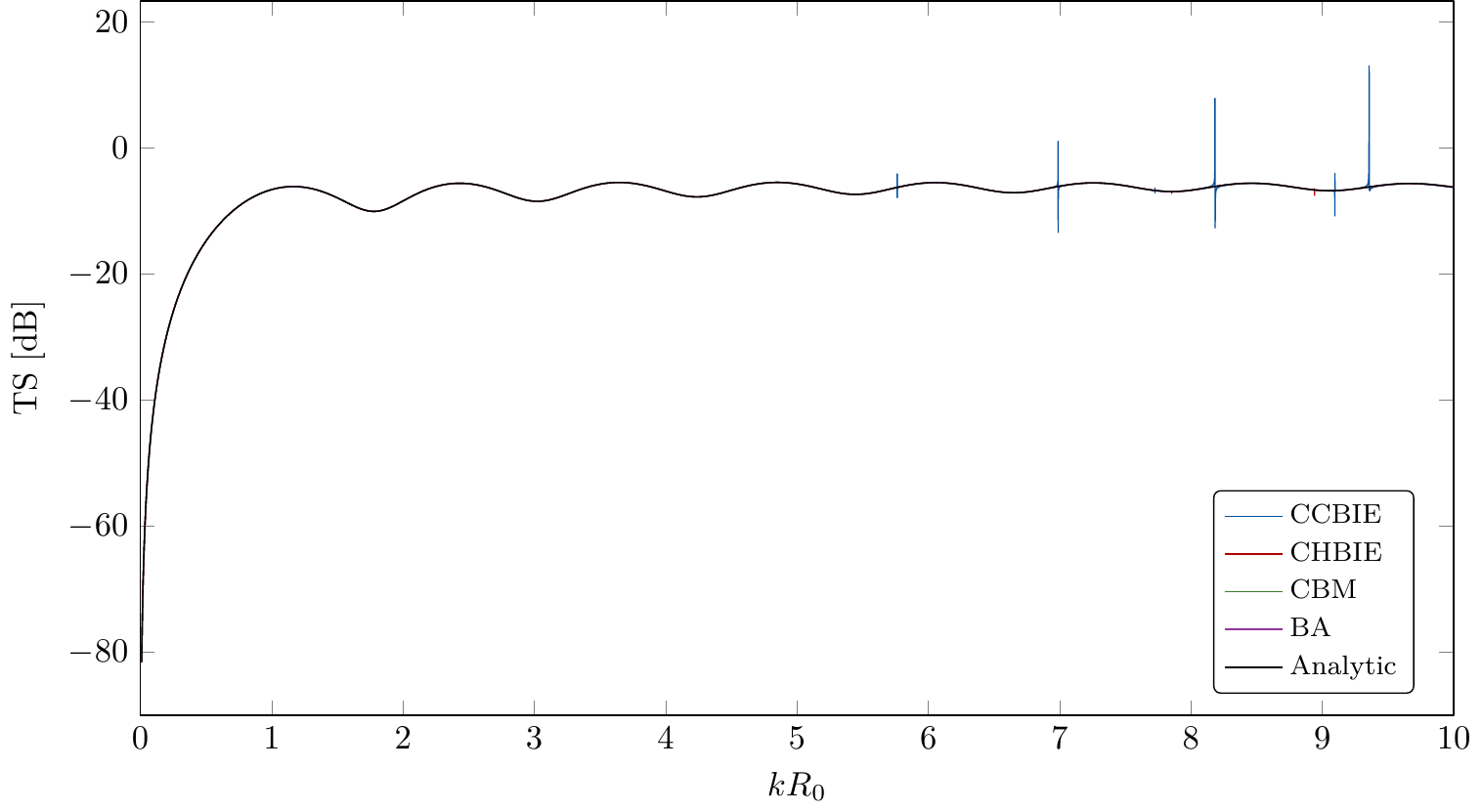}
		\caption{Target strength of backscattered far field.}
		\label{Fig3:eigenFreqDirichlet}
	\end{subfigure} 
	\par\bigskip
	\par\bigskip
	\begin{subfigure}[t]{\textwidth}
		\includegraphics[width=\textwidth]{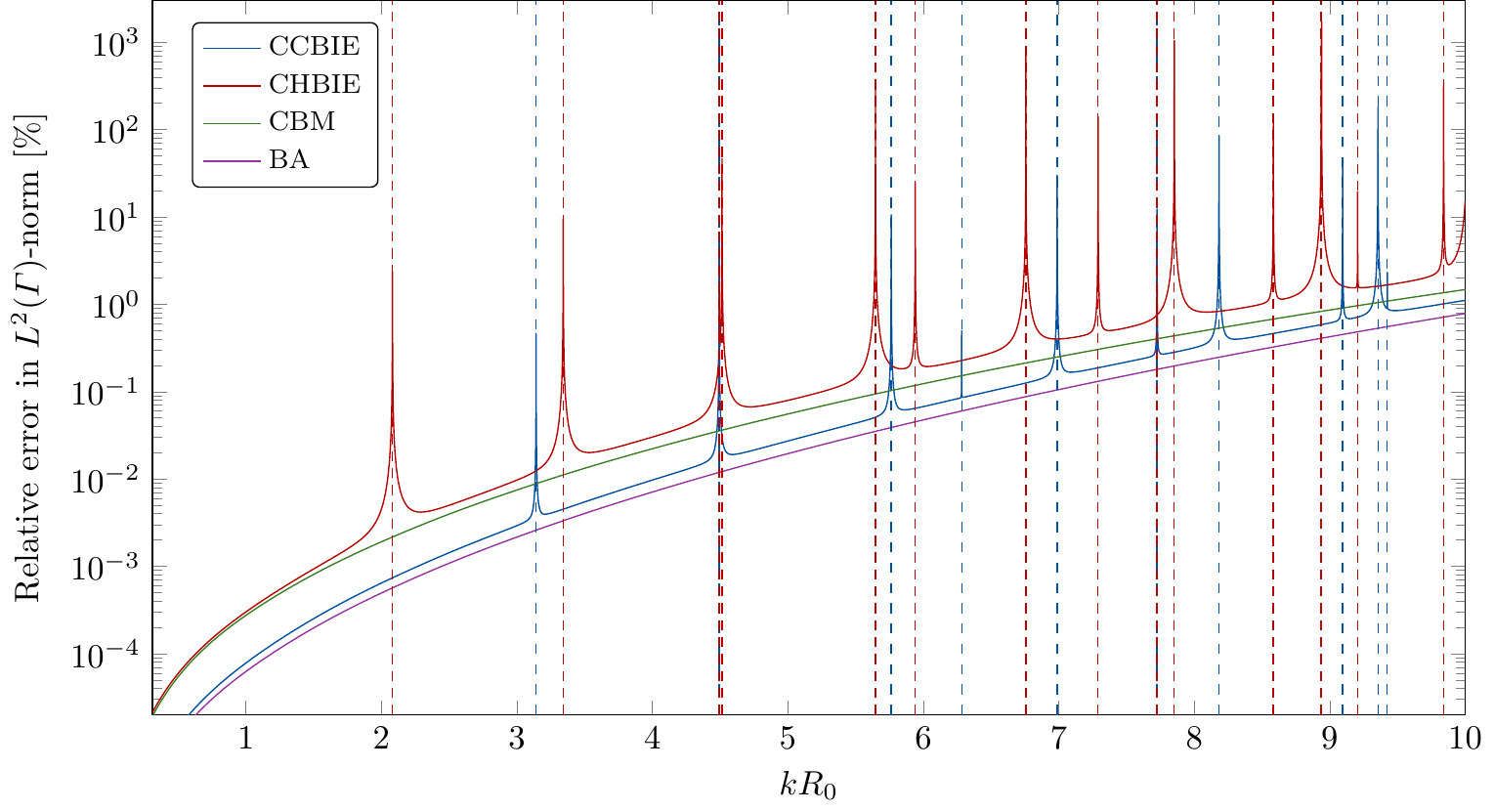}
		\caption{Surface error on $\Gamma$.}
		\label{Fig3:eigenFreqDirichletError}
	\end{subfigure} 
	\caption{\textbf{Rigid scattering on a sphere}: The plots show the instabilities around eigenfrequencies of the corresponding interior Dirichlet problem. All computations are done using the parametrization in \Cref{Fig3:parm2} refined uniformly three times with NURBS degree 4 (resulting in 384 elements and 728 degrees of freedom).}
\end{figure}

\begin{table}
	\centering
	\caption{The non-zero dimensionless eigenvalues below $kR_0=10$ for the interior Dirichlet problem~\cite{Zheng2015itb}.}
	\label{Tab3:eigenFreqDirichlet}
	\begin{tabular}{c l}
		\toprule
		$n$ & Roots of $\besselj_n(kR_0)$ \\
		\hline
		0 & $\PI$, $2\PI$, $3\PI$, \dots \\
	    1 & 4.49340945790907, 7.72525183693771, \dots \\
	  	2 & 5.76345919689455, 9.09501133047635, \dots \\
	    3 & 6.98793200050052, \dots \\
	    4 & 8.18256145257124, \dots \\
	    5 & 9.35581211104275, \dots \\
		\bottomrule
	\end{tabular}
\end{table}
\begin{table}
	\centering
	\caption{The non-zero dimensionless eigenvalues below $kR_0=10$ for the interior Neumann problem~\cite{Zheng2015itb}.}
	\label{Tab3:eigenFreqNeumann}
	\begin{tabular}{c l}
		\toprule
		$n$ & Roots of $\besselj_n'(kR_0)$\\
		\hline
		0 & 4.49340945790907, 7.72525183693771, \dots \\
	    1 & 2.08157597781810, 5.94036999057271, 9.20584014293667, \dots \\
	  	2 & 3.34209365736570, 7.28993230409335, \dots \\
	    3 & 4.51409964703228, 8.58375495636577, \dots \\
	    4 & 5.64670362043680, 9.84044604304014, \dots \\
	    5 & 6.75645633020413, \dots \\
	    6 & 7.85107767947440, \dots \\
	    7 & 8.93483887835284, \dots \\
		\bottomrule
	\end{tabular}
\end{table}

As we can see from \Cref{Fig3:eigenFreqDirichlet}, the dimensionless fictitious eigenfrequencies in \Cref{Tab3:eigenFreqDirichlet,Tab3:eigenFreqNeumann} appear quite clearly for the CBIE and the HBIE, respectively, while the eigenvalues for the Burton--Miller formulation are shifted away from the real axis into the complex plane \cite{Zheng2015itb}. The fictitious eigenfrequencies are of course not present in the best approximation (BA) solution.


\subsection{Torus interior acoustic problem}
\begin{figure}
	\centering
	\includegraphics[width=\textwidth]{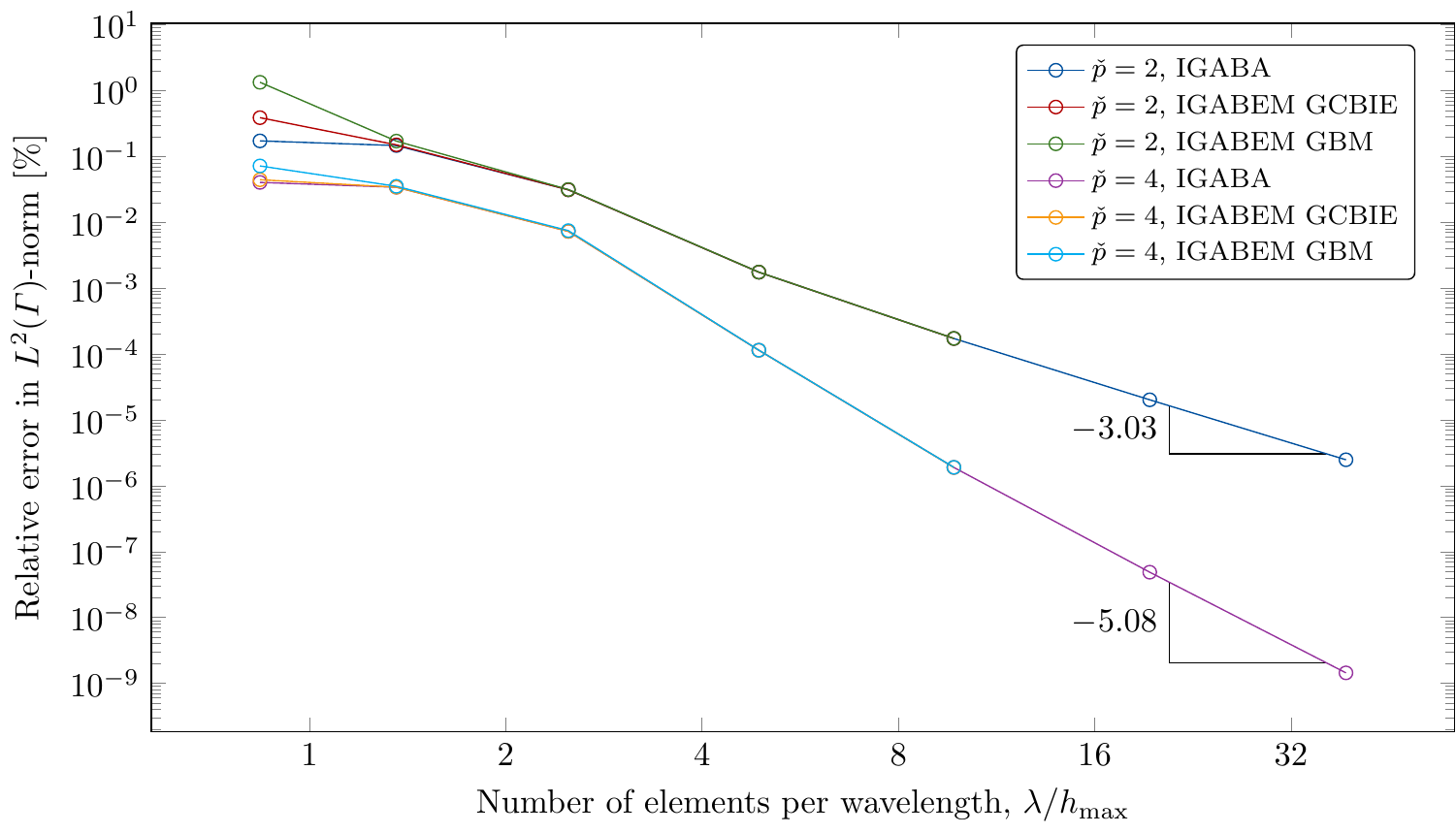}
	\caption{\textbf{Torus interior acoustic problem}: Convergence analysis at $k=\SI{2}{m^{-1}}$.}
	\label{Fig3:torus}
\end{figure}

Consider the Torus problem presented in~\cite{Simpson2014aib}. This example sets the stage for optimal conditions for the a priori error estimate in~\Cref{Eq3:aprioriErrorEstimate} to be fulfilled. The geometry of the torus (with parametrization described in~\Cref{Sec3:torus}) has $G^\infty$ continuity and contains no polar singularities in the exact NURBS parametrization illustrated in~\Cref{Fig3:Torus} (as opposed to the sphere parametrization in~\Cref{Fig3:parm1}). The torus considered here has major radius $r_{\mathrm{o}} = 2$ and minor radius $r_{\mathrm{i}}=1$. Consider the exact solution
\begin{equation*}
	p(\vec{x}) = \sin\frac{kx_1}{\sqrt{3}}\sin\frac{kx_2}{\sqrt{3}}\sin\frac{kx_3}{\sqrt{3}}
\end{equation*}
with corresponding Neumann boundary conditions at the boundary $\Gamma$
\begin{equation*}
	\pderiv{p}{n} = \frac{k}{\sqrt{3}}\begin{bmatrix}
	\cos\frac{kx_1}{\sqrt{3}}\sin\frac{kx_2}{\sqrt{3}}\sin\frac{kx_3}{\sqrt{3}}\\
	\sin\frac{kx_1}{\sqrt{3}}\cos\frac{kx_2}{\sqrt{3}}\sin\frac{kx_3}{\sqrt{3}}\\
	\sin\frac{kx_1}{\sqrt{3}}\sin\frac{kx_2}{\sqrt{3}}\cos\frac{kx_3}{\sqrt{3}}
	\end{bmatrix}\cdot\vec{n}.
\end{equation*}
From \Cref{Fig3:torus}, the sharpness ($C_1 \approx 1$) of the a priori error estimate in \Cref{Eq3:aprioriErrorEstimate} is demonstrated. The convergence rates for the best approximation (IGABA) are revealed quite clearly here. 

\begin{figure}
	\centering
	\includegraphics[width=\textwidth]{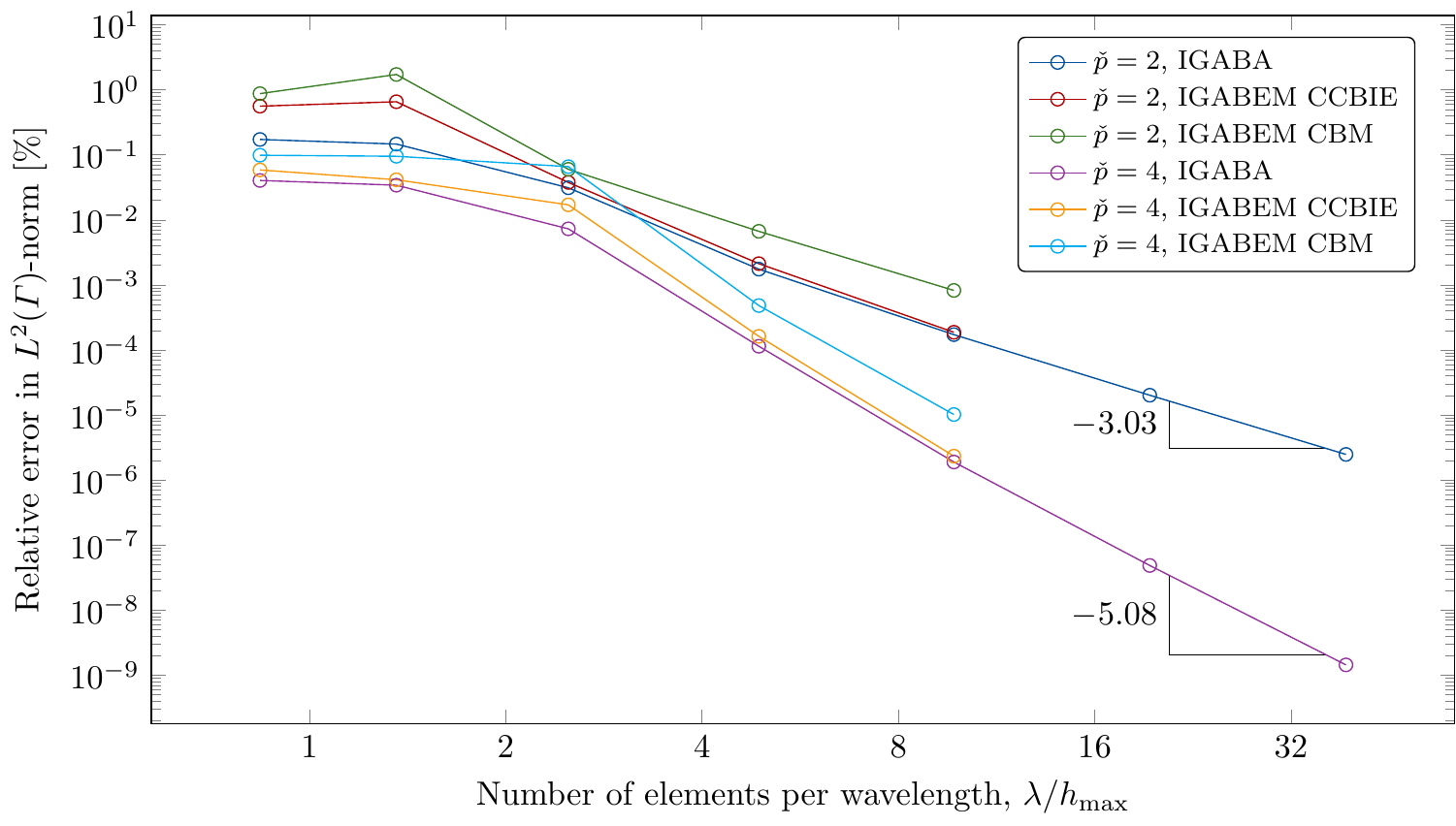}
	\caption{\textbf{Torus interior acoustic problem}: Convergence analysis at $k=\SI{2}{m^{-1}}$.}
	\label{Fig3:torus_2}
\end{figure}

Results for the same study using collocation formulation are given in \Cref{Fig3:torus_2}. The CCBIE formulation obtains very good results as it approaches the best approximation during refinement. Correct convergence rates are also obtained for the CBM formulation, but with a somewhat higher constant $C_1$ in \Cref{Eq3:aprioriErrorEstimate}.

In~\cite{Simpson2014aib} Simpson projects the Neumann data onto the same basis used for the solution space. The accuracy for collocation formulations may be increased in some cases using this projection, but for Galerkin formulations projecting the Neumann data yields worse results. Moreover, if $\Gamma$ is $G^0$ sub optimal results are obtained also for the collocation formulations.

\subsection{Manufactured solutions for complex geometries}
In this section we shall consider the method of manufactured solutions (MMS). The idea behind MMS is explained in detail in \cite{Roy2005roc}. 

By construction of the fundamental solution of the Helmholtz equation ($\Phi_k(\vec{x},\vec{y})$ in \Cref{Eq3:FreeSpaceGrensFunction}), the function $p(\vec{x}) = \Phi_k(\vec{x},\vec{y})$ is a solution to \Cref{Eq3:HelmholtzEqn,Eq3:HelmholtzEqnNeumannCond,Eq3:sommerfeldCond} whenever $\vec{y}\in\R^3\setminus\overline{\Omega^+}$ and for the Neumann boundary condition $g(\vec{x})=\partial_n\Phi_k(\vec{x},\vec{y})$ on $\Gamma$.  Hence, we have an exact manufactured solution for the exterior Helmholtz problem for arbitrary geometries $\Gamma$ which encloses the point $\vec{y}$. It is emphasized that this solution is non-physical for non-spherical geometries $\Gamma$ (for the sphere, the solution represents a pulsating sphere \cite{Simpson2014aib}). General solutions may be constructed by separation of variables (cf.~\cite[p. 26]{Ihlenburg1998fea})
\begin{equation}
	p(\vec{x}) = \sum_{n=0}^\infty\sum_{m=-n}^n C_{nm} \hankel_n^{(1)}(kR) \legendre_n^{|m|}(\cos\vartheta)\euler^{\imag m\varphi} 
\end{equation}
with
\begin{equation*}
	R = |\vec{x}-\vec{y}|,\quad \vartheta=\arccos\left(\frac{x_3-y_3}{R}\right),\quad\varphi = \operatorname{atan2}(x_2-y_2,x_1-y_1)
\end{equation*}
where $\hankel_n^{(1)}$ is the $n^{\mathrm{th}}$ spherical Hankel function of first kind and $\legendre_n^m$ are the associated Legendre functions. In fact, the solution $p(\vec{x}) = \Phi_k(\vec{x},\vec{y})$ is a special case of this general form with 
\begin{equation}
	C_{nm} = \begin{cases}
		\frac{\imag k}{4\PI} & n = 0,\,\,m=0\\
		0 & \text{otherwise}.
		\end{cases}
\end{equation}
Inspired by the method of fundamental solutions~\cite{Fairweather2003tmo}, we can also use the solution
\begin{equation}\label{Eq3:manuMFS}
	p(\vec{x}) = \sum_{n=1}^N C_n \Phi_k(\vec{x},\vec{y}_n)
\end{equation}
for a set of $N$ source points $\{\vec{y}_n\}_{n=1}^N$. To increase the complexity of the solution, we use $C_n=\cos(n-1)$ in this work.


The complexity of this problem setup does not scale with the complexity of the model as it is independent of $\Gamma$. However, it preserves two important properties of acoustic scattering, namely the radial decay and the oscillatory nature. Thus, this problem setup represents a general way of constructing manufactured solutions that can be utilized to verify the correctness of the implemented code for solving the Helmholtz equation. Moreover, as the boundary condition is the only condition that is altered from the original problem, one can solve the original system of equation with an extra appended column vector on the right-hand side (corresponding to the problem of finding the manufactured solution) with a small computational effort. This gives some control over the correctness of the computed solution to the original problem. Since the fictitious eigenfrequencies are the same for both solutions, one can compute the error for the manufactured solution to give an indication whether the solution is polluted by such a frequency. If this is the case, one should resort to the somewhat more costly Burton--Miller formulation.

Note that from the first limit of \Cref{Eq3:Phi_k_limits}, the far field of \Cref{Eq3:manuMFS} is given by
\begin{equation*}
	p_0(\hat{\vec{x}})=\frac{1}{4\PI}\sum_{n=1}^N C_n \euler^{-\imag k \hat{\vec{x}}\cdot\vec{y}_n}.
\end{equation*}

Whenever $\partial_n p_{\mathrm{tot}}\neq 0$ we must deal with an integral which is weakly singular, and the manufactured solution thus does not give the optimal test for the rigid body scattering problem as the CBIE formulation is free from such integrals in this case.

\subsubsection{Manufactured solution with a cube}
\begin{figure}
	\centering
	\begin{subfigure}[t]{0.3\textwidth}
		\includegraphics[width=0.9\textwidth]{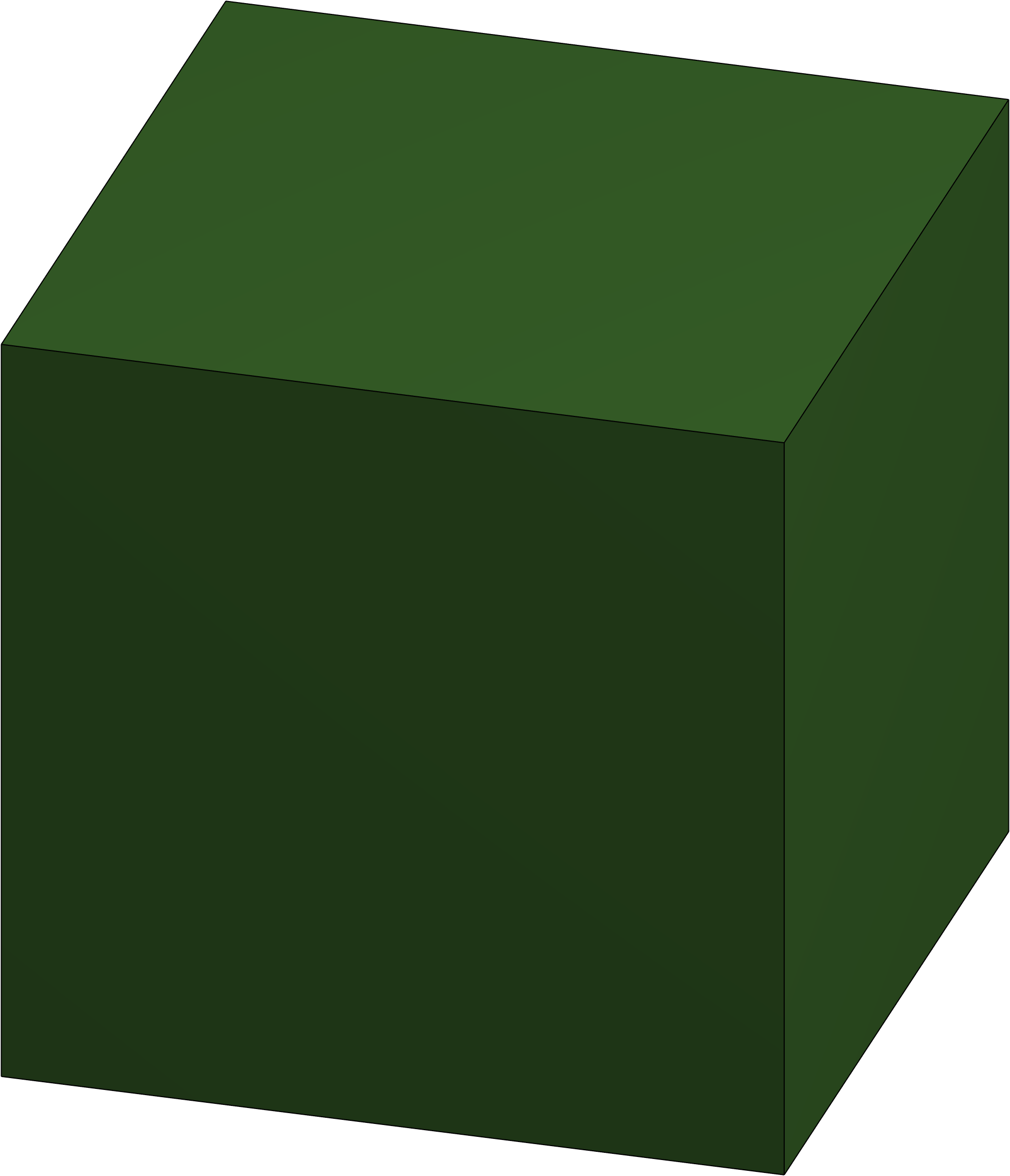}
		\caption{Mesh 1.}
		\label{Fig3:Cube_mesh1}
	\end{subfigure} 
	~
	\begin{subfigure}[t]{0.3\textwidth}
		\includegraphics[width=0.9\textwidth]{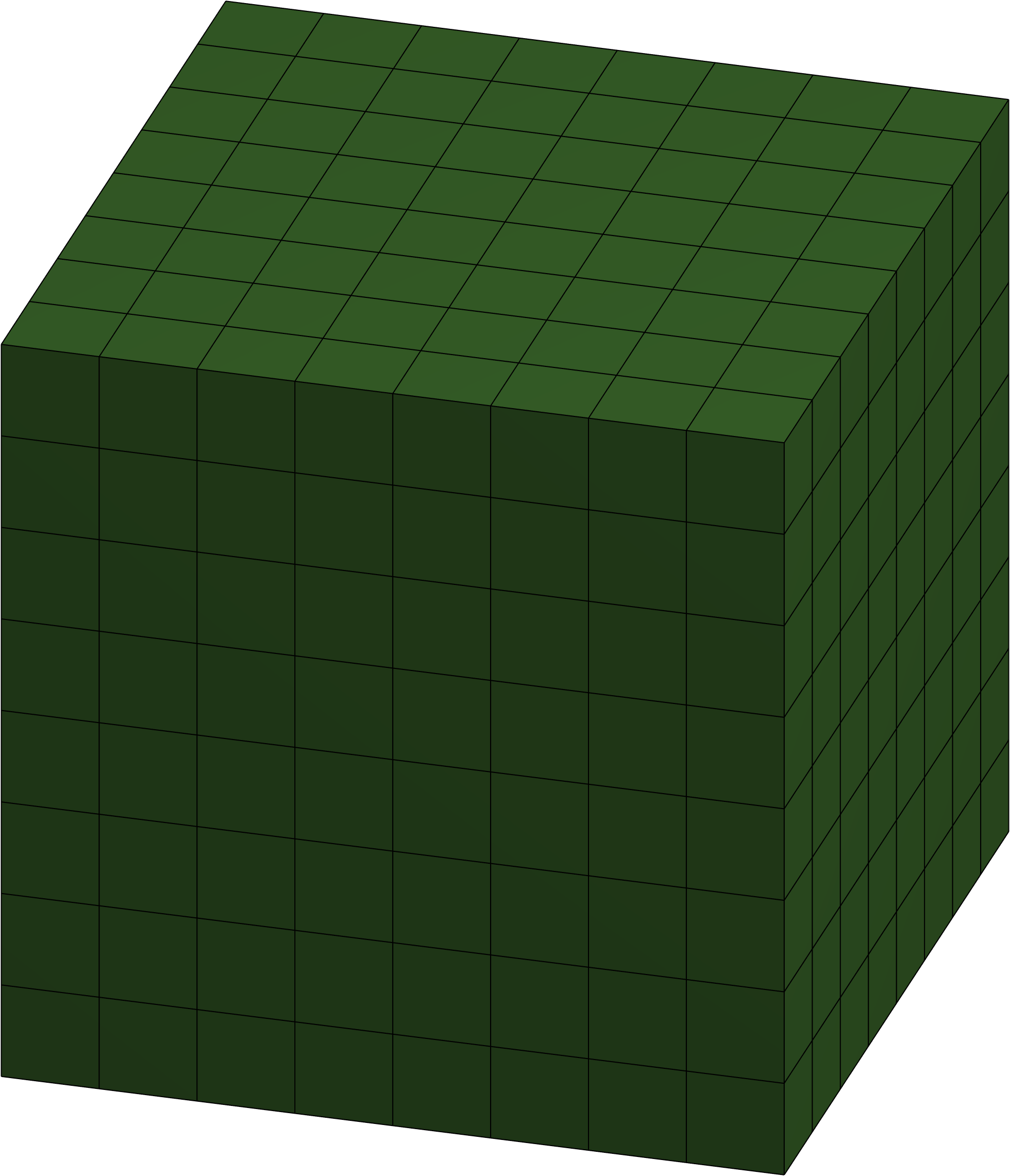}
		\caption{Mesh 4.}
		\label{Fig3:Cube_mesh4}
	\end{subfigure} 
	\caption{Parametrization of a cube using 6 patches of degree $\check{p}\geq 1$.}
	\label{Fig3:CubeParametrizations}
\end{figure}
Consider a cube of side length $a$ centered at the origin. Its interior Dirichlet problem has eigenfunctions (cf. \cite[p. 52]{Schenck1968iif})
\begin{equation*}
	p(\vec{x}) = \prod_{i=1}^d\sin\frac{n_i\PI (x_i+a/2)}{a},\quad \vec{x}\in\Omega^-
\end{equation*}
and the interior Neumann problem has eigenfunctions
\begin{equation*}
	p(\vec{x}) = \prod_{i=1}^d\cos\frac{n_i\PI (x_i+a/2)}{a},\quad \vec{x}\in\Omega^-
\end{equation*}
where
\begin{equation*}
	\sum_{i=1}^d n_i^2 = \left(\frac{ka}{\PI}\right)^2\quad\text{and}\quad\Omega^-=\left[-\frac{a}{2},\frac{a}{2}\right]^d.
\end{equation*}
The dimensionless eigenfrequencies are thus given by
\begin{equation*}
	ka=\PI\sqrt{\sum_{i=1}^d n_i^2},
\end{equation*}
where $n_i\in\N^*$ for the interior Dirichlet problem and $n_i\in\N$ for the interior Neumann problem. For the exterior problem these eigenfrequencies correspond to the fictitious eigenfrequencies for the CBIE formulation and the HBIE formulation, respectively. The dimensionless fictitious eigenfrequencies below $ka=10$ are then $\PI\sqrt{3}$, $\PI\sqrt{6}$ and $3\PI$ for the CBIE formulation, and $\PI\sqrt{n}$ with $n\in\{0,1,2,3,4,5,6,8,9,10\}$ for the HBIE formulation.

\begin{figure}
	\centering
	\includegraphics[width=\textwidth]{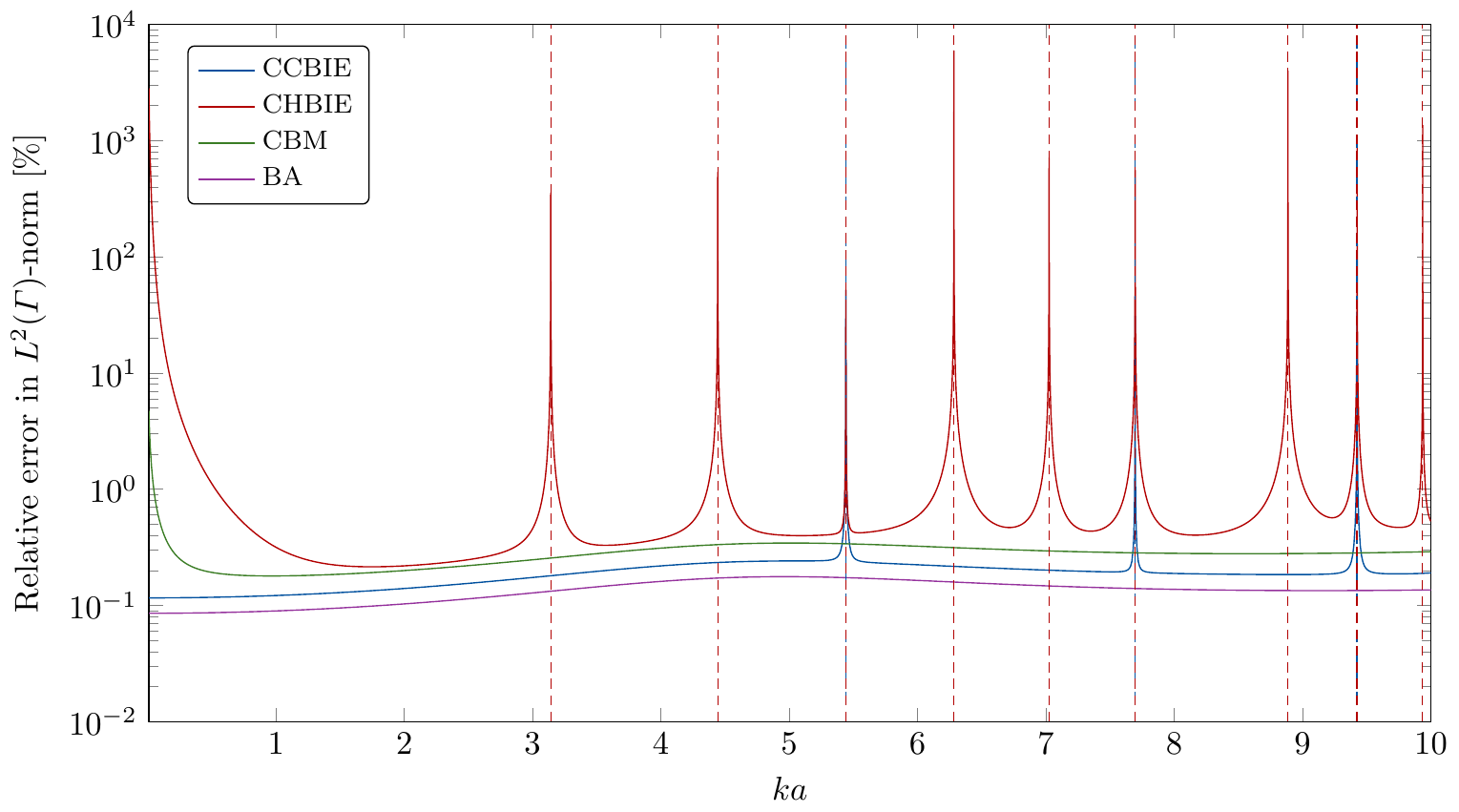}
	\caption{\textbf{Manufactured solution with a cube}: The plots show the instabilities around eigenfrequencies of the corresponding interior Dirichlet problem. All computations are done using the parametrization in \Cref{Fig3:Cube_mesh1} refined uniformly three times with NURBS degree 4 (resulting in 384 elements and 728 degrees of freedom) as highlighted in \Cref{Fig3:Cube_mesh4}.}
	\label{Fig3:Cube_P_sweep}
\end{figure}

\begin{figure}
	\centering
	\begin{subfigure}[t]{\textwidth}
		\includegraphics[width=\textwidth]{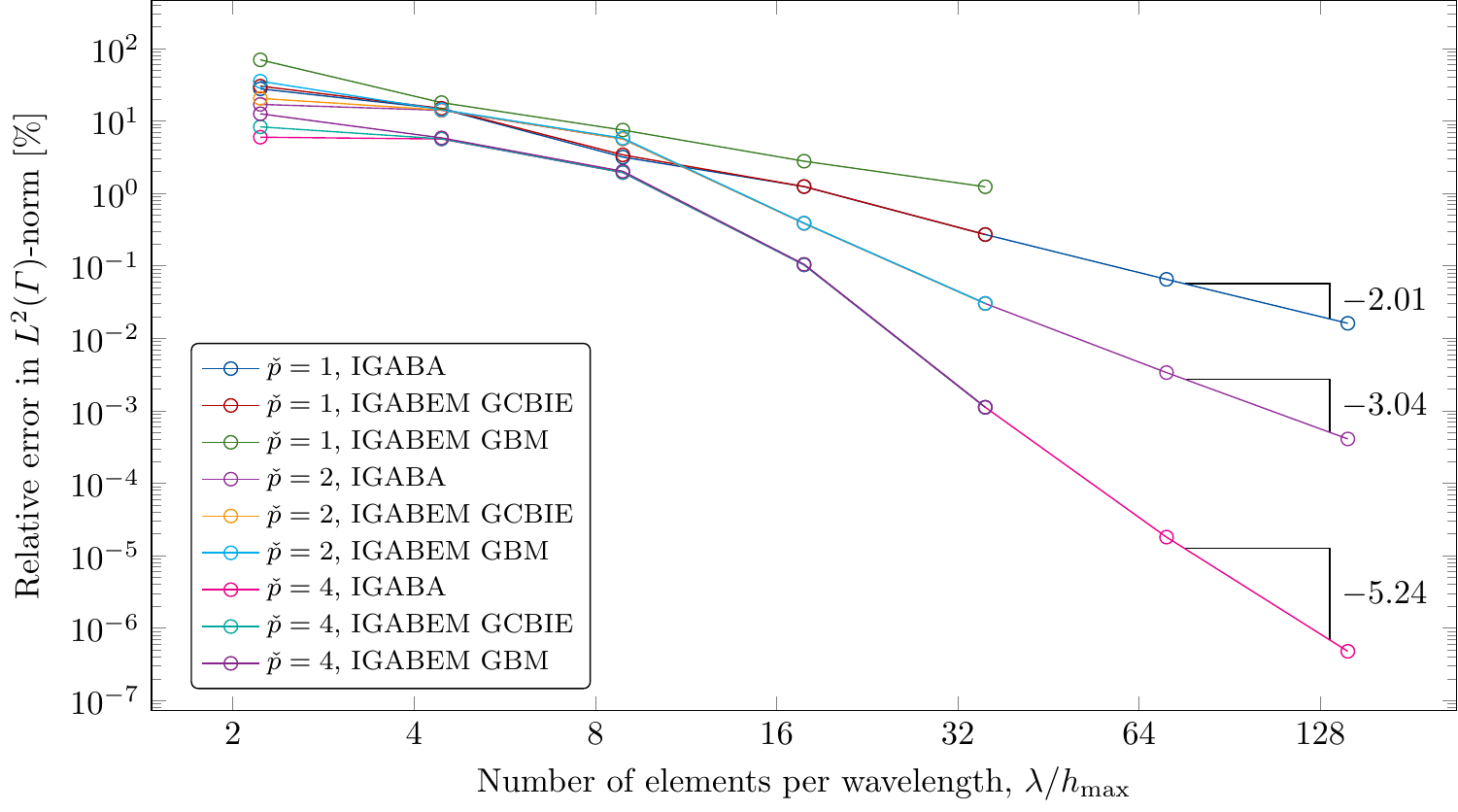}
		\caption{Galerkin formulations}
		\label{Fig3:Cube_P1}
	\end{subfigure} 
	\par\bigskip
	\par\bigskip
	\begin{subfigure}[t]{\textwidth}
		\includegraphics[width=\textwidth]{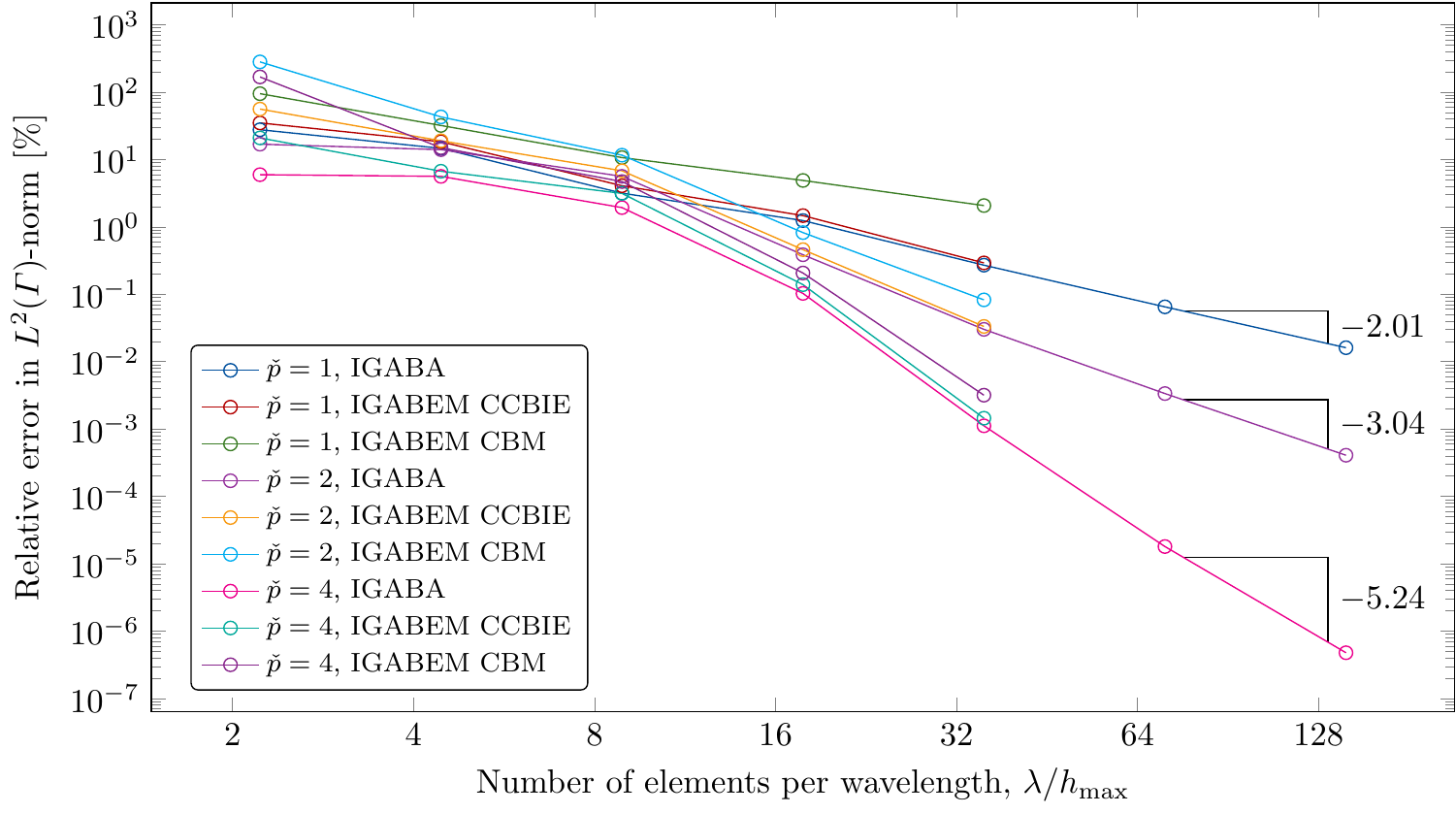}
		\caption{Collocation formulations}
		\label{Fig3:Cube_P2}
	\end{subfigure} 
	\caption{\textbf{Manufactured solution with a cube}: Convergence analysis at $k=\SI{2}{m^{-1}}$.}
	\label{Fig3:Cube_P}
\end{figure}

Consider the manufactured solution~\Cref{Eq3:manuMFS} with $N=3^3=27$ source points 
\begin{equation*}
	\vec{y}_n = \frac{a}{4}[c_i, c_j, c_l],\quad n=i+3(j-1)+3^2(l-1),\quad i,j,l=1,2,3
\end{equation*}
where $c_1=-1$, $c_2=0$ and $c_3=1$. In~\Cref{Fig3:Cube_P_sweep} we again show a frequency sweep to illustrate the instability around the fictitious eigenfrequencies of the CBIE and HBIE formulations.
From \Cref{Fig3:Cube_P1}, the sharpness of the a priori error estimate in \Cref{Eq3:aprioriErrorEstimate} is again demonstrated. Remarkably, the $G^0$ continuity of the cube poses no problems for the Galerkin Burton--Miller formulation using $\check{p}\geq 2$, despite the problematic mathematical nature of the formulations with basis functions that are $C^0$ continuous~\cite{Liu1999anf}. Poor results are obtained for the BM formulation using $\check{p}=1$ for both collocation and Galerkin formulation. This is in stark contrast to the CBIE which performs optimally for $\check{p}=1$ in both cases. The CCBIE obtains good results in all cases and outperforms the CBM formulation.

\subsubsection{Manufactured solution with the BeTSSi submarine}
\label{Subsec3:manufactured}
\begin{figure}
	\centering
	\begin{subfigure}[t]{\textwidth}
		\includegraphics[width=\textwidth]{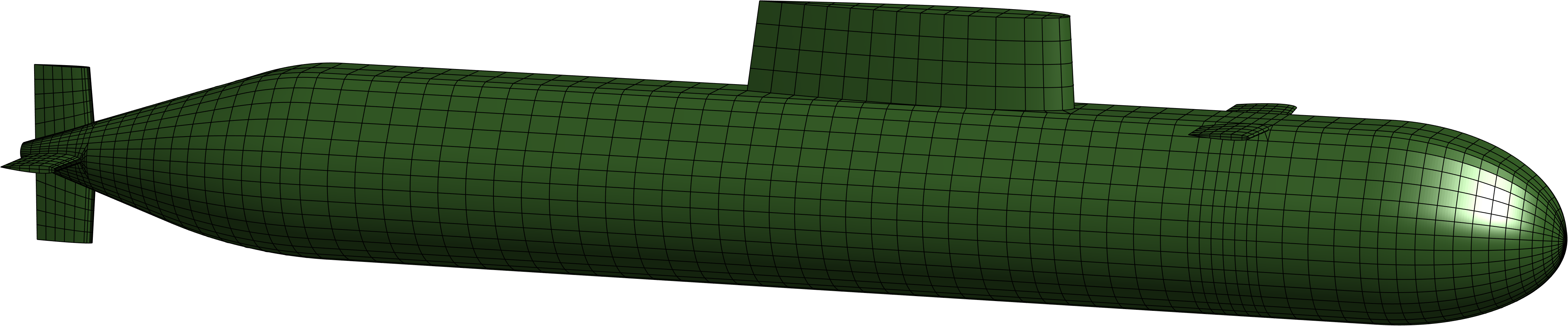}
		\caption{${\cal M}_{1,6,5}^{\textsc{igabem}}$ - 3718 elements}
	\end{subfigure} 
	\par\bigskip
	\par\bigskip
	\begin{subfigure}[t]{\textwidth}
		\includegraphics[width=\textwidth]{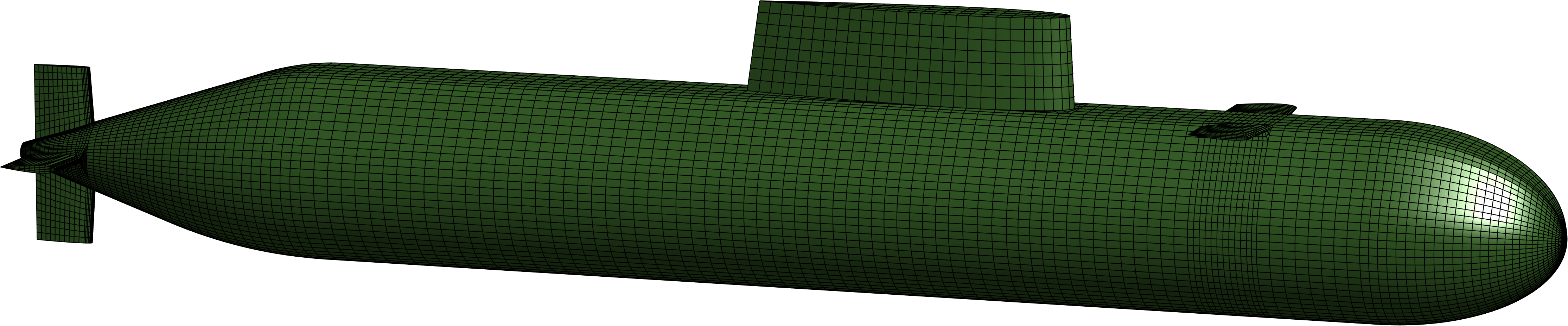}
		\caption{${\cal M}_{2,6,5}^{\textsc{igabem}}$ - 14872 elements}
	\end{subfigure} 
	\par\bigskip
	\par\bigskip
	\begin{subfigure}[t]{\textwidth}
		\includegraphics[width=\textwidth]{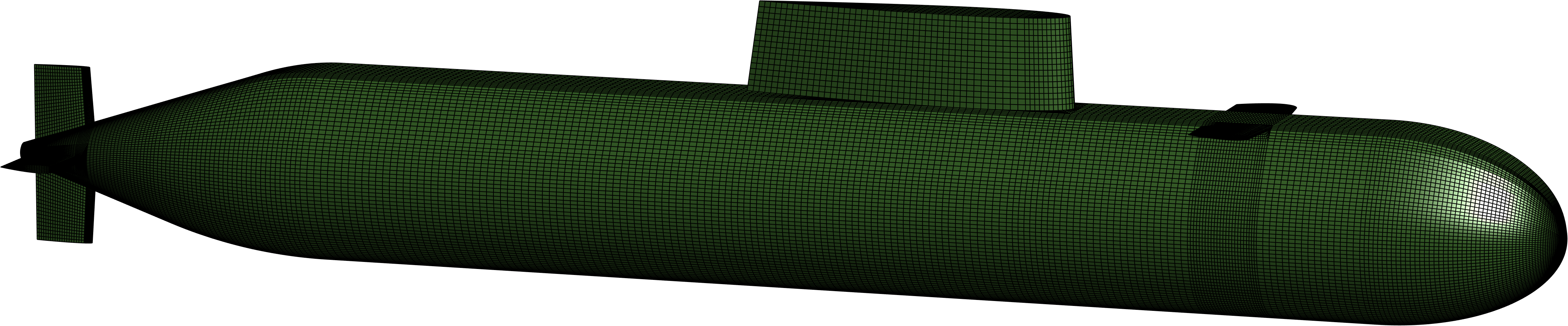}
		\caption{${\cal M}_{3,6,5}^{\textsc{igabem}}$ - 59488 elements}
	\end{subfigure} 
	\caption{\textbf{The BeTSSi submarine}: Computational IGA meshes for $\Gamma_{\check{p}}$ with $\check{p}=6$.}
	\label{Fig3:BeTSSimeshes}
\end{figure}
Consider now the BeTSSi submarine described in \Cref{Sec3:BeTSSi_description}. The BeTSSi meshes considered in this work are denoted by ${\cal M}_{m,\check{p},\check{k}}^{\textsc{igabem}}$, where $m$ is the mesh number, and are illustrated in \Cref{Fig3:BeTSSimeshes} where $m=1$ is the coarsest mesh, and $m=2$ and $m=3$ are uniformly refined meshes iterated on the coarsest mesh. Again, $\check{p}$ denotes the polynomial order and $\check{k}$ the continuity.

\begin{figure}
	\centering
	\begin{subfigure}[t]{\textwidth}
		\includegraphics[width=\textwidth]{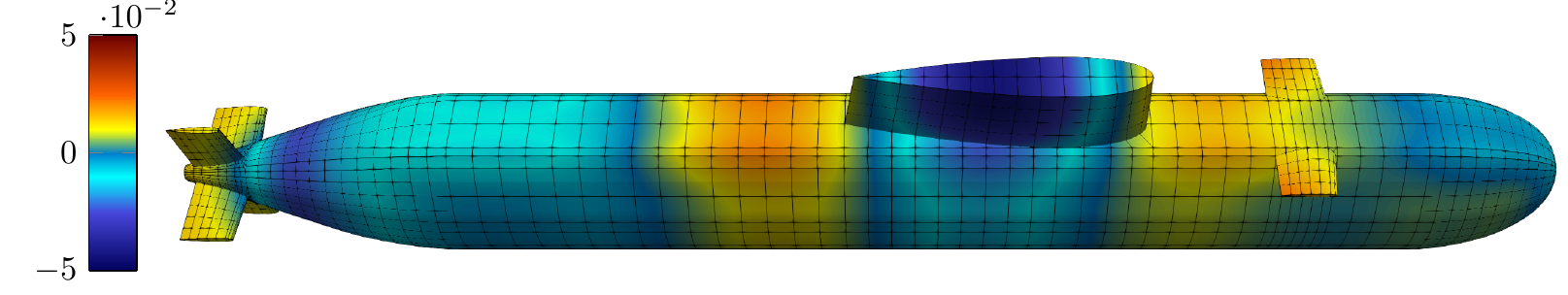}
		\caption{$f=\SI{100}{Hz}$}
	\end{subfigure} 
	\par\bigskip
	\begin{subfigure}[t]{\textwidth}
		\includegraphics[width=\textwidth]{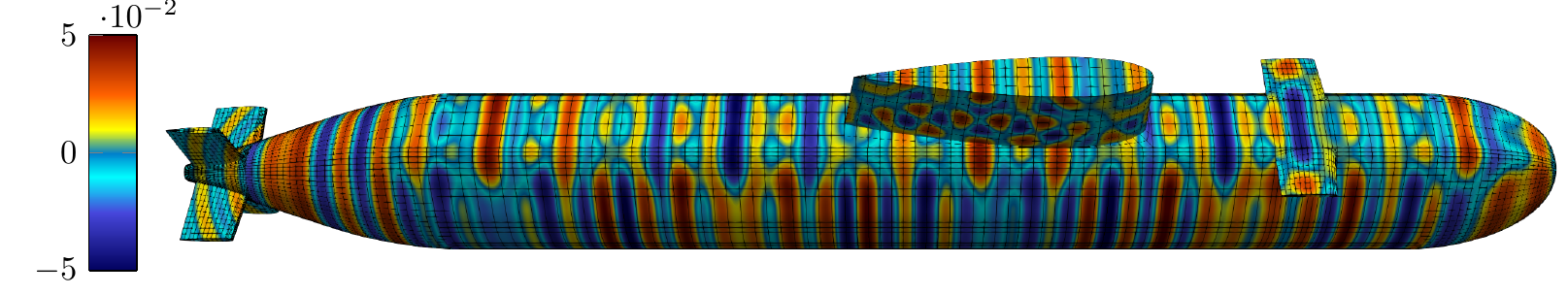}
		\caption{$f=\SI{1000}{Hz}$}
	\end{subfigure} 
	\caption{\textbf{Manufactured solution with the BeTSSi submarine}: Analytic manufactured solution. Mesh 1 and mesh 2 are added to visualize elements to wavelength ratio for $\SI{100}{Hz}$ and $\SI{1000}{Hz}$, respectively.}
	\label{Fig3:BCA_P_analytic}
\end{figure}

\begin{figure}
	\centering
	\begin{subfigure}[t]{\textwidth}
		\includegraphics[width=\textwidth]{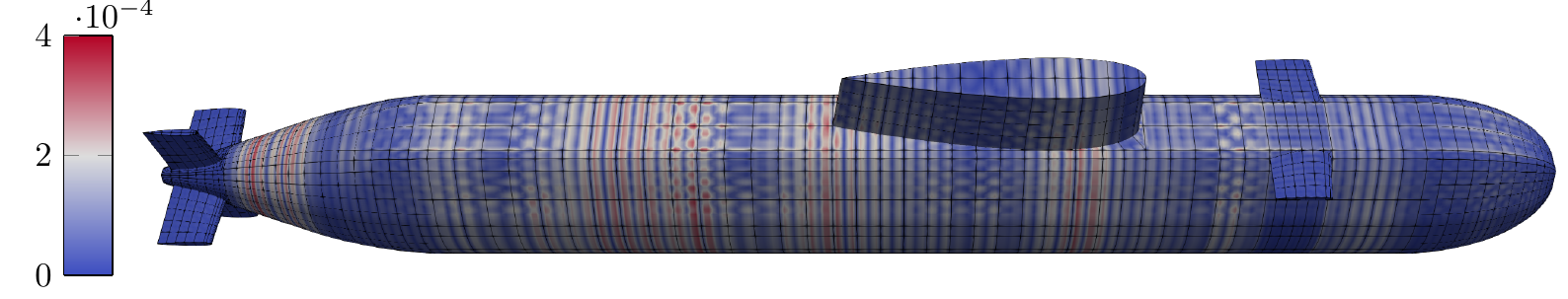}
		\caption{IGABA: $\displaystyle \max_{\vec{x}\in\Gamma_2}\frac{|p-p_h|}{|p|} = \num{3.8e-4}$}
	\end{subfigure} 
	\par\bigskip
	\begin{subfigure}[t]{\textwidth}
		\includegraphics[width=\textwidth]{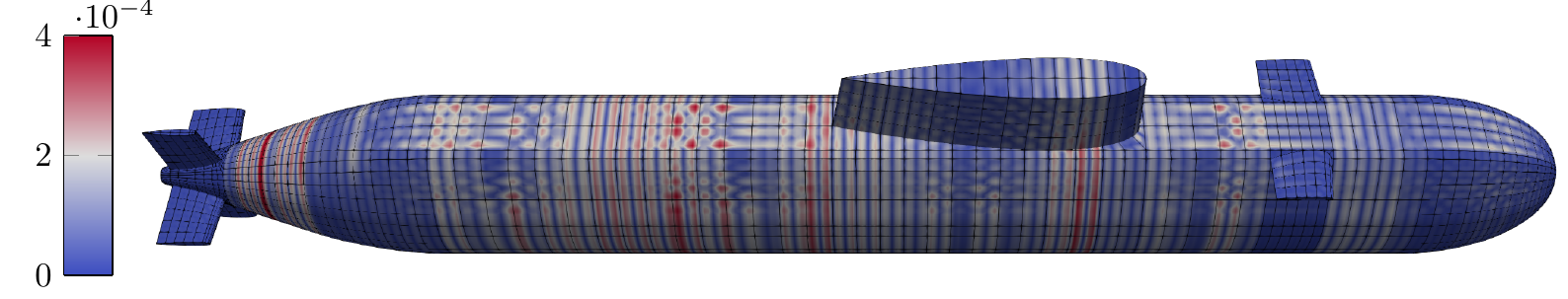}
		\caption{IGABEM CCBIE: $\displaystyle \max_{\vec{x}\in\Gamma_2}\frac{|p-p_h|}{|p|} = \num{4.0e-4}$}
	\end{subfigure} 
	\caption{\textbf{Manufactured solution with the BeTSSi submarine}: Relative error on the surface of the scatterer at $f=\SI{100}{Hz}$ on the mesh ${\cal M}_{1,2,1}^{\textsc{igabem}}$.}
	\label{Fig3:BCA_P_p2M1f100}
\end{figure}

\begin{figure}
	\centering
	\begin{subfigure}[t]{\textwidth}
		\includegraphics[width=\textwidth]{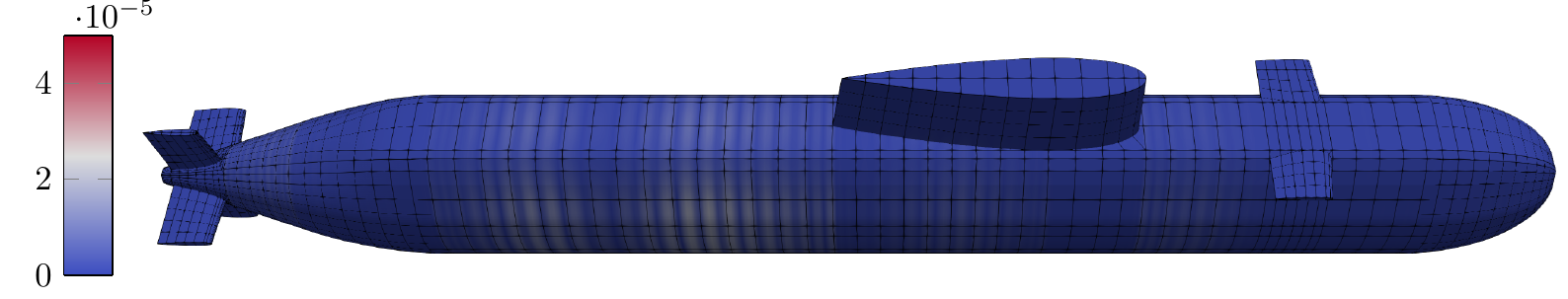}
		\caption{IGABA: $\displaystyle \max_{\vec{x}\in\Gamma_5}\frac{|p-p_h|}{|p|} = \num{2.1e-5}$}
	\end{subfigure} 
	\par\bigskip
	\begin{subfigure}[t]{\textwidth}
		\includegraphics[width=\textwidth]{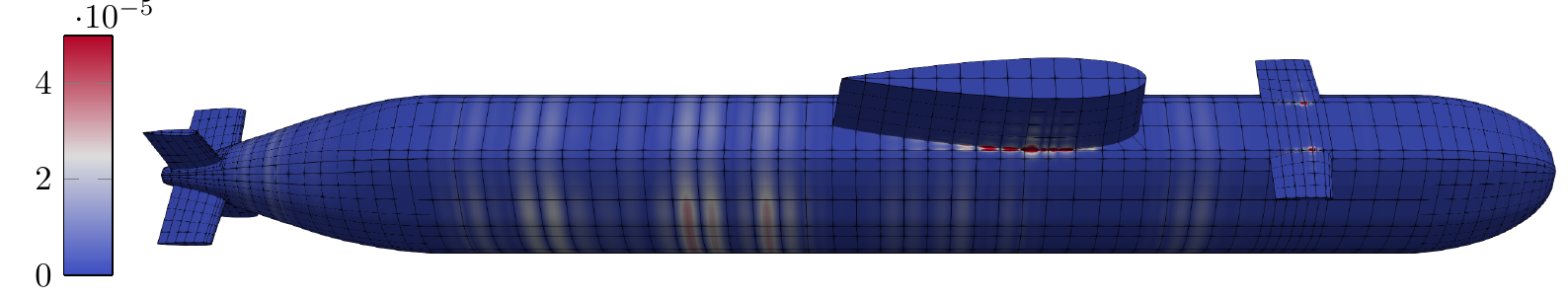}
		\caption{IGABEM CCBIE: $\displaystyle \max_{\vec{x}\in\Gamma_5}\frac{|p-p_h|}{|p|} = \num{4.4e-4}$}
	\end{subfigure} 
	\caption{\textbf{Manufactured solution with the BeTSSi submarine}: Relative error on the surface of the scatterer at $f=\SI{100}{Hz}$ on the mesh ${\cal M}_{1,5,4}^{\textsc{igabem}}$.}
	\label{Fig3:BCA_P_p5M1f100}
\end{figure}

\begin{figure}
	\centering
	\begin{subfigure}[t]{\textwidth}
		\includegraphics[width=\textwidth]{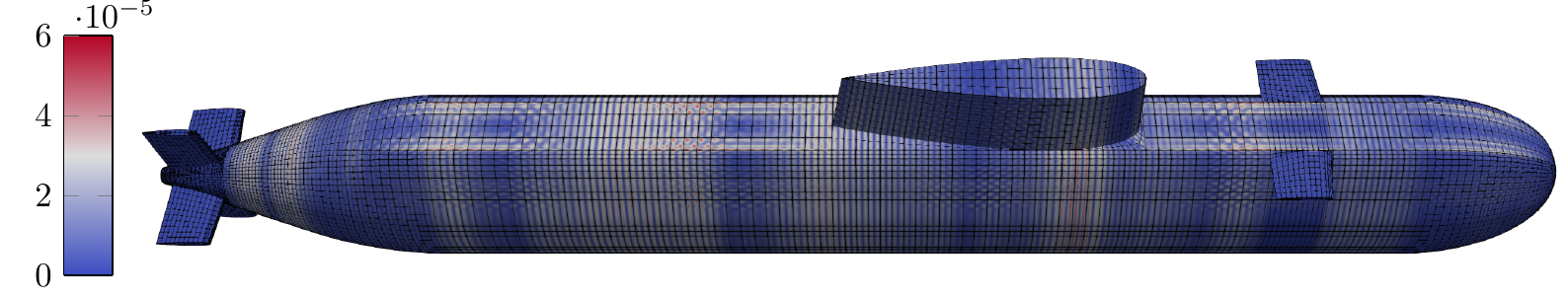}
		\caption{IGABA: $\displaystyle \max_{\vec{x}\in\Gamma_2}\frac{|p-p_h|}{|p|} = \num{6.1e-5}$}
	\end{subfigure} 
	\par\bigskip
	\begin{subfigure}[t]{\textwidth}
		\includegraphics[width=\textwidth]{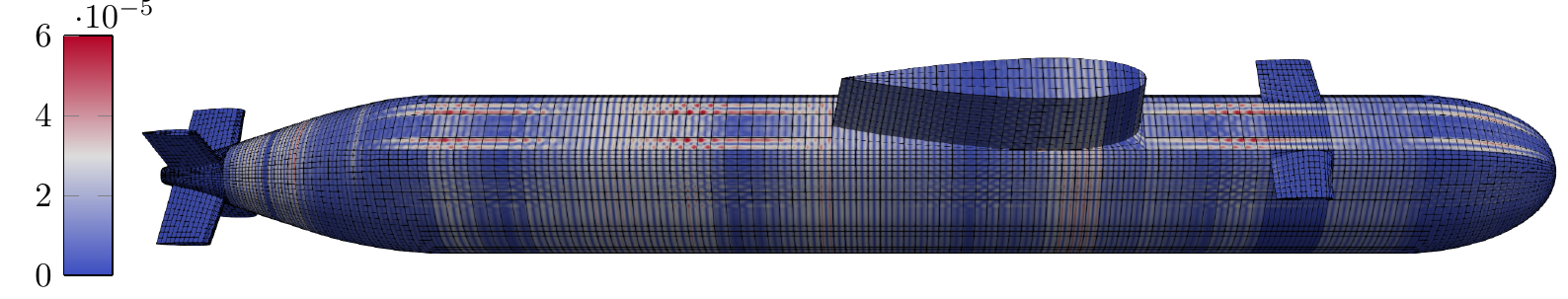}
		\caption{IGABEM CCBIE: $\displaystyle \max_{\vec{x}\in\Gamma_2}\frac{|p-p_h|}{|p|} = \num{2.3e-4}$}
	\end{subfigure} 
	\caption{\textbf{Manufactured solution with the BeTSSi submarine}: Relative error on the surface of the scatterer at $f=\SI{100}{Hz}$ on the mesh ${\cal M}_{2,2,1}^{\textsc{igabem}}$.}
	\label{Fig3:BCA_P_p2M2f100}
\end{figure}

\begin{figure}
	\centering
	\begin{subfigure}[t]{\textwidth}
		\includegraphics[width=\textwidth]{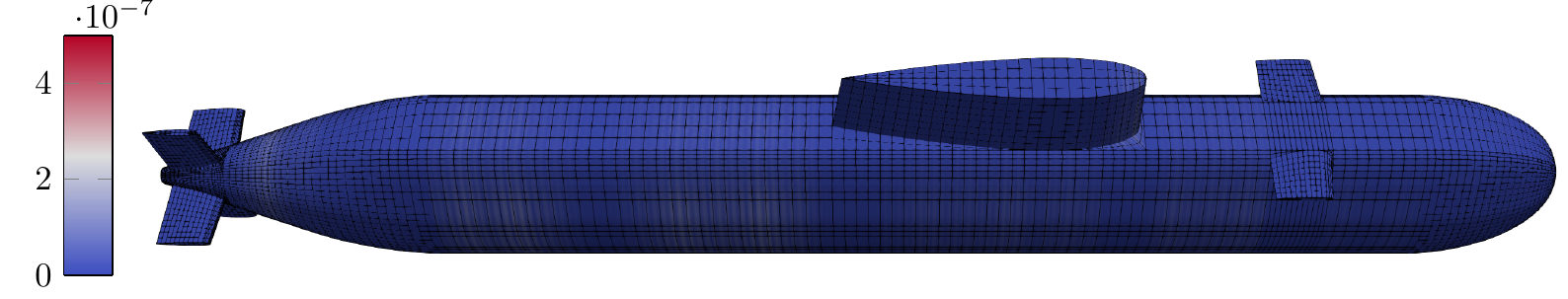}
		\caption{IGABA: $\displaystyle \max_{\vec{x}\in\Gamma_5}\frac{|p-p_h|}{|p|} = \num{1.7e-7}$}
	\end{subfigure} 
	\par\bigskip
	\begin{subfigure}[t]{\textwidth}
		\includegraphics[width=\textwidth]{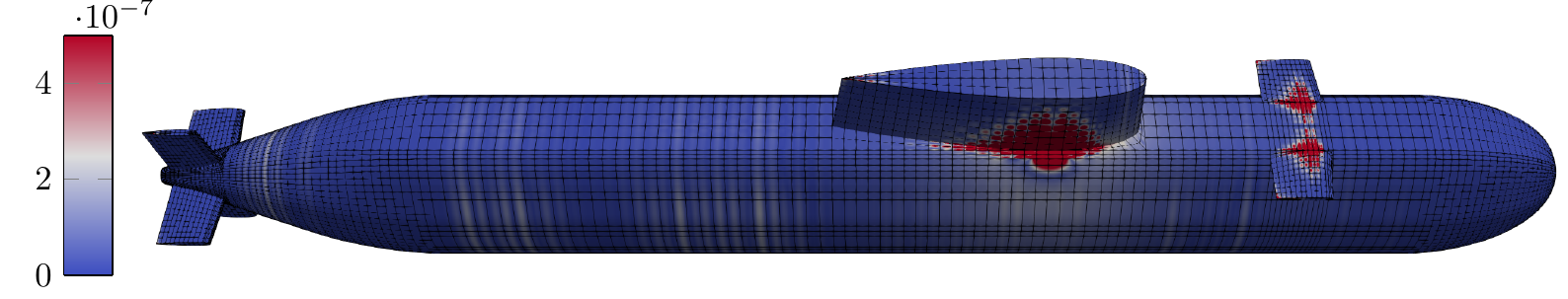}
		\caption{IGABEM CCBIE: $\displaystyle \max_{\vec{x}\in\Gamma_5}\frac{|p-p_h|}{|p|} = \num{7.8e-4}$}
	\end{subfigure} 
	\caption{\textbf{Manufactured solution with the BeTSSi submarine}: Relative error on the surface of the scatterer at $f=\SI{100}{Hz}$ on the mesh ${\cal M}_{2,5,4}^{\textsc{igabem}}$.}
	\label{Fig3:BCA_P_p5M2f100}
\end{figure}

\begin{figure}
	\centering
	\begin{subfigure}[t]{\textwidth}
		\includegraphics[width=\textwidth]{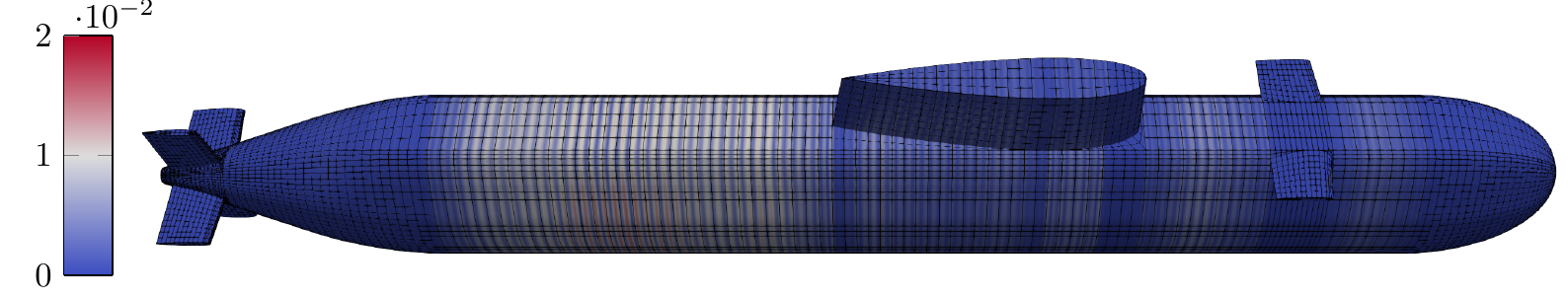}
		\caption{IGABA: $\displaystyle \max_{\vec{x}\in\Gamma_5}\frac{|p-p_h|}{|p|} = 0.014$}
		\label{Fig3:BCA_P_p5M2f1000_BA}
	\end{subfigure} 
	\par\bigskip
	\begin{subfigure}[t]{\textwidth}
		\includegraphics[width=\textwidth]{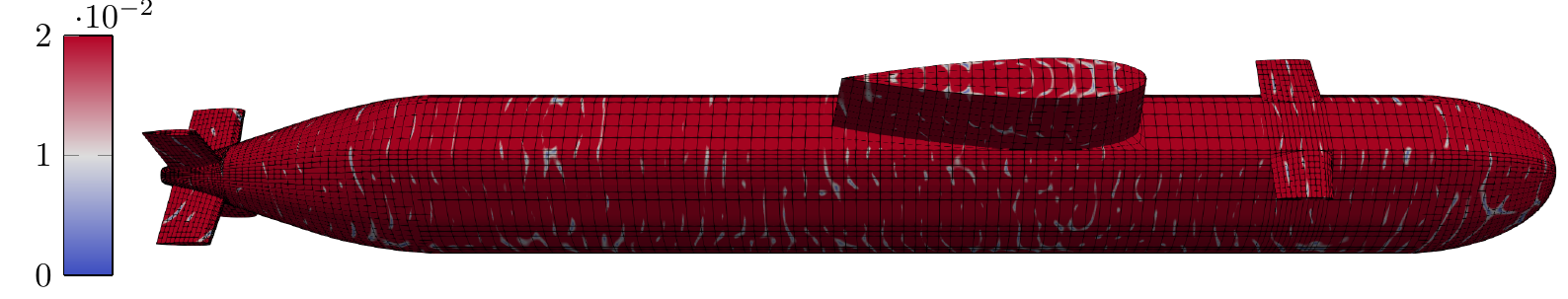}
		\caption{IGABEM CCBIE: $\displaystyle \max_{\vec{x}\in\Gamma_5}\frac{|p-p_h|}{|p|} = 0.40$}
	\end{subfigure} 
	\par\bigskip
	\begin{subfigure}[t]{\textwidth}
		\includegraphics[width=\textwidth]{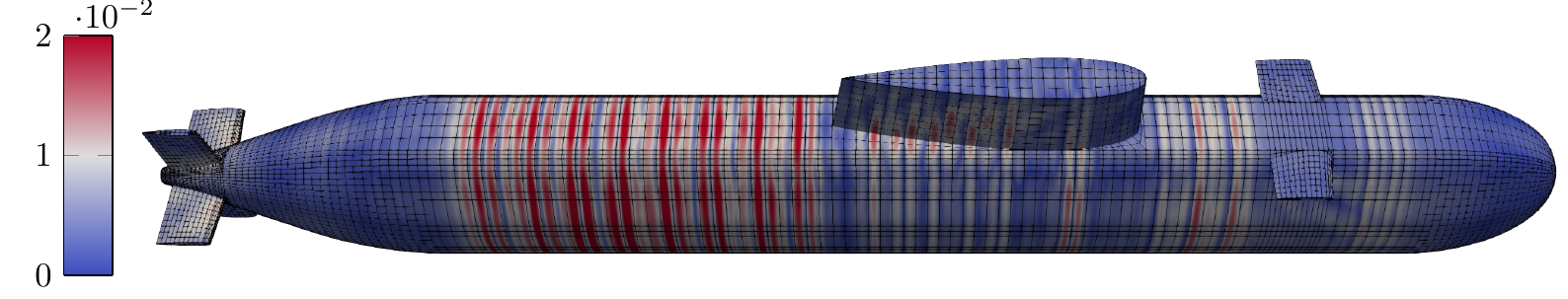}
		\caption{IGABEM CBM: $\displaystyle \max_{\vec{x}\in\Gamma_5}\frac{|p-p_h|}{|p|} = 0.030$}
	\end{subfigure} 
	\caption{\textbf{Manufactured solution with the BeTSSi submarine}: Relative error on the surface of the scatterer at $f=\SI{1000}{Hz}$ on the mesh ${\cal M}_{2,5,4}^{\textsc{igabem}}$.}
	\label{Fig3:BCA_P_p5M2f1000}
\end{figure}

Consider the manufactured solution~\Cref{Eq3:manuMFS} on the BeTSSi submarine with $N=16$ and where 16 source points are uniformly placed at the $x$-axis starting at $x=b$ and ending at $x=-L-2b$ (parameters taken from \Cref{Tab3:BeTSSiParameters}). The analytic real part of the pressure, $\Re{p}$, is visualized on the surface of the scatterer in~\Cref{Fig3:BCA_P_analytic}. 
A simulation at $f=\SI{100}{Hz}$ on mesh ${\cal M}_{1,2,1}^{\textsc{igabem}}$ yields the error plots in~\Cref{Fig3:BCA_P_p2M1f100}, which show good agreement between the best approximation and the BEM simulation. For more refined meshes in \Cref{Fig3:BCA_P_p5M1f100,Fig3:BCA_P_p2M2f100,Fig3:BCA_P_p5M2f100} (especially \Cref{Fig3:BCA_P_p5M2f100}) the numerical quadrature around the source points is too inaccurate. At this level of numerical accuracy, one quickly runs into issues due to round-off errors. The non-Lipschitz domains do not in and of itself pose any analysis suitable issues as described in~\Cref{Sec3:BeTSSi_approximation}, so the effect seen here is due to the numerical integration in the boundary element method. At $f=\SI{1000}{Hz}$ it is clear from \Cref{Fig3:BCA_P_p5M2f1000} that the IGABEM CCBIE simulation is polluted from a fictitious eigenfrequency. The remedy for this is to use the CBM formulation which obtains results with maximal error roughly twice the size of the best approximation. The meshes for the BeTSSi submarine in~\Cref{Fig3:BeTSSimeshes} might give the impression of evenly distributed control points in some areas, in particular the area behind the sail ($-L<x<x_{\mathrm{s}}-l_{\mathrm{ls}}$). In this case there are additional knot insertions around the submarine to obtain the $C^0$ lines, which results in ``bands'' of slightly larger errors along the submarine. This effect will be larger for higher polynomial orders, particularly for mesh ${\cal M}_{2,5,4}^{\textsc{igabem}}$ in~\Cref{Fig3:BCA_P_p5M2f1000_BA}.

\begin{figure}
	\centering    
	\begin{subfigure}[t]{0.49\textwidth}
		\centering
		\includegraphics{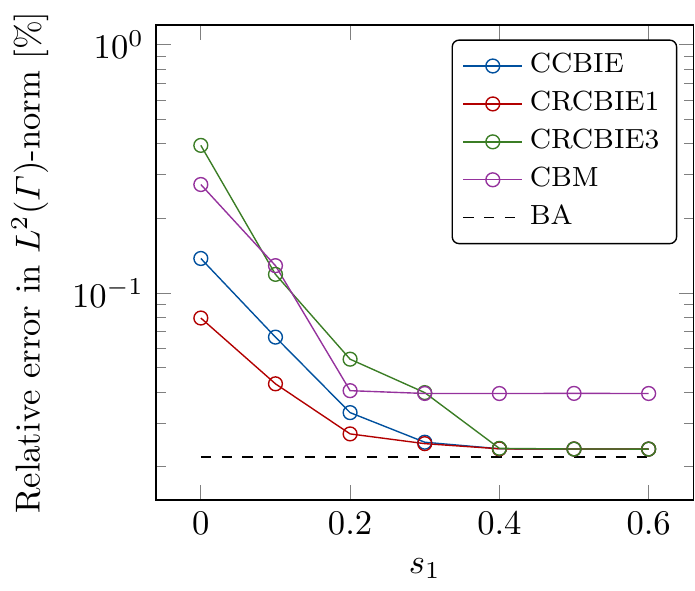}
		\caption{$\check{p} = 2$}
	\end{subfigure}%
	\hspace*{0.02\textwidth}%
	\begin{subfigure}[t]{0.49\textwidth}
		\centering
		\includegraphics{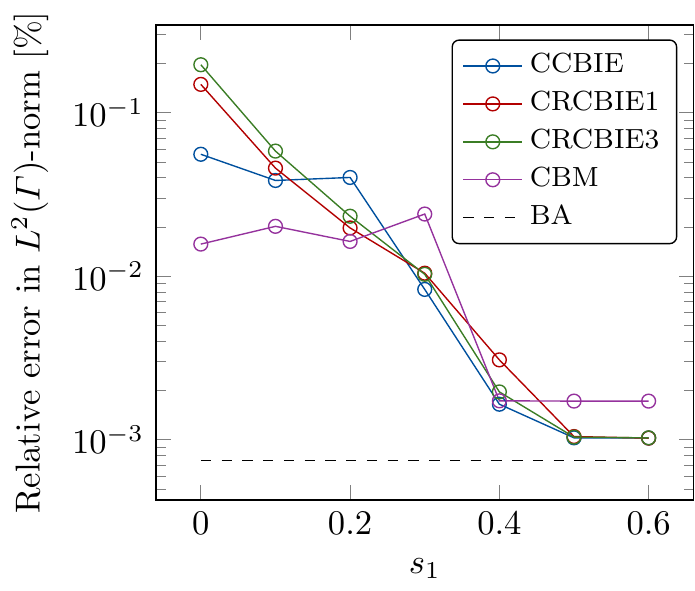}
		\caption{$\check{p} = 5$}
	\end{subfigure}
	\caption{\textbf{Manufactured solution with the BeTSSi submarine}: Surface error as a function of the parameter $s_1$, on the mesh ${\cal M}_{1,\check{p},\check{p}-1}^{\textsc{igabem}}$.}
	\label{Fig3:BCA_P_qp}
\end{figure}

To assess the parameter $s_1$ in~\Cref{Eq3:numberOfSubElements}, a low frequency of $\SI{100}{Hz}$ is now considered. In \Cref{Fig3:BCA_P_qp} we illustrate the effect of different choices of the parameter $s_1$ for the more complex geometry of the BeTSSi submarine. Again, the optimal choice for $s_1$ is polynomial dependent. Moreover, even the regularized formulations CRCBIE1 and CRCBIE3 must have $s_1>0$ contrary to what was proposed in~\cite{Sun2015bri} (stating that the singular free integrals ``can be evaluated by any convenient integration quadrature''). Whenever care is not taken for the numerical quadrature, incorrect conclusions may arise. This example illustrates the power of the manufactured solution as it enables computation of the best approximation such that the numerical integration may be controlled.


\subsection{Rigid scattering on the BeTSSi submarine}
\begin{figure}
	\centering
	\includegraphics[width=\textwidth]{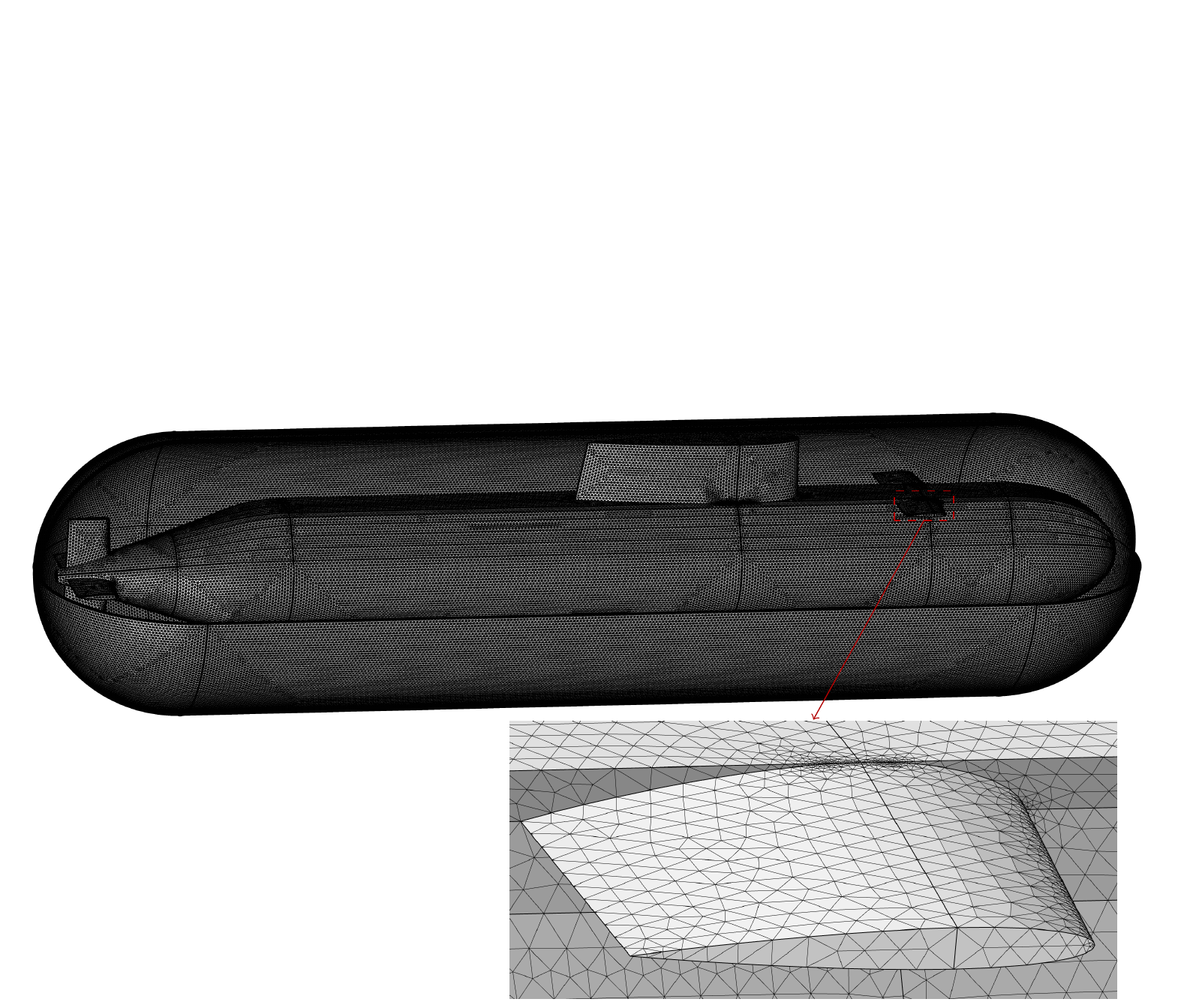}
	\caption{\textbf{Rigid scattering on the BeTSSi submarine}: Mesh used in \COMSOL simulations. The mesh consists of \num{27614929} second order finite elements including the elements in the PML (resulting in \num{43431671} degrees of freedom). This corresponds to 80 and 8 elements per wavelength for $\SI{100}{Hz}$ and $\SI{1000}{Hz}$, respectively.  The PML domain consists of a cylinder with two spherical end caps and are discretized by 10 layers of prismatic elements. The domain inside this PML is discretized with tetrahedral elements. The distance between the PML and the scatterer at the $x$-axis is $t_{\mathrm{a}}=\SI{1}{m}$ at both ends. The thickness of the PML is the same as the maximal tetrahedral diameter $h_{\mathrm{max}}=\SI{0.1875}{m}$. The PML cylinder starts at $x=-L-g_2-g_3+a$ and ends at $x=0$. The radius of the PML cylinder and the PML spherical end caps are $r_{\mathrm{a}} = a+t_{\mathrm{a}}$. The PML uses a polynomial coordinate stretching type with scaling factor and scaling curvature equal to 1. The simulations use \COMSOL version 5.4 with the acoustics module (to enable the PML method) and the design module (to import the CAD model).}
	\label{Fig3:COMSOL}
\end{figure}

\begin{figure}
	\centering
	\includegraphics[width=\textwidth]{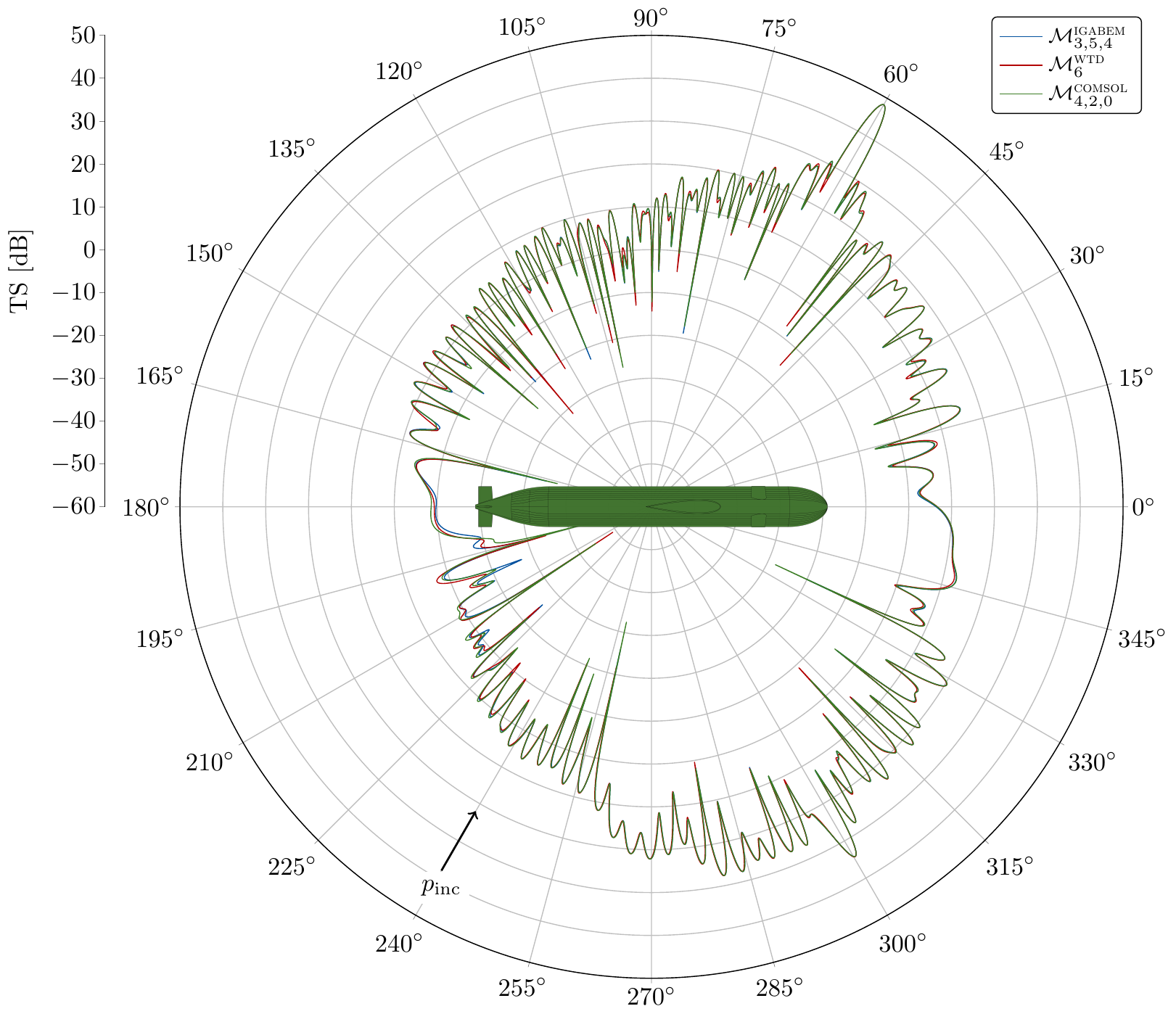}
	\caption{\textbf{Rigid scattering on the BeTSSi submarine}: Polar plot of the bistatic target strength ($\TS$) plotted against the azimuth angle $\varphi$ at $f=\SI{1000}{Hz}$. Direction of incident wave, $p_{\mathrm{inc}}$ is given by \Cref{Eq3:d_s} with $\alpha_{\mathrm{s}}=\ang{240}$ and $\beta_{\mathrm{s}}=\ang{0}$. The IGA mesh here used is ${\cal M}_{3,5,4}^{\textsc{igabem}}$. The \COMSOL simulation used \num{3.8} hours on mesh ${\cal M}_{4,2,0}^{\textsc{comsol}}$. The WTD 71 simulation was made using a direct BEM collocation method with the Burton--Miller formulation on mesh ${\cal M}_{5}^{\textsc{wtd}}$ described in \Cref{Sec3:BeTSSi_triangulation} with constant basis functions over each element.}
	\label{Eq3:polar_BI}
\end{figure}

\begin{figure}
	\centering
	\begin{subfigure}[t]{\textwidth}
		\includegraphics[width=\textwidth]{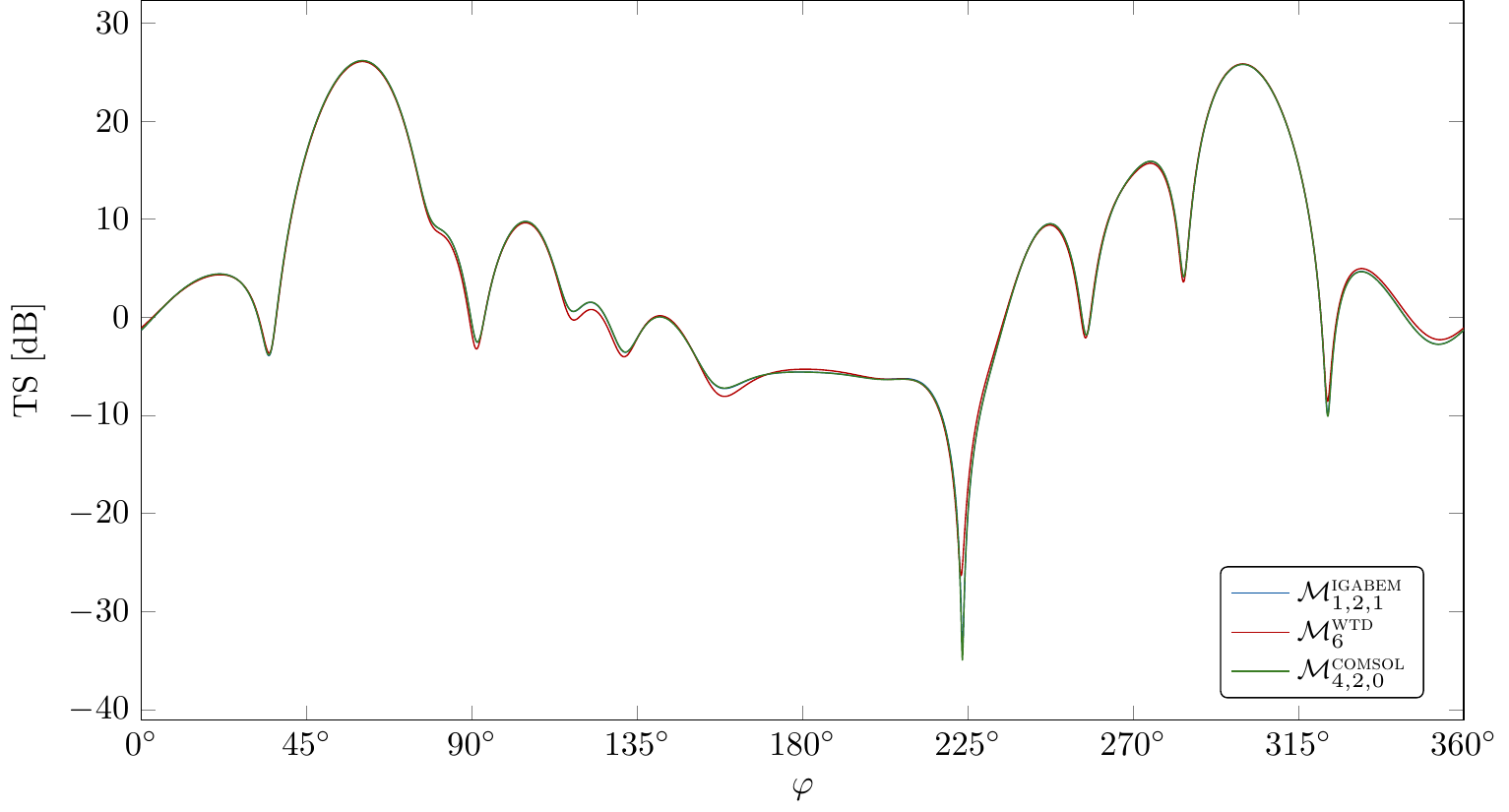}
		\caption{$f=\SI{100}{Hz}$}
		\label{Fig3:xy_BI_100}
	\end{subfigure} 
	\par\bigskip
	\par\bigskip
	\begin{subfigure}[t]{\textwidth}
		\includegraphics[width=\textwidth]{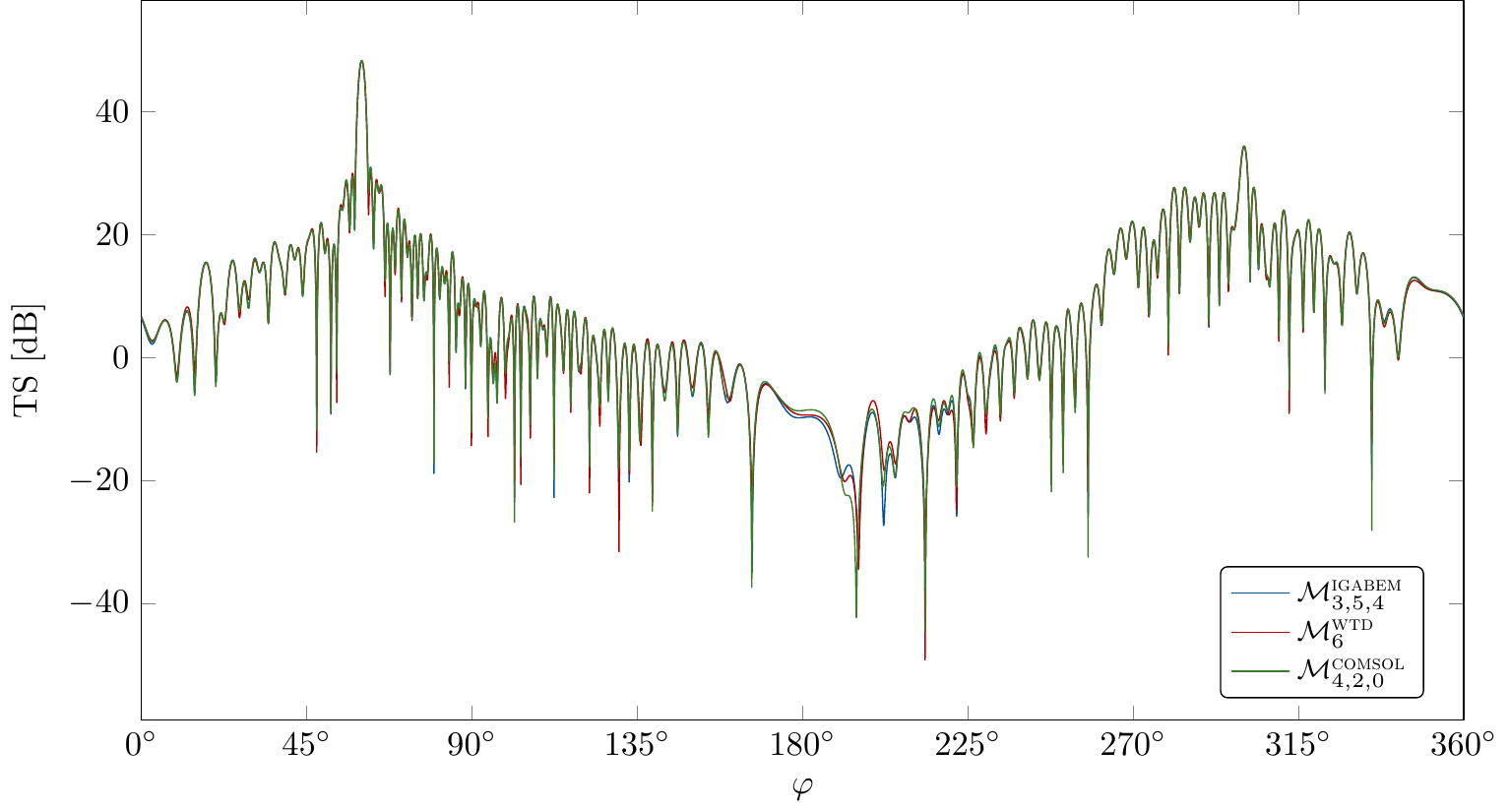}
		\caption{$f=\SI{1000}{Hz}$}
		\label{Fig3:xy_BI_1000}
	\end{subfigure} 
	\caption{\textbf{Rigid scattering on the BeTSSi submarine}: The bistatic target strength ($\TS$) plotted against the azimuth angle $\varphi$.}
\end{figure}

\begin{figure}
	\centering
	\includegraphics[width=\textwidth]{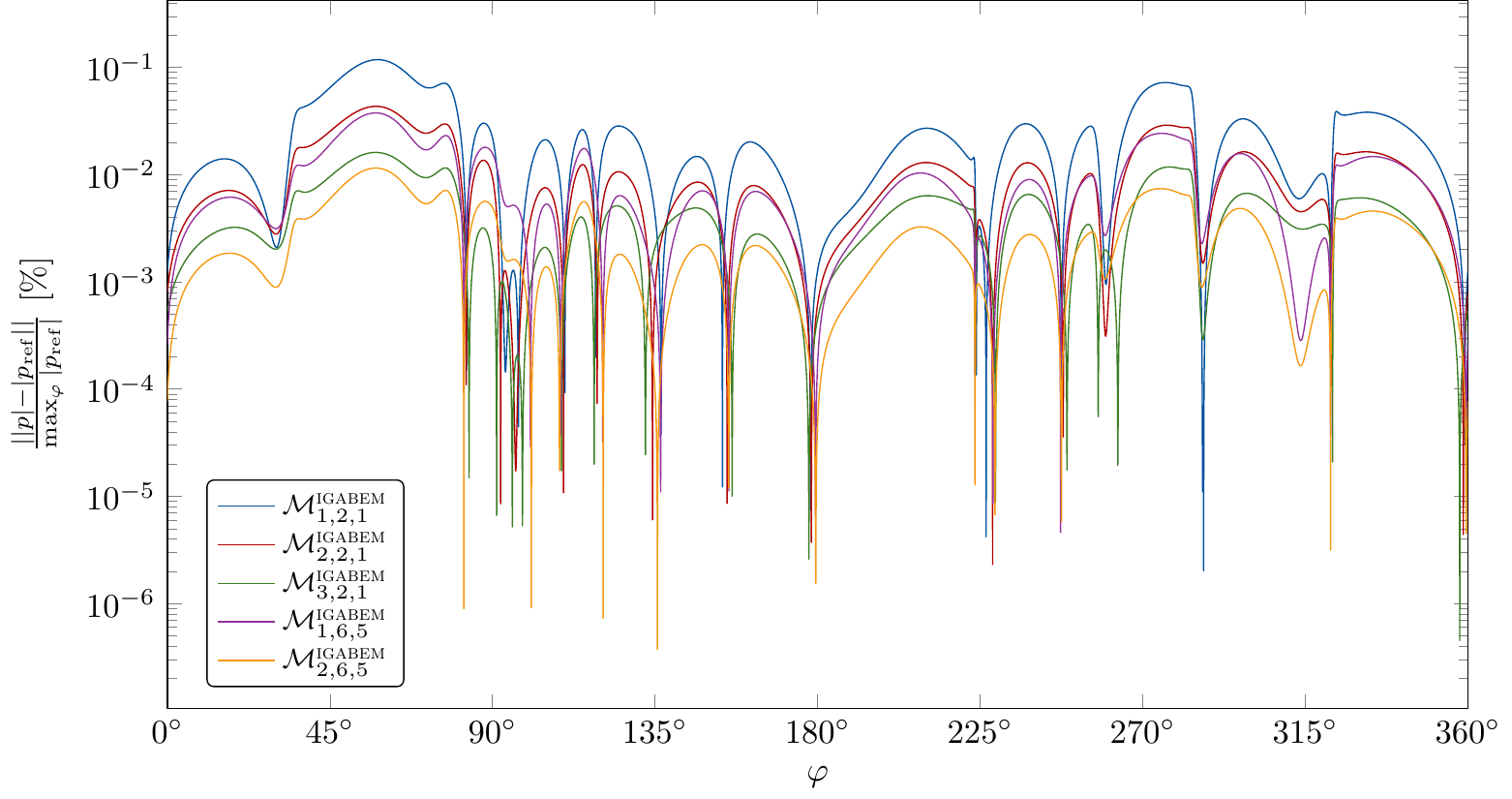}
	\caption{\textbf{Rigid scattering on the BeTSSi submarine}: The relative error in the far field absolute pressure plotted against the azimuth angle $\varphi$ at $f = \SI{100}{Hz}$, with the simulations from ${\cal M}_{3,6,5}^{\textsc{igabem}}$ as reference solution.}
	\label{Fig3:error_BI_100_IGA}
\end{figure}

\begin{figure}
	\centering
	\includegraphics[width=\textwidth]{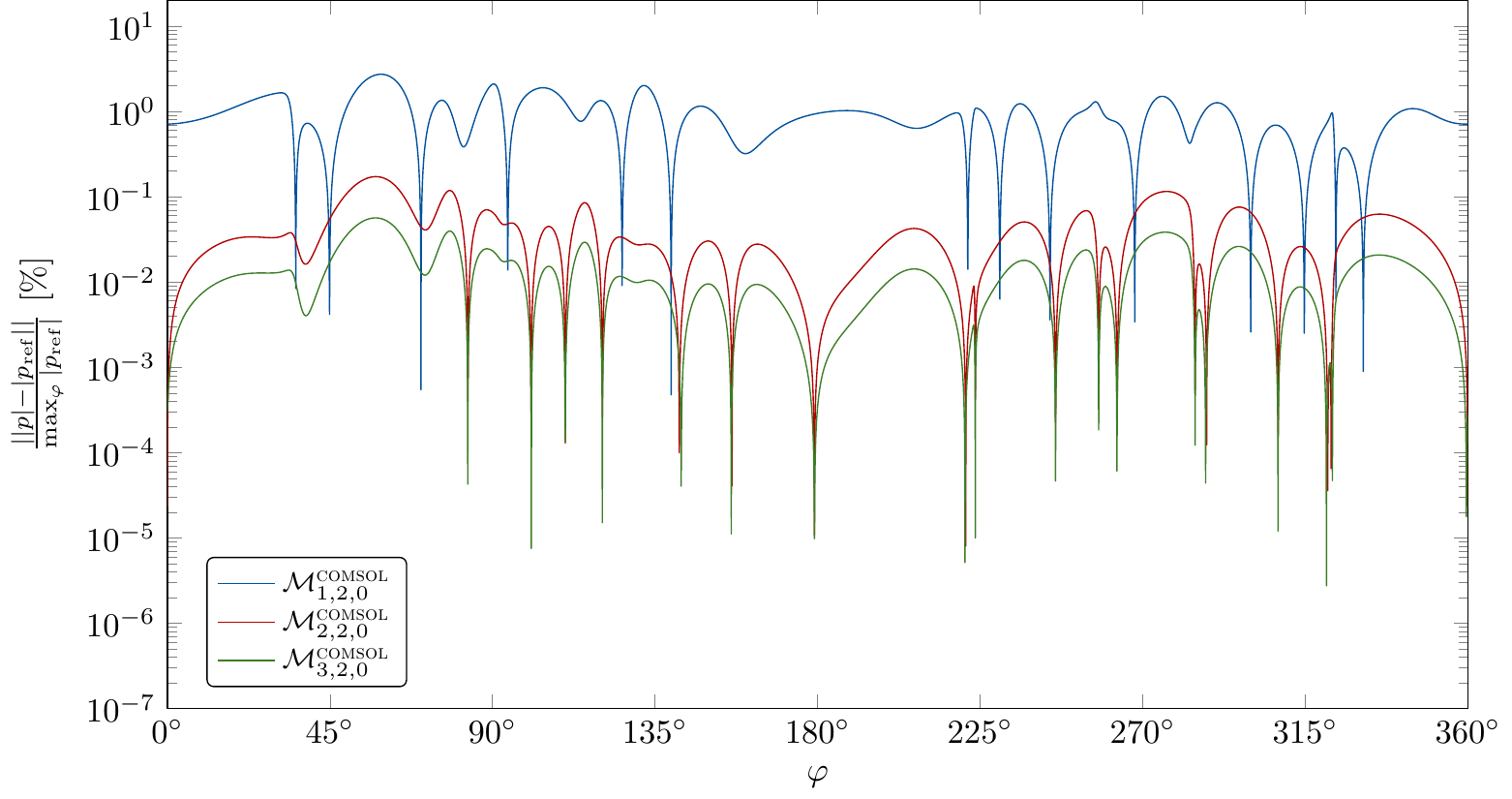}
	\caption{\textbf{Rigid scattering on the BeTSSi submarine}: The relative error in the far field absolute pressure plotted against the azimuth angle $\varphi$ at $f = \SI{100}{Hz}$, with the simulations from ${\cal M}_{4,2,0}^{\textsc{comsol}}$ as reference solution.}
	\label{Fig3:error_BI_100_COMSOL}
\end{figure}

\begin{figure}
	\centering
	\includegraphics[width=\textwidth]{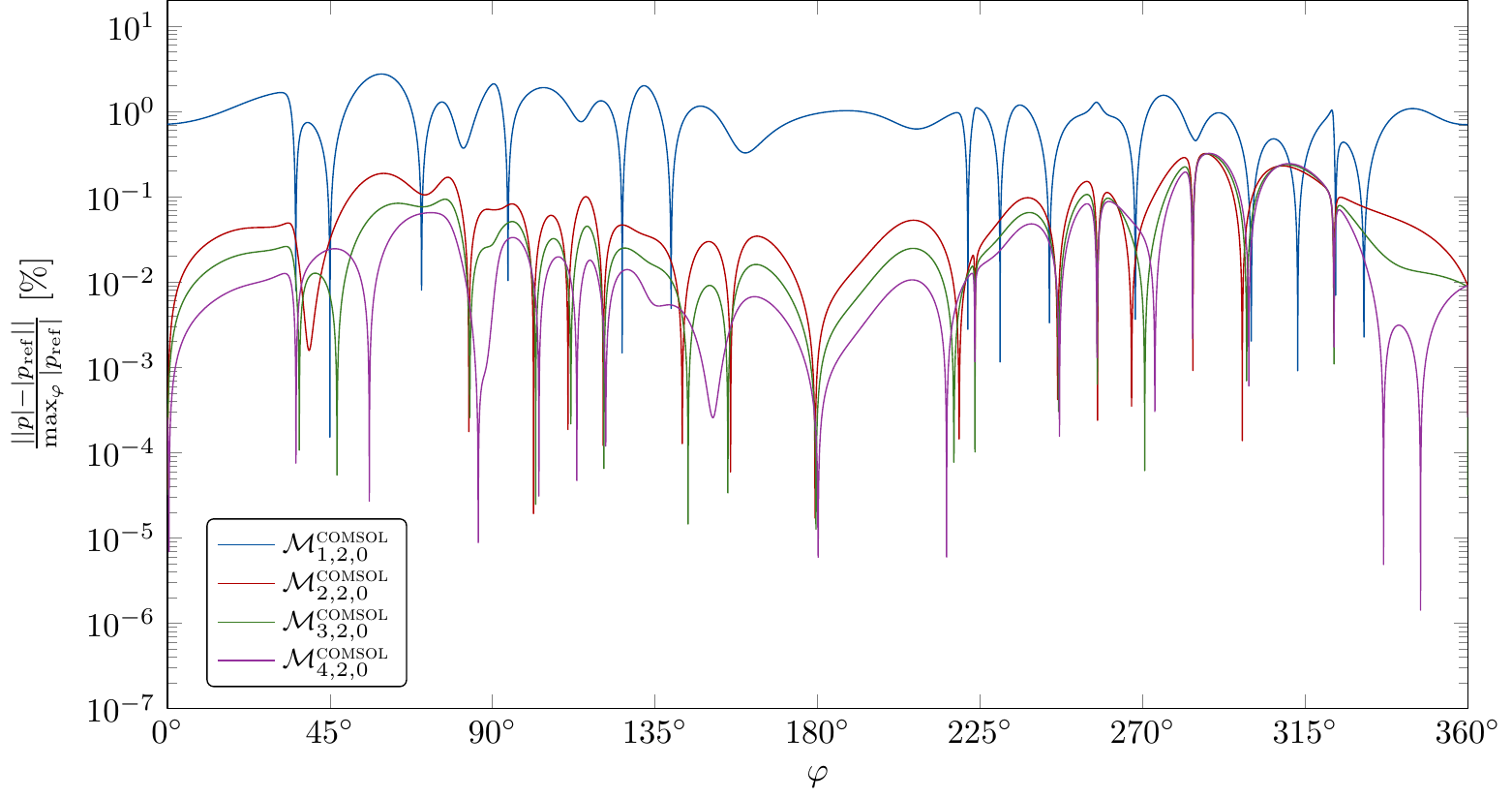}
	\caption{\textbf{Rigid scattering on the BeTSSi submarine}: The relative error in the far field absolute pressure plotted against the azimuth angle $\varphi$ at $f = \SI{100}{Hz}$, with the simulations from ${\cal M}_{3,6,5}^{\textsc{igabem}}$ as reference solution.}
	\label{Fig3:error_BI_100_COMSOL_IGA}
\end{figure}

Consider now a plane wave scattered by a rigid BeTSSi submarine. Throughout this section (motivated by the previous section) we use the CCBIE formulation at $f=\SI{100}{Hz}$ and the CBM formulation at $f=\SI{1000}{Hz}$. To verify our simulations, we compare with corresponding simulations done in \COMSOL, with mesh and parameters as illustrated and described in~\Cref{Fig3:COMSOL}. Comparisons are also made with simulations done by WTD 71\footnote{Wehrtechnische Dienststelle f\"{u}r Schiffe und Marinewaffen, Maritime Technologie und Forschung.}.
The polar plot in \Cref{Eq3:polar_BI} illustrates bistatic scattering where the incident wave is fixed, and the observation points for the far field computations sweep the aspect angles. A very good match is obtained, although some discrepancies are observed around the aft angles (around $\alpha=\ang{180}$). One can argue that the logarithmic scale of the target strength (TS) yields a somewhat misguided conception of the numerical error in the pressure. The pressure at these angles is very low such that the global relative error in the pressure is not as bad as the plot may suggest.
In \Cref{Fig3:xy_BI_100} and \Cref{Fig3:xy_BI_1000} the corresponding $xy$-plots are given at $\SI{100}{Hz}$ and $\SI{1000}{Hz}$, respectively.
In~\Cref{Fig3:xy_BI_100} (at $\SI{100}{Hz}$) the IGA and \COMSOL simulations are visually indistinguishable, such that error plots are in order. Let the simulation from ${\cal M}_{3,6,5}^{\textsc{igabem}}$, ${\cal M}_{4,2,0}^{\textsc{comsol}}$ and ${\cal M}_{6}^{\textsc{wtd}}$ be a reference solution for IGABEM, \COMSOL and WTD71, respectively. In~\Cref{Fig3:error_BI_100_IGA} we compare the IGA results for lower resolved meshes. Convergence throughout the aspect angles is observed. In~\Cref{Fig3:error_BI_100_COMSOL} a corresponding comparison is done with the \COMSOL simulations. Better convergence rates for higher polynomial degrees in the IGA simulations are not present. This is probably due to the problem of numerical integration over the non-Lipschitz domains as discussed in~\Cref{Subsec3:manufactured}. Another reason could be the need for adaptive refinement, for example using LR B-splines~\cite{Johannessen2014iau} based on a posteriori error estimates, e.g.\ by exploiting $k$-refinement as presented in~\cite{Kumar2015sap}.
This might also be the reason that the \COMSOL simulations converge to a different solution around $\varphi=\ang{280}$ as illustrated in~\Cref{Fig3:error_BI_100_COMSOL_IGA}.
In~\Cref{Tab3:BeTSSiComputationalData} we present the computational complexity of the different simulations. The number of degrees of freedom per wavelength is denoted by $\tau$. We shall use another definition of $\tau$ compared to the definition found in~\cite[p. 767]{Peake2015eib}\footnote{Here, $\tau$ is defined as $\tau = \lambda\sqrt{n_{\mathrm{dof}}/|\Gamma|}$.}, namely the minimal number of degrees of freedom per wavelength (instead of an average). This is arguably a better definition as it more precisely captures how well the frequency is resolved. We compute $\tau$ by
\begin{equation*}
	\tau=\frac{\lambda}{d_{\mathrm{max}}},\quad d_{\mathrm{max}} = \max_{\vec{x}\in X}\min_{\vec{y}\in X\setminus\vec{x}}\|\vec{x}-\vec{y}\|
\end{equation*}
where $X$ is the set of nodes in the mesh. For IGA these nodes are chosen to be the Greville points in the physical domain (as the control points do not lie on the geometry). For the \COMSOL simulations we get $\tau=\frac{\lambda}{h_{\mathrm{max}}/2}$ and for constant triangular elements (WTD 71 simulations) we get ${\tau=\frac{\lambda}{2h_{\mathrm{max}}/3}}$. Considering the error as a function of $\tau$, IGA outperforms the simulations from both \COMSOL and WTD 71. Even considering the error as a function of time usage, the IGA simulations obtain comparable results despite the sub-optimal implementation discussed earlier.
\begin{table}
	\centering
	\caption{\textbf{Rigid scattering on the BeTSSi submarine}: Data for the meshes used in the BeTSSi simulations at $f=\SI{100}{Hz}$. The error is a relative $l^2$-error of the absolute far field pressure with the simulation from ${\cal M}_{3,6,5}^{\textsc{igabem}}$, ${\cal M}_{4,2,0}^{\textsc{comsol}}$ and ${\cal M}_{6}^{\textsc{wtd}}$ as a reference solution for IGABEM, \COMSOL and WTD71, respectively. The IGABEM and \COMSOL simulations were computed on 28 Intel CPUs ($2\times 24$-core Xeon 2.6 GHz) with 768 GB RAM available and the WTD71 simulations were computed on a 32 core Xeon computer with 2.3 GHz.}
	\label{Tab3:BeTSSiComputationalData}
	\begin{tabular}{c S[table-format = 8.0] S[table-format = 8.0] S[table-format = 1.2,round-mode=places,round-precision=2] S[table-format = 3.1,round-mode=places,round-precision=1] S[table-format = 2.4,round-mode=places,round-precision=4] S[table-format = 6.0]}
		\toprule
		Mesh & {$n_{\mathrm{el}}$} & {$n_{\mathrm{dofs}}$} & {$h_{\mathrm{max}}$[$\si{m}$]}  & {$\tau$ [$\si{m^{-1}}$]}  & {Error [\%]}  & {$t_{\mathrm{tot}}$ [$\si{s}$]}\\
		\hline
		${\cal M}_{1,2,1}^{\textsc{igabem}}$ & 3718 & 6725 & 1.65059 & 17.0461 & 0.117586 & 227\\
		${\cal M}_{2,2,1}^{\textsc{igabem}}$ & 14872 & 20521 & 0.827789 & 30.5556 & 0.046594 & 2611\\
		${\cal M}_{3,2,1}^{\textsc{igabem}}$ & 59488 & 70421 & 0.433037 & 52.9412 & 0.0184557 & 34244\\
		${\cal M}_{1,6,5}^{\textsc{igabem}}$ & 3718 & 27537 & 1.65059 & 25.5368 & 0.0393652 & 1789\\
		${\cal M}_{2,6,5}^{\textsc{igabem}}$ & 14872 & 52293 & 0.82779 & 34.0537 & 0.012171 & 11860\\
		${\cal M}_{3,6,5}^{\textsc{igabem}}$ & 59488 & 124113 & 0.433251 & 61.5283 & {-} & 108741\\
%
		${\cal M}_{1,2,0}^{\textsc{comsol}}$ & 100436 & 250638 & 2.21 & 13.6 & 3.3481 & 10\\
		${\cal M}_{2,2,0}^{\textsc{comsol}}$ & 550300 & 1167195 & 1.14 & 26.3 & 0.1703 & 38\\
		${\cal M}_{3,2,0}^{\textsc{comsol}}$ & 3729303 & 6654972 & 0.60 & 50.0 & 0.0569 & 375\\
		${\cal M}_{4,2,0}^{\textsc{comsol}}$ & 27614929 & 43431671 & 0.32 & 93.8 & {-} & 5650\\
		${\cal M}_{1}^{\textsc{wtd}}$ 		& 4140 & 4140 & 1.89333 & 11.88 & 2.4013 & 2\\
		${\cal M}_{2}^{\textsc{wtd}}$ 		& 10406 & 10406 & 1.00496 & 22.3890 & 1.8815 & 8\\
		${\cal M}_{3}^{\textsc{wtd}}$ 		& 31104 & 31104 & 0.498592 & 45.1271 & 1.2824 & 25\\
		${\cal M}_{4}^{\textsc{wtd}}$ 		& 106888 & 106888 & 0.256928 & 87.5732 & 1.0328 & 38\\
		${\cal M}_{5}^{\textsc{wtd}}$ 		& 400886 & 400886 & 0.130041 & 173.0224 & 0.6598 & 112\\
		${\cal M}_{6}^{\textsc{wtd}}$ 		& 1584014 & 1584014 & 0.0691461 & 325.3980 & {-} & 400\\
		\bottomrule
	\end{tabular}
\end{table}

\begin{figure}
	\centering
	\includegraphics[width=\textwidth]{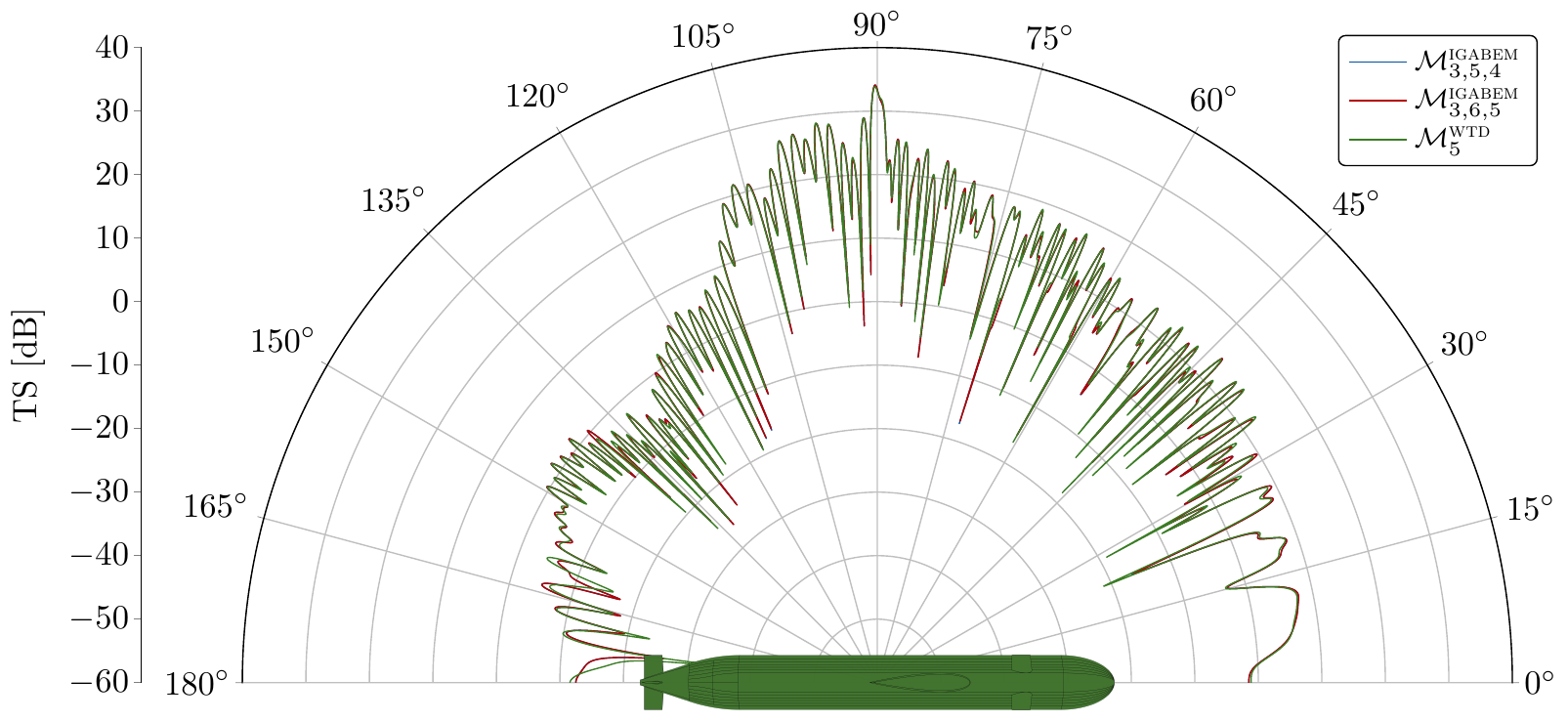}
	\caption{\textbf{Rigid scattering on the BeTSSi submarine}: Polar plot of the monostatic target strength ($\TS$) at $f=\SI{1000}{Hz}$ plotted against the azimuth angle $\varphi$. All simulations use the CBM formulation.} 
	\label{Fig3:BCA_MS}
\end{figure}

A monostatic\footnote{The incident wave has the same origin as the far field point in a monostatic sweep.} polar plot is shown in \Cref{Fig3:BCA_MS} at $f=\SI{1000}{Hz}$. The results for ${\cal M}_{3,5,4}^{\textsc{igabem}}$ and ${\cal M}_{3,6,5}^{\textsc{igabem}}$ are practically indistinguishable in this plot. A comparison is made with a simulation done by WTD 71 showing good agreement. The $l^2$-error of the absolute far field pressure for ${\cal M}_{3,5,4}^{\textsc{igabem}}$ (with ${\cal M}_{3,6,5}^{\textsc{igabem}}$ as reference solution) is about 0.052\%. The corresponding error for the WTD simulation is 5.5\%. Using a direct solver for the IGA simulations, monostatic scattering can easily be solved with multiple right-hand sides (in the present case 3601 column vectors that correspond to 3601 distinct azimuth angles $\varphi\in[0,\ang{180}]$ with steps of $\ang{0.5}$). The time consumption for monostatic scattering is then increased by less than 1\% compared to bistatic scattering since the most computationally complex operation here is to build the system of equations. The WTD 71 simulation solves the 3601 cases individually, resulting in a time consumption increase of about 1392\% (the computations used \num{43.3} hours on a 32 core Xeon computer with 2.3 GHz). The reason that number is not 7201\% (WTD 71 timings are here for all angles in $[0,\ang{360}]$) is because WTD 71 uses a precondition matrix based on the result from 5 neighboring monostatic angles.

\begin{figure}
	\centering    
	\begin{subfigure}[b]{\textwidth}
		\centering
		\includegraphics[width=\textwidth]{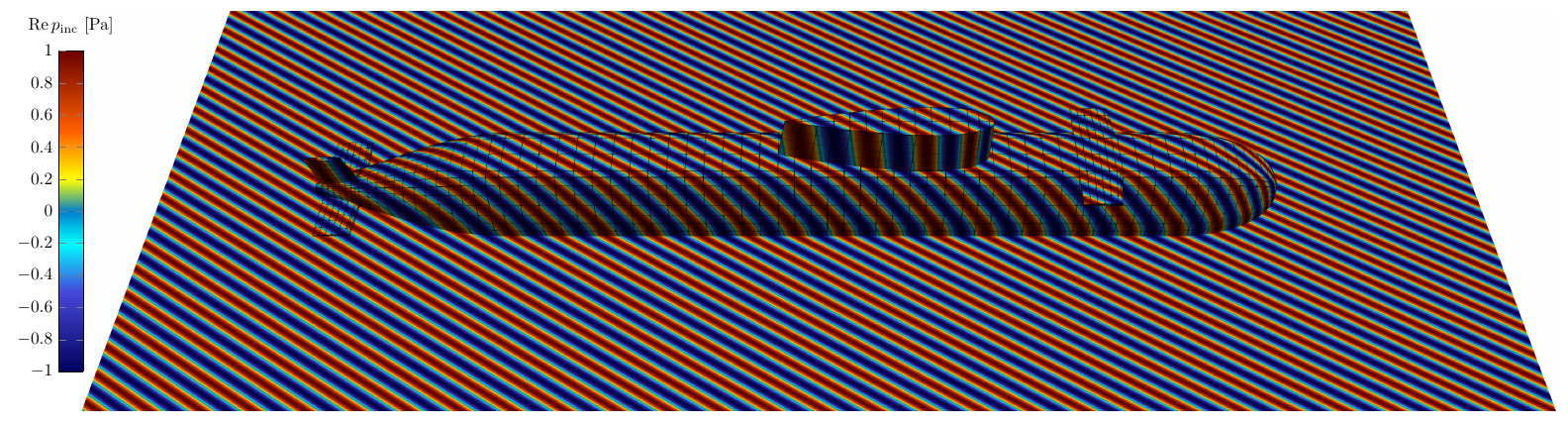}
		\caption{Real part of the incident wave $p_{\mathrm{inc}}(\vec{x})=P_{\mathrm{inc}}\euler^{\imag k\vec{d}_{\mathrm{s}}\cdot \vec{x}}$.}
	\end{subfigure}
	\par\smallskip
	\begin{subfigure}[b]{\textwidth}
		\centering
		\includegraphics[width=\textwidth]{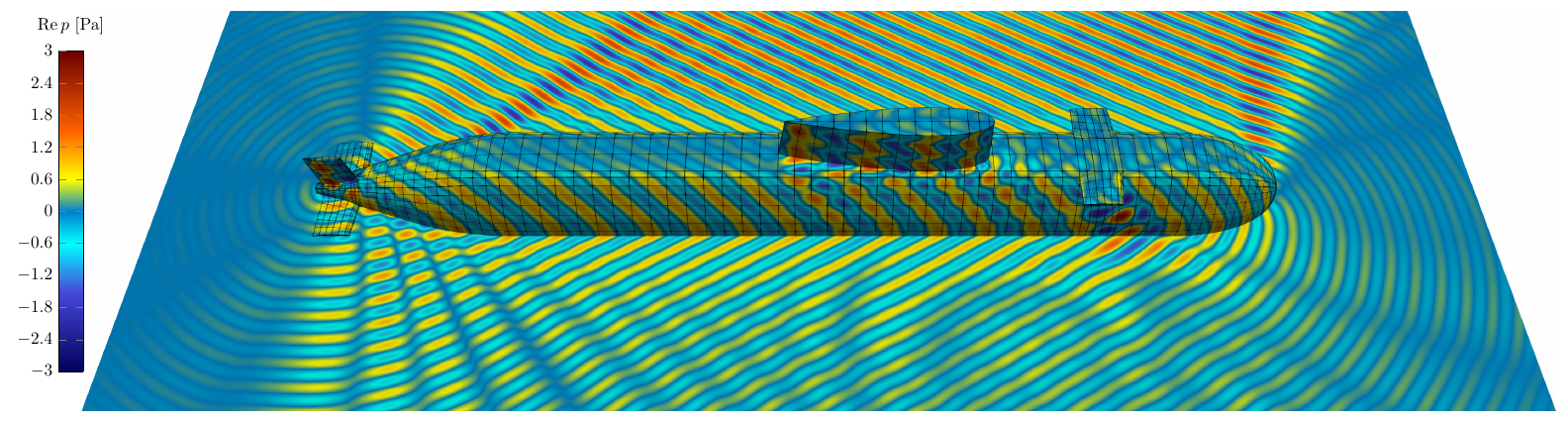}
		\caption{Real part of the scattered pressure $p(\vec{x})$.}
	\end{subfigure}
	\par\smallskip
	\begin{subfigure}[b]{\textwidth}
		\centering
		\includegraphics[width=\textwidth]{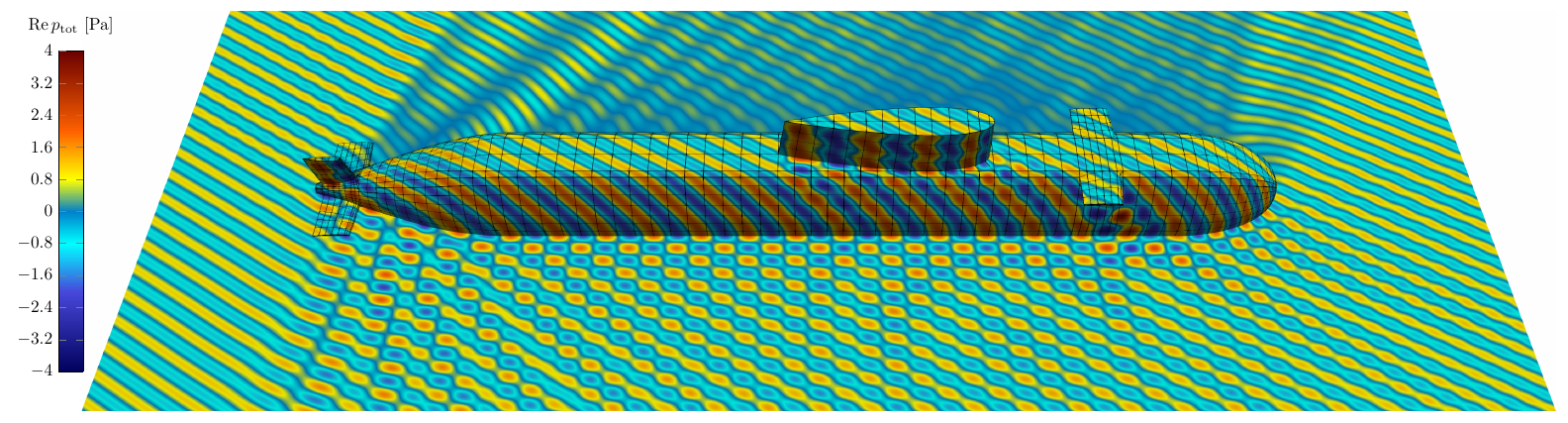}
		\caption{Real part of the total pressure $p_{\mathrm{tot}}(\vec{x})=p_{\mathrm{inc}}(\vec{x})+p(\vec{x})$.}
	\end{subfigure}
	\par\smallskip
	\begin{subfigure}[b]{\textwidth}
		\centering
		\includegraphics[width=\textwidth]{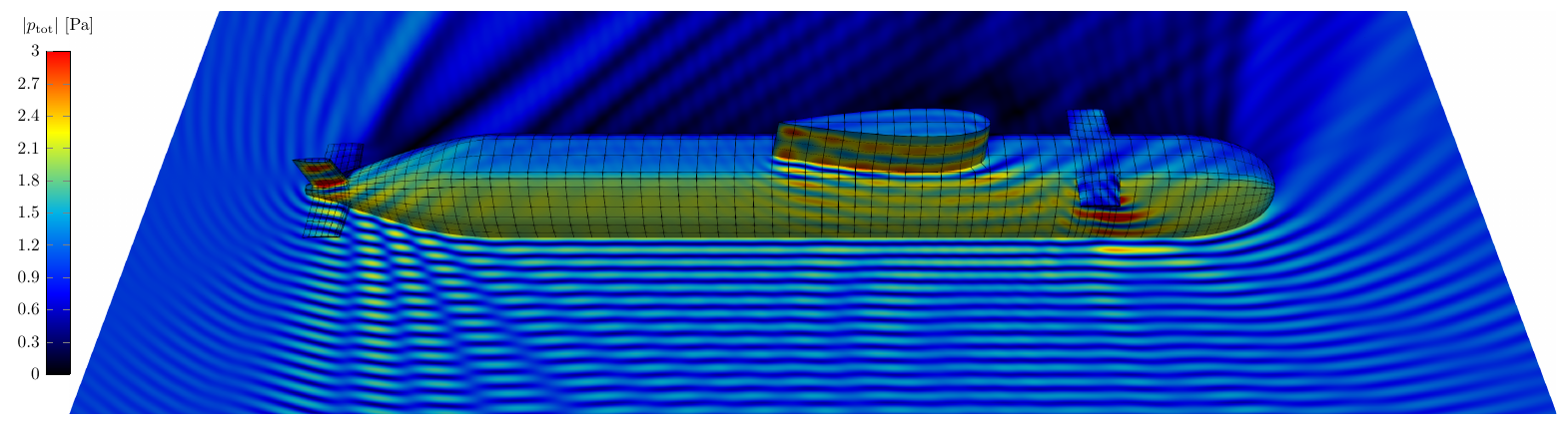}
		\caption{Modulus of the total pressure $p_{\mathrm{tot}}(\vec{x})=p_{\mathrm{inc}}(\vec{x})+p(\vec{x})$.}
		\label{Fig3:BC_NearField_abs}
	\end{subfigure}
	\caption{\textbf{Rigid scattering on the BeTSSi submarine}: The simulation at $f=\SI{1000}{Hz}$ is visualized in the $xy$-plane (and on the scatterer), and is computed on mesh ${\cal M}_{3,5,4}^{\textsc{igabem}}$. For visualization purposes, the mesh ${\cal M}_{1,5,4}^{\textsc{igabem}}$ is here visualized.}
	\label{Fig3:BC_NearField}
\end{figure}

Finally, the near field at $f=\SI{1000}{Hz}$ is visualized in \Cref{Fig3:BC_NearField}. From \Cref{Fig3:BC_NearField_abs} one can observe that the incident wave is reflected multiple times beneath the right depth rudder.

\clearpage
\section{Conclusions}
\label{Sec3:conclusions}
This article addresses acoustic scattering characterized by sound waves reflected by man-made elastic objects. The present approach is characterized by:
\begin{itemize}
	\item The scatterer is discretized using isogeometric analysis (IGA), which enables discretization directly from the basis functions used in the computer aided design (CAD) description of the model.
	\item Both collocation and Galerkin method are considered in combination with several boundary integral equation (BIE) formulations including the conventional (CBIE) formulation and the Burton--Miller (BM) formulation.
	\item The method of manufactured solution is used as a quality insurance.
\end{itemize}

The main finding of the present study is that the use of IGA significantly increases the accuracy compared to the use of $C^0$ finite element analysis (FEA) due to increased inter-element continuity of the spline basis functions.

Furthermore, the following observations are made
\begin{itemize}
	\item IGA's ability to represent the geometry exactly was observed to be of less importance for accuracy when comparing to higher order ($\hat{p}\geq 2$) isoparametric FEA. However, a more significant improvement offered by IGA is due to higher continuity of the spline basis functions in the solution space.
	\item For linear approximation of the geometry using classical boundary element method (BEM) the convergence order is reduced for higher order sub parametric elements.
	\item For resolved meshes, the IGA framework enables roughly the same accuracy per element (compared to higher order isoparametric FEA) even though the number of degrees of freedom is significantly reduced.
	\item IGA is more computationally efficient than FEA to obtain highly accurate solutions. That is, when the mesh is sufficiently resolved, a given accuracy is obtained computationally faster using IGA.
	\item Reduced accuracy is obtained for collocation simulations compared to Galerkin simulations, especially for the hypersingular BIE (HBIE) formulation and BM formulation. Better located collocation points may remedy this difference and is suggested as future work.
	\item The method of manufactured solution enables a convenient method of checking the mesh quality and to some extent the numerical accuracy of the rigid body scattering problem. It can be used to check the presence of fictitious eigenfrequencies.
	\item The improved adaptive integration procedure presented in this work uses significantly less quadrature points than the integration procedure presented in~\cite{Simpson2014aib} for a given accuracy.
	\item The presence of non-Lipschitz domain does not in principle cause problems for the analysis suitability of the problem as the best approximation is not significantly affected by such areas. However, for the boundary element method, the integral over singular kernels in such domain may cause problems. This is especially the case for highly accurate solution as round-off errors may become significant.
	\item Regularizing the weakly singular integrands in the BIEs does not eliminate the need for special quadrature rules around the source points. The small reduction in the number of quadrature points needed for the three versions of the regularized conventional BIE (RCBIE1, RCBIE2 and RCBIE3) formulations compared to the CBIE formulation is arguable not significant.
	\item Using the collocation method, an advantage for the CBIE formulation compared with the regularized formulations (RCBIE1, RCBIE2 and RCBIE3) is that there is no need to compute the normal vector at the collocation point for the CBIE formulation which could be problematic if the geometric mapping is singular at that point (as is the case for the north and south pole of the parametrization in \Cref{Fig3:parm1} and several locations for the BeTSSi submarine).
	\item The Galerkin method obtains results remarkably close to the best approximation combined with any formulation, illustrating the sharpness of the a priori error estimate in~\Cref{Eq3:aprioriErrorEstimate}.
\end{itemize}
The Burton-Miller formulation yields somewhat reduced accuracy in combination with the collocation method, which is the cost of removing fictitious eigenfrequencies. Another popular alternative is the combined Helmholtz integral formulation (CHIEF) framework which does not have this reduction in accuracy but has other downsides. By adding more constraints to the linear system of equations, the CHIEF method can remove fictitious eigenfrequencies with the cost of having to solve an over determined linear system of equations (using for example least squares). The main disadvantage with the CHIEF framework, however, is arguably the difficulty of finding interior points at which to evaluate the BIEs. This is especially problematic for high frequencies. An approach for solving this issue was made in~\cite{Wu1991awr}. The results in this work may be improved even further with the discontinuous IGABEM~\cite{Sun2019dib}.

The boundary element method is the method of choice in the BeTSSi community for obtaining accurate results for the BeTSSi submarine, mainly to avoid surface-to-volume parametrization. Although IGABEM seems to be a prominent framework to solve acoustic scattering problems, there are still issues on the BeTSSi submarine that was not resolved in this paper, in particular the integration procedure over non-Lipschitz areas on the BeTSSi submarine.

\section*{Acknowledgements}
This work was supported by the Department of Mathematical Sciences at the Norwegian University of Science and Technology and by the Norwegian Defence Research Establishment.

The publication of the BeTSSi models (storage of large data files) was provided by UNINETT Sigma2 --- the National Infrastructure for High Performance Computing and Data Storage in Norway.

The authors would like to thank Jan Ehrlich and Ingo Schaefer (WTD 71) for their simulations on the BeTSSi submarine and additional fruitful discussions.

The authors would also like to thank the reviewers for detailed response and many constructive comments.

\clearpage
\appendix
\section{NURBS parametrization of the sphere}
Two standard ways of parametrizing a sphere using NURBS are given below for the unit sphere (a simple scaling generalizes this for spheres of arbitrary radii). The first is represented by 8 elements in a single patch (only one element is given below, as the others are obtained by symmetry), and the second is represented by 6 patches (only one patch is given below, as the others are obtained by symmetry).
\subsection{Parametrization 1}
\label{Sec3:NURBSsphere1}
\begin{figure}
	\centering
	\begin{subfigure}{0.40\textwidth}
		\centering
		\includegraphics[width=0.9\textwidth]{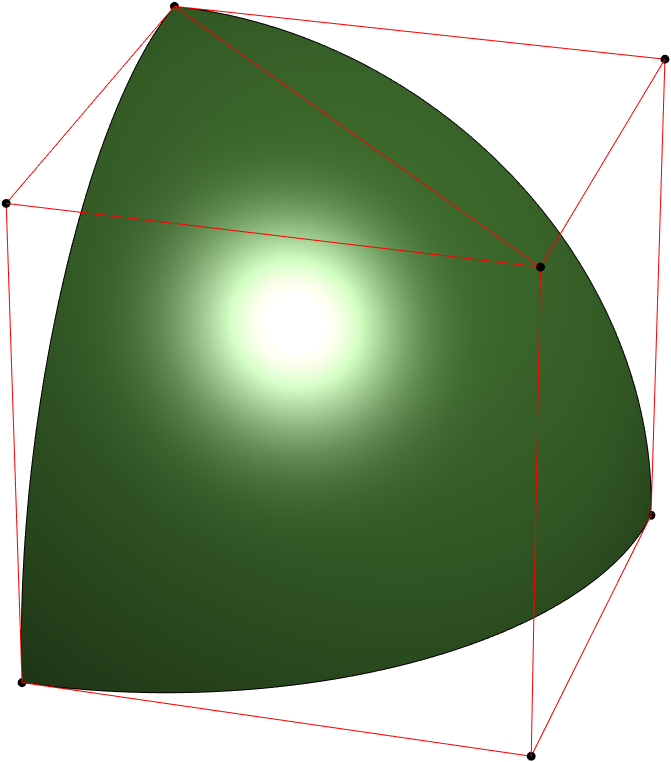}
		\caption{Parametrization 1}
		\label{Fig3:parametrizationOfSphere1} 
	\end{subfigure}
	~
	\begin{subfigure}{0.55\textwidth}
		\centering
		\includegraphics[width=0.9\textwidth]{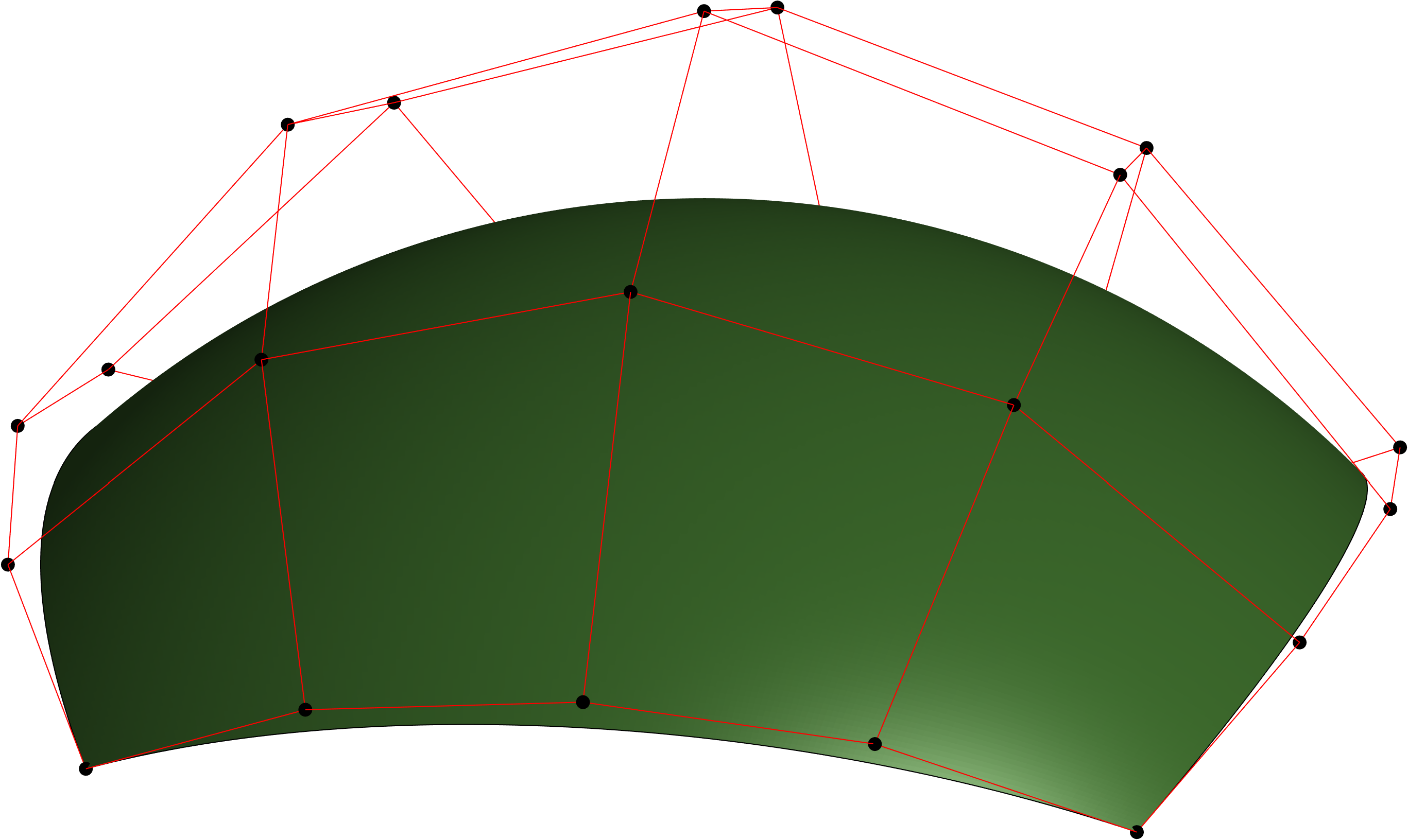}
		\caption{Parametrization 2}
		\label{Fig3:parametrizationOfSphere2} 
	\end{subfigure}
	\caption{\textbf{NURBS parametrization of the sphere}: Two NURBS parametrizations of the sphere. The control polygon is also shown.}       
\end{figure}

The sphere can be exactly parametrized by 8 NURBS elements of degree 2. One of these elements with corresponding control points is illustrated in \Cref{Fig3:parametrizationOfSphere1}. The weights and control points are given in \Cref{Tab3:sphere1} (a parametrization of all elements in a single patch can be found in~\cite[p. 168]{Venas2015iao}).

\begin{table}
	\centering
	\caption{\textbf{Parametrization 1}: Weights and control points for an element of a unit sphere.}
	\label{Tab3:sphere1}
	\begin{tabular}{c c c c c c}
		\toprule
		$i$		& 	$j$	& 	$x_{i,j}$ 	& $y_{i,j}$ 	& $z_{i,j}$ 	& $w_{i,j}$\\
		\hline
		$1$		&	$1$	&	1	& 0		& 0		&	1				\\
		$2$		&	$1$	&	1	& 1		& 0 	&	$1/\sqrt{2}$	\\
		$3$		&	$1$	&	0	& 1		& 0		&	1				\\ \\
		
		$1$		&	$2$	&	1	& 0		& 1		&	$1/\sqrt{2}$				\\
		$2$		&	$2$	&	1	& 1		& 1 	&	$1/2$	\\
		$3$		&	$2$	&	0	& 1		& 1		&	$1/\sqrt{2}$				\\ \\
		
		$1$		&	$3$	&	0	& 0		& 1		&	1				\\
		$2$		&	$3$	&	0	& 0		& 1 	&	$1/\sqrt{2}$	\\
		$3$		&	$3$	&	0	& 0		& 1		&	1				\\
		\bottomrule
	\end{tabular}
\end{table}

\subsection{Parametrization 2}
\label{Sec3:NURBSsphere2}
The sphere can be exactly parametrized \cite[p. 11]{Cobb1988tts} by 6 NURBS patches of degree 4. One of these patches with corresponding control points is illustrated in \Cref{Fig3:parametrizationOfSphere2}. Some of the weights and weighted control points are given in \Cref{Tab3:sphere2}. The remaining data is found by symmetry about the planes $x=0$, $y=0$, $y=x$ and $y=-x$. In particular (by symmetry about the $y=x$ plane)
\begin{equation*}
	x_{i,j} = y_{j,i},\quad y_{i,j} = x_{j,i},\quad z_{i,j} = z_{j,i},\quad w_{i,j} = w_{j,i}
\end{equation*}
for the pairs $(i,j) \in\{(1,2),(1,3),(2,3)\}$, and (by symmetry about the $y=0$ plane)
\begin{equation*}
	x_{i,j} = -x_{6-i,j},\quad y_{i,j} = y_{6-i,j},\quad z_{i,j} = z_{6-i,j},\quad w_{i,j} = w_{6-i,j}
\end{equation*}
for $i=4,5$ and $j=1,2,3$, and then (by symmetry about the $x=0$ plane)
\begin{equation*}
	x_{i,j} = x_{i,6-j},\quad y_{i,j} = -y_{i,6-j},\quad z_{i,j} = z_{i,6-j},\quad w_{i,j} = w_{i,6-j}
\end{equation*}
for $i=1,2,3,4,5$ and $j=4,5$.

\begin{table}
	\centering
	\caption{\textbf{Parametrization 2}: Weights and weighted control points for a tile of a unit sphere.}
	\label{Tab3:sphere2}
	\begin{tabular}{c c c c c c}
		\toprule
		$i$		& 	$j$	& 	$w_{i,j}x_{i,j}$ 	& $w_{i,j}y_{i,j}$ 	& $w_{i,j}z_{i,j}$ 	& $w_{i,j}$\\
		\hline
		$1$		&	$1$	&	$4(1-\sqrt{3})$ 	& $4(1-\sqrt{3})$ 			& $4(\sqrt{3}-1)$		&	$4(3-\sqrt{3})$				\\
		$2$		&	$1$	&	$-\sqrt{2}$ 		& $\sqrt{2}(\sqrt{3}-4)$ 	& $\sqrt{2}(4-\sqrt{3})$&	$\sqrt{2}(3\sqrt{3}-2)$				\\
		$3$		&	$1$	&	$0$ 				& $4(1-2\sqrt{3})/3$ 		& $ 4(2\sqrt{3}-1)/3$	&	$4(5-\sqrt{3})/3$				\\ \\

		$2$		&	$2$	&	$-(3\sqrt{3}-2)/2$ 	& $(2-3\sqrt{3})/2$ 		& $(\sqrt{3}+6)/2$		&	$(\sqrt{3}+6)/2$				\\
		$3$		&	$2$	&	$0$ 				& $\sqrt{2}(2\sqrt{3}-7)/3$ & $5\sqrt{6}/3$		&	$\sqrt{2}(\sqrt{3}+6)/3$\\ \\

		$3$		&	$3$	&	$0$ 				& $0$ 						& $4(5-\sqrt{3})/3$		&	$4(5\sqrt{3}-1)/9$				\\
		\bottomrule
	\end{tabular}
\end{table}

\section{NURBS parametrization of the torus}
\label{Sec3:torus}
\begin{figure}
	\centering
	\includegraphics[width=0.7\textwidth]{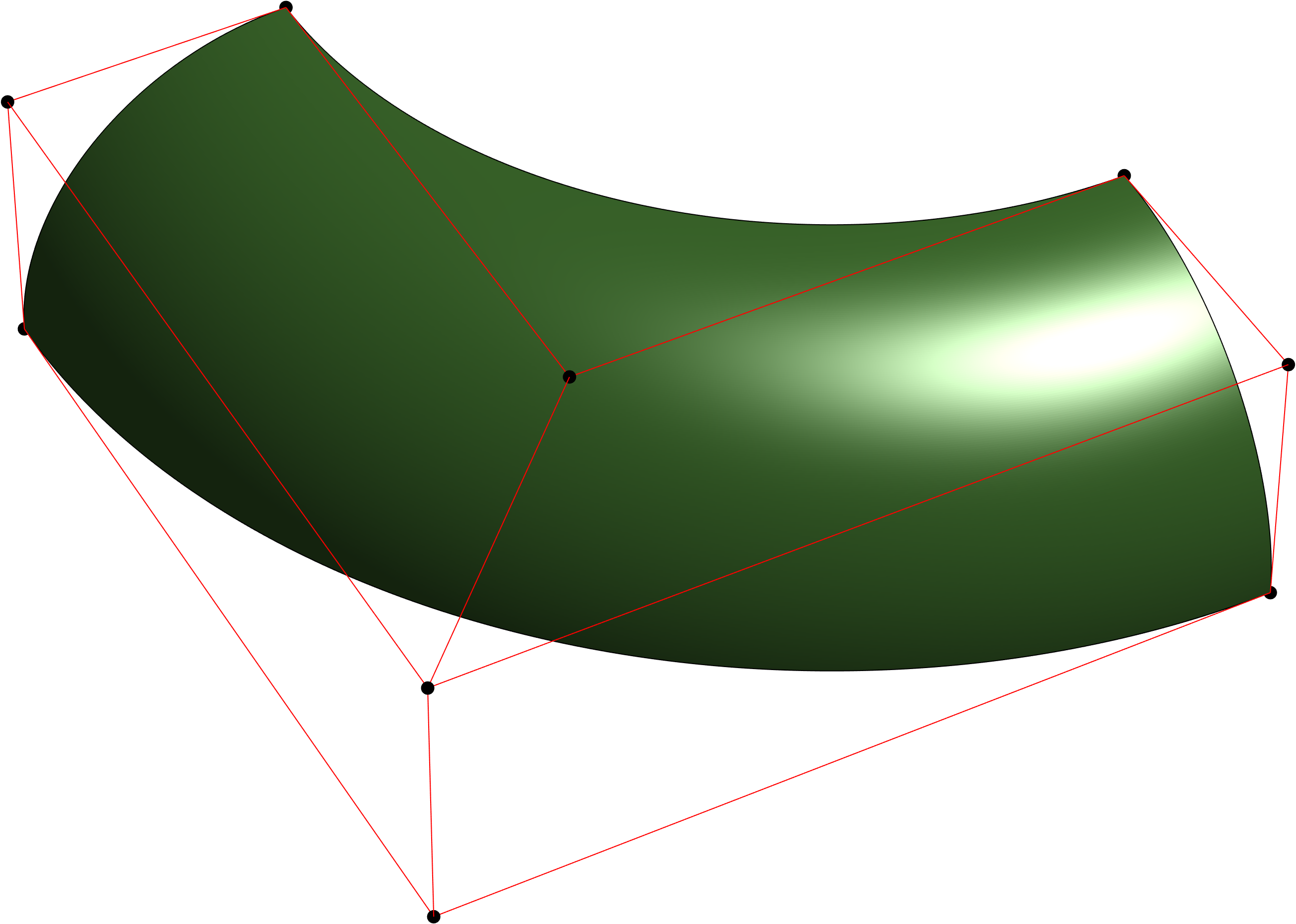}
	\caption{\textbf{NURBS parametrization of the torus}: A NURBS parametrization of a 1/16 of a torus. The control polygon is also shown.}       
	\label{Fig3:parametrizationOfTorus} 
\end{figure}

A torus with major radius $r_{\mathrm{o}}$ and minor radius $r_{\mathrm{i}}$ can be represented by a single NURBS patch with 16 elements (as visualized in \Cref{Fig3:Torus}). One of these elements is shown in~\Cref{Fig3:parametrizationOfTorus} with corresponding control polygon. The weights and control points are given in \Cref{Tab3:torus}.

\begin{table}
	\centering
	\caption{\textbf{NURBS parametrization of the torus}: Weights and control points for the torus.}
	\label{Tab3:torus}
	\begin{tabular}{c c c c c c}
		\toprule
		$i$		& 	$j$	& 	$x_{i,j}$ 	& $y_{i,j}$ 	& $z_{i,j}$ 	& $w_{i,j}$\\
		\hline
		$1$		&	$1$	&	$r_{\mathrm{o}}+r_{\mathrm{i}}$	& 0		& 0		&	1				\\
		$2$		&	$1$	&	$r_{\mathrm{o}}+r_{\mathrm{i}}$	& $r_{\mathrm{o}}+r_{\mathrm{i}}$		& 0 	&	$1/\sqrt{2}$	\\
		$3$		&	$1$	&	0	& $r_{\mathrm{o}}+r_{\mathrm{i}}$		& 0		&	1				\\ \\
		
		$1$		&	$2$	&	$r_{\mathrm{o}}+r_{\mathrm{i}}$	& 0		& $r_{\mathrm{i}}$		&	$1/\sqrt{2}$				\\
		$2$		&	$2$	&	$r_{\mathrm{o}}+r_{\mathrm{i}}$	& $r_{\mathrm{o}}+r_{\mathrm{i}}$		& $r_{\mathrm{i}}$ 	&	$1/2$	\\
		$3$		&	$2$	&	0	& $r_{\mathrm{o}}+r_{\mathrm{i}}$		& $r_{\mathrm{i}}$		&	$1/\sqrt{2}$				\\ \\
		
		$1$		&	$3$	&	$r_{\mathrm{o}}$	& 0		& $r_{\mathrm{i}}$		&	1				\\
		$2$		&	$3$	&	$r_{\mathrm{o}}$	& $r_{\mathrm{o}}$		& $r_{\mathrm{i}}$ 	&	$1/\sqrt{2}$	\\
		$3$		&	$3$	&	0	& $r_{\mathrm{o}}$	& $r_{\mathrm{i}}$		&	1				\\
		\bottomrule
	\end{tabular}
\end{table}
\section{The BeTSSi submarine model}
\label{Sec3:BeTSSi_description}
\begin{figure}
	\centering
	\includegraphics[width=\textwidth]{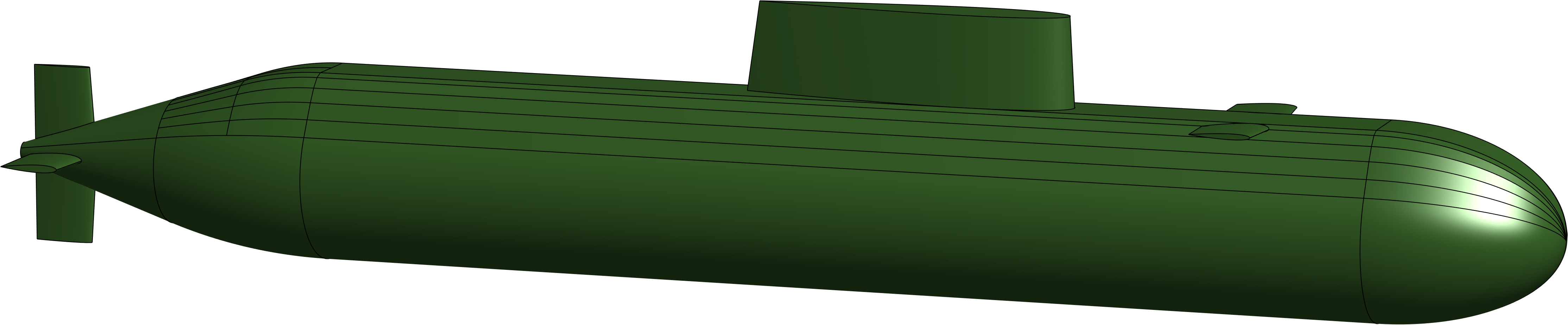}
	\caption{Outer pressure hull for BeTSSi submarine.}
	\label{Fig3:BeTSSi_BC}
\end{figure}

In this section the BeTSSi \cite{Nolte2014bib} submarine model (depicted in \Cref{Fig3:BeTSSi_BC}) will be presented. The BeTSSi submarine contains many standard designing features including circles, ellipses, straight panels, cylinders and cones. In addition, several NACA profiles are present giving a very nice benchmark model for sub-surface scattering. For the analysis part, it contains challenges such as trimming curves and non-Lipschitz domains~\cite{Lipton2010roi}. All in all, a challenging benchmark without being too complex.

The original BeTSSi submarine model presented in  \cite{Nolte2014bib} contains several discrepancies that is arguably not optimal for a benchmark model. 
First, the NACA profiles used to create the sail and the rudders are only given with 5 digits of accuracy. This in turn, results in for example the sail not being tangent to the side lines of the deck with an error of around $\SI{1}{mm}$. This creates problems for the meshing procedure as this results in either very small elements in this area, or element with high aspect ratios. 
Second, the exact geometry for the upper transition from the deck to the rotationally symmetric cone tail, is hidden by an ``internal routine in ANSYS''. Not only is this hard to reproduce for anyone without an ANSYS license, but the available CAD file for this model does not represent the transition to the lower part exactly (as this curve should be a circular arc and is not represented by a NURBS curve). In order to create a watertight model, the available CAD file approximates the lower transition such that the side curves match. 

The relevant BeTSSi parameters for the work presented herein are given in \Cref{Tab3:BeTSSiParameters}.
\begin{table}
	\centering
	\caption{\textbf{BeTSSi submarine:} Free parameters for the BeTSSi submarine benchmark.}
	\label{Tab3:BeTSSiParameters}
	\begin{tabular}{l l}
		\toprule
		Parameter & Description\\
		\midrule
		$\alpha=\ang{18}$ & Arc angle of transition to the tail cone\\
		$\beta=\ang{240}$ & Rotational angle for the axisymmetric lower part of the pressure hull\\
		$g_2=\SI{6.5}{m}$ & Distance in the $x$-direction of transition to the tail cone\\
		$g_3=\SI{6.5}{m}$ & Distance in the $x$-direction of the tail cone\\
		$L=\SI{42}{m}$ & Length of the deck\\
		$a=\SI{7}{m}$ & Semi-major axis of bow\\
		$b=\SI{3.5}{m}$ & Semi-minor axis of bow\\
		$c=\SI{4}{m}$ & Height from the $x$-axis to the deck\\
		$s=\SI{1.2}{m}$ & Half of the width of the deck\\
		$l_{\mathrm{ls}}=\SI{13}{m}$ & Length of the lower cross-section of the sail\\
		$l_{\mathrm{lm}}=\SI{2.6}{m}$ & Length of the lower cross-section of the main rudders\\
		$l_{\mathrm{ld}}=\SI{2.6}{m}$ & Length of the lower cross-section of the depth rudders\\
		$l_{\mathrm{us}}=\SI{12.3}{m}$ & Length of the upper cross-section of the sail\\
		$l_{\mathrm{um}}=\SI{2.35}{m}$ & Length of the upper cross-section of the main rudders\\
		$l_{\mathrm{ud}}=\SI{2.35}{m}$ & Length of the upper cross-section of the depth rudders\\
		$b_{\mathrm{lm}}=\SI{0.4}{m}$ & Width of the lower cross-section of the main rudders\\
		$b_{\mathrm{us}}=\SI{2}{m}$ & Width of the upper cross-section of the sail\\
		$b_{\mathrm{um}}=\SI{0.3}{m}$ & Width of the upper cross-section of the main rudders\\
		$b_{\mathrm{ud}}=\SI{0.22}{m}$ & Width of the upper cross-section of the depth rudders\\
		$\delta_{\mathrm{s}}=\SI{0.2}{m}$ & Parameter for shifting the upper and lower cross-section of the sail\\
		$h_{\mathrm{s}}=\SI{3.5}{m}$ & Height of the sail\\
		$h_{\mathrm{m}}=\SI{3.5}{m}$ & Height of the main rudders\\
		$x_{\mathrm{s}}=-\SI{12}{m}$ & Positioning of the sail\\
		$x_{\mathrm{m}}=-\SI{51.9}{m}$ & Positioning of the main rudders\\
		$x_{\mathrm{d}}=-\SI{4}{m}$ & Positioning of the depth rudders\\
		\bottomrule
	\end{tabular}
\end{table}

\subsection{Main body}
\begin{figure}
	\centering
	\includegraphics[scale=1]{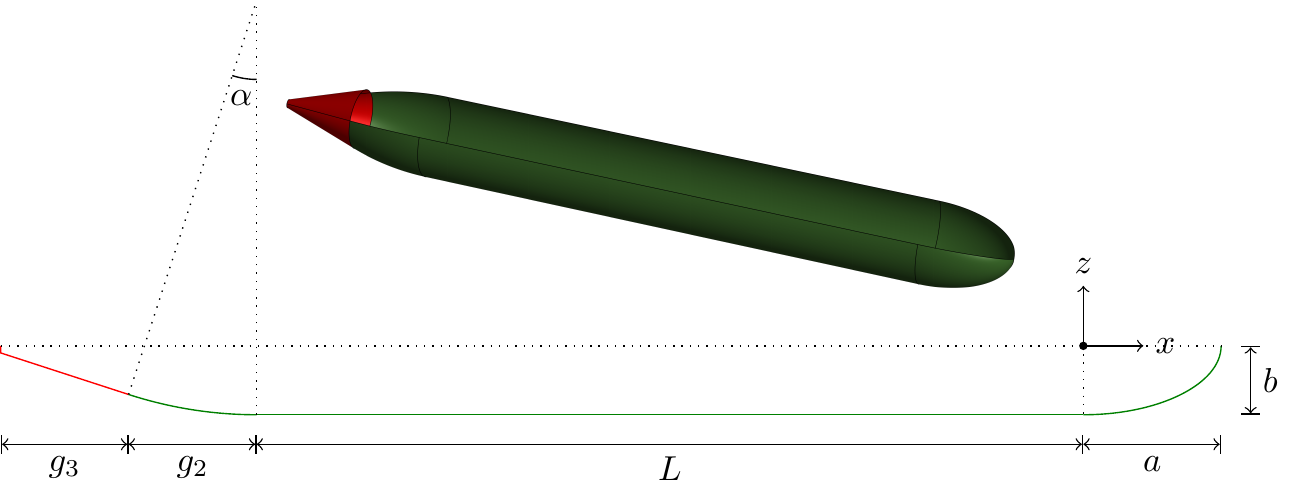}
	\caption{The sideline of the lower part of the BeTSSi submarine. The side lines are formed (from the right) by an ellipse with semi-major axis $a$ and semi-minor axis $b$, followed by a straight line of length $L$, then an arc of angle $\alpha$ and finally two straight lines. The latter two straight lines (in red) are rotated about the $x$-axis and the remaining part (in green) are rotated an angle $\beta$ around the $x$-axis.}
	\label{Fig3:bettsi_bottom}
\end{figure}

\begin{figure}
	\centering
	\includegraphics[scale=1]{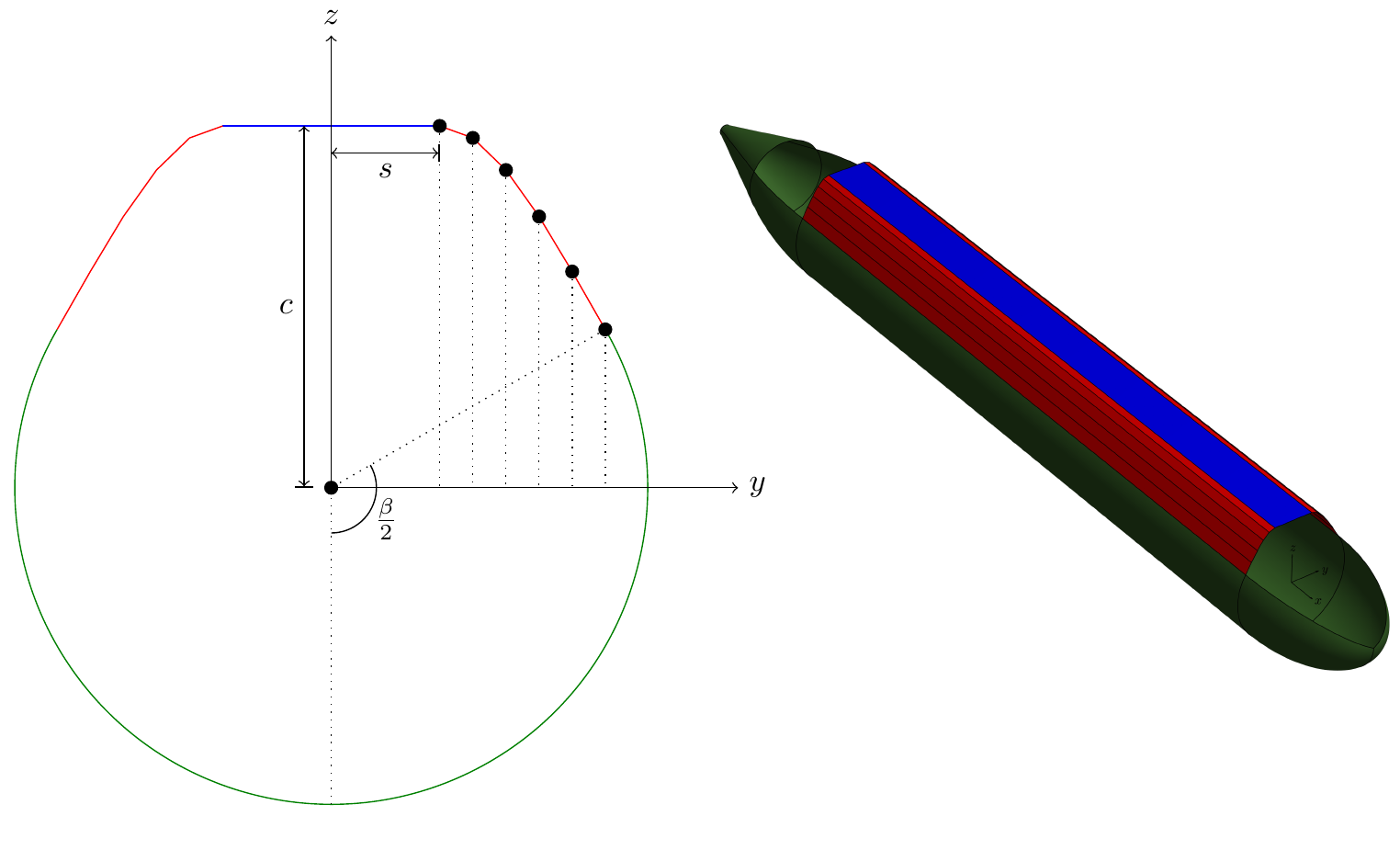}
	\caption{The transition (red line) from the axisymmetric hull (green line) to the deck (blue line) is given by sampling a cubic polynomial, $P_{\mathrm{p}}(y)$, at 6 equidistant points in the $y$-direction and connecting the resulting points with straight lines (corresponding 6 points are found for negative values $y$-values, $(0,y,P_{\mathrm{p}}(|y|))$).}
	\label{Fig3:bettsi_top}
\end{figure}

\begin{figure}
	\centering    
	\begin{subfigure}{0.49\textwidth}
		\centering
		\includegraphics[width=\textwidth]{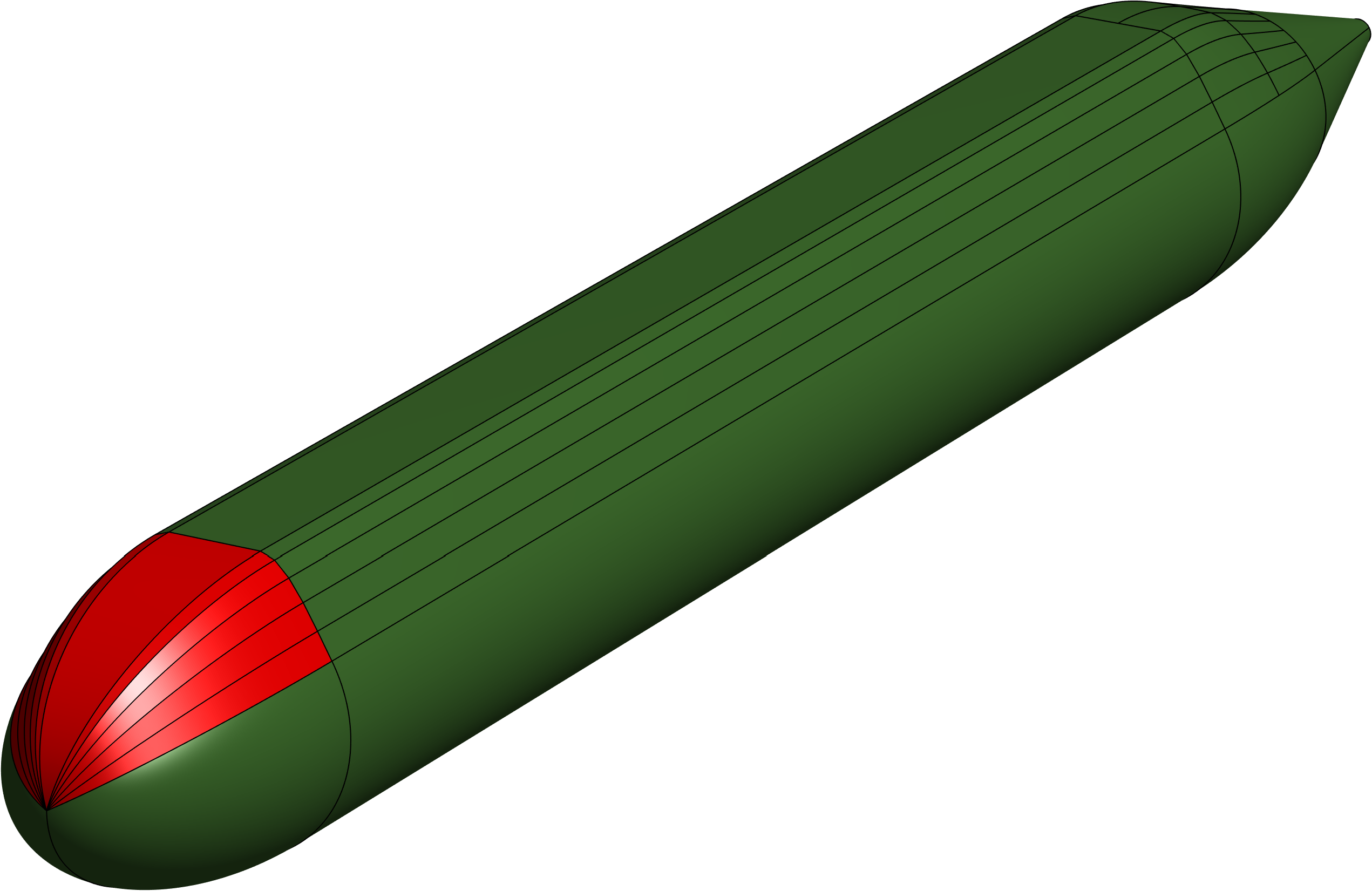}
		\caption{Illustration of the upper bow part.}
		\label{Fig3:bettsi_upperBow}
	\end{subfigure}
	~    
	\begin{subfigure}{0.49\textwidth}
		\centering
		\includegraphics[width=\textwidth]{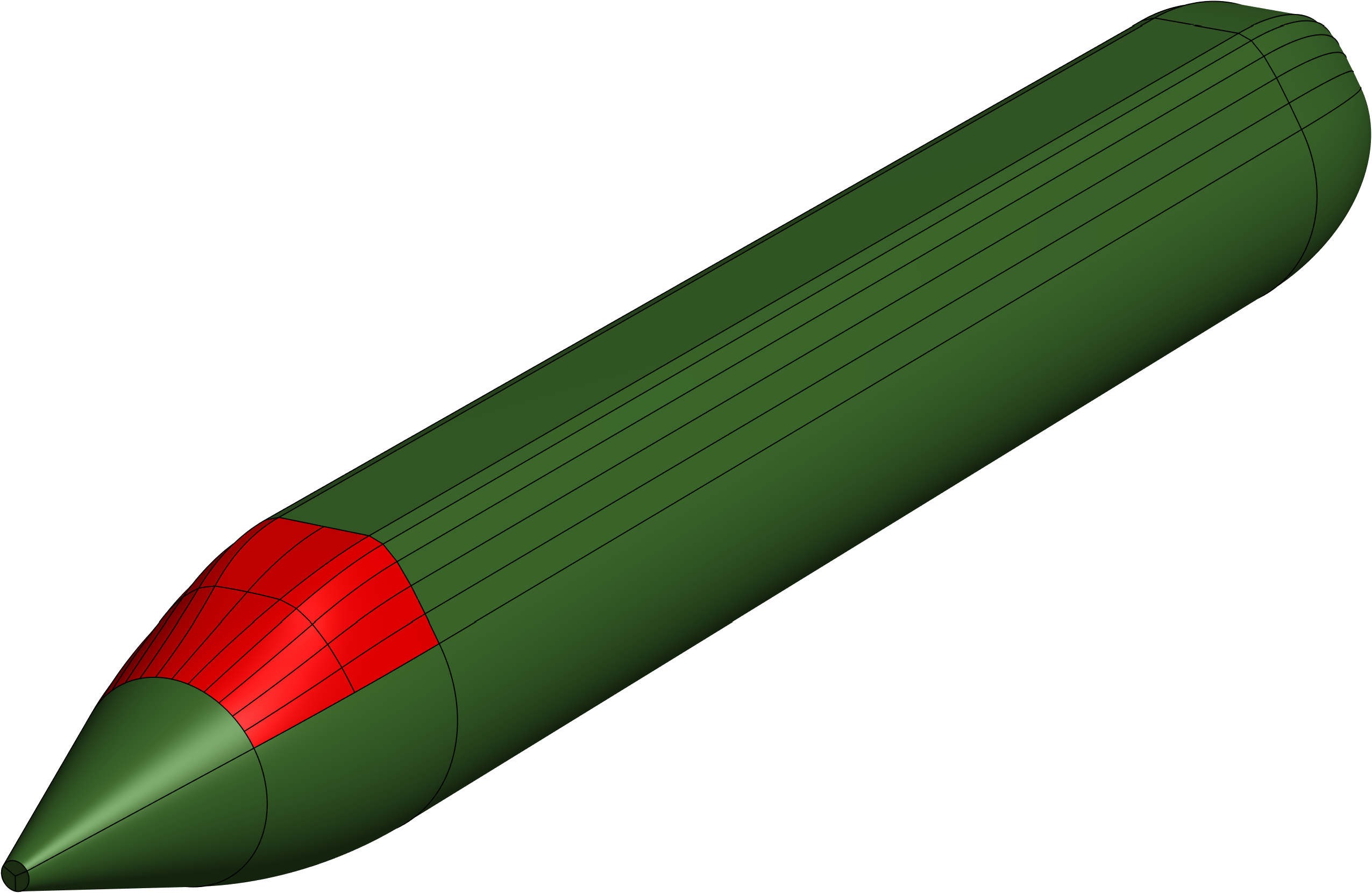}
		\caption{Illustration of the upper transition part.}
		\label{Fig3:bettsi_upperPartOfTailSection}
	\end{subfigure}
	\caption{Main body of BeTSSi submarine.}
\end{figure}

The model is symmetric about the $xz$-plane and has rotational symmetry for the lower part as described in \Cref{Fig3:bettsi_bottom}.
The transition from this axisymmetric part to the deck is described in \Cref{Fig3:bettsi_top}. This transition as well as the deck itself, contains a set of rectangular panels of length $L$.
The cubic polynomial $P_{\mathrm{p}}(y)$, is uniquely defined by the requirement that it defines a smooth transition between the hull and the deck. More precisely, the following requirement must be satisfied: 
\begin{alignat*}{3}
	P_{\mathrm{p}}(s) &= c,\quad  &&P_{\mathrm{p}}\left(b\sin\frac{\beta}{2}\right) = -b\cos\frac{\beta}{2}\\
	P_{\mathrm{p}}'(s) &= 0,\quad &&P_{\mathrm{p}}'\left(b\sin\frac{\beta}{2}\right) = \tan\frac{\beta}{2}
\end{alignat*}
which gives the polynomial
\begin{equation*}
	P_{\mathrm{p}}(y) = c+C_1(y-s)^2+C_2(y-s)^3
\end{equation*}
where
\begin{equation*}
	C_1 = -\frac{3C_4+C_3\tan\frac{\beta}{2}}{C_3^2}, \quad
	C_2 = \frac{2C_4+C_3\tan\frac{\beta}{2}}{C_3^3}, \quad
	C_3 = b\sin\frac{\beta}{2}-s, \quad
	C_4 = c+b\cos\frac{\beta}{2}.
\end{equation*}
The upper part of the bow (highlighted in \Cref{Fig3:bettsi_upperBow}) is obtained by linear lofting of elliptic curves from the 12 points described in \Cref{Fig3:bettsi_top} to the tip of the bow. 

\begin{figure}
	\centering    
	\begin{subfigure}{0.49\textwidth}
		\centering
		\includegraphics[width=0.9\textwidth]{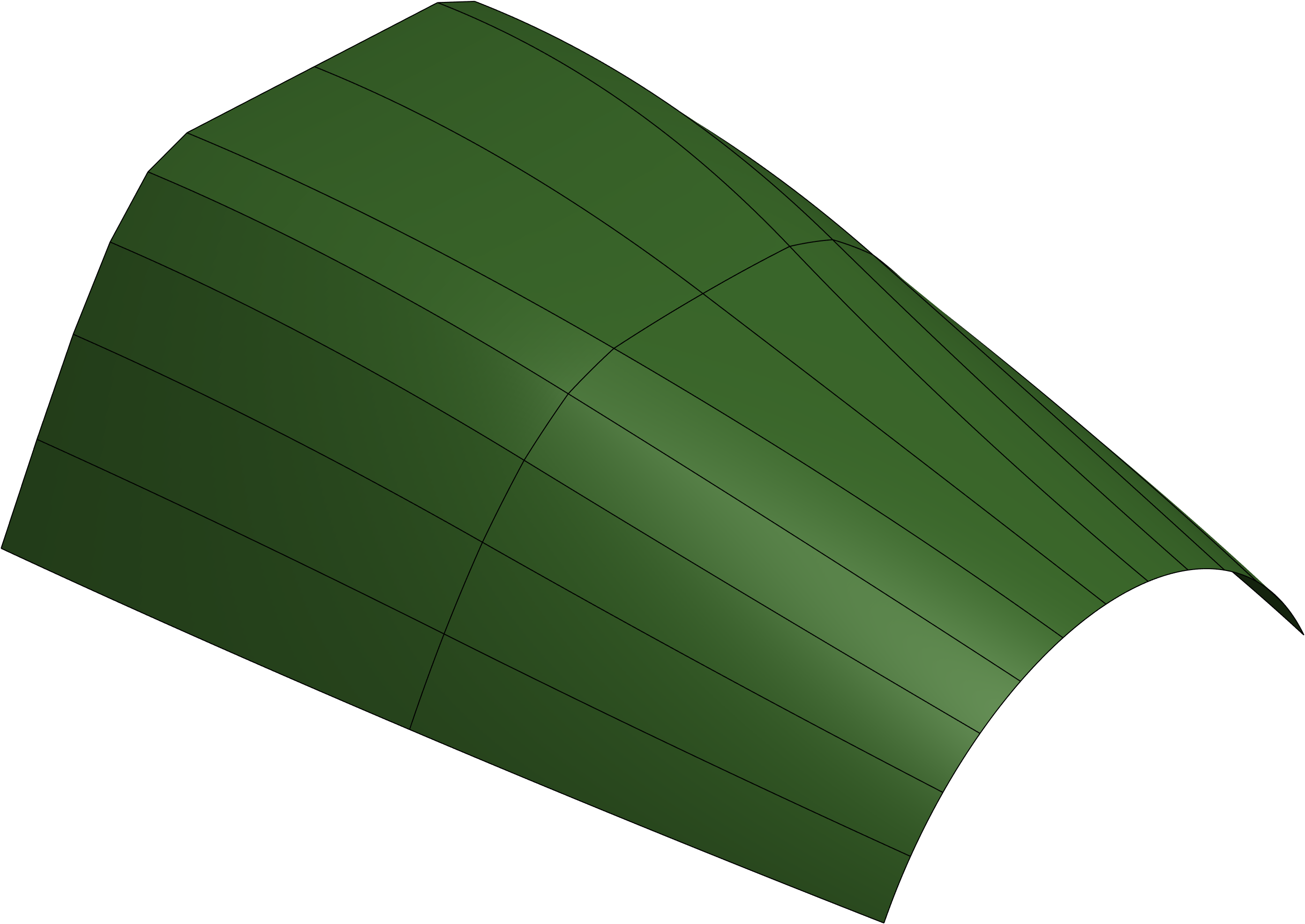}
		\caption{Illustration of the mesh.}
		\label{Fig3:BeTSSi_BC_tailSection}
	\end{subfigure}
	~    
	\begin{subfigure}{0.49\textwidth}
		\centering
		\includegraphics{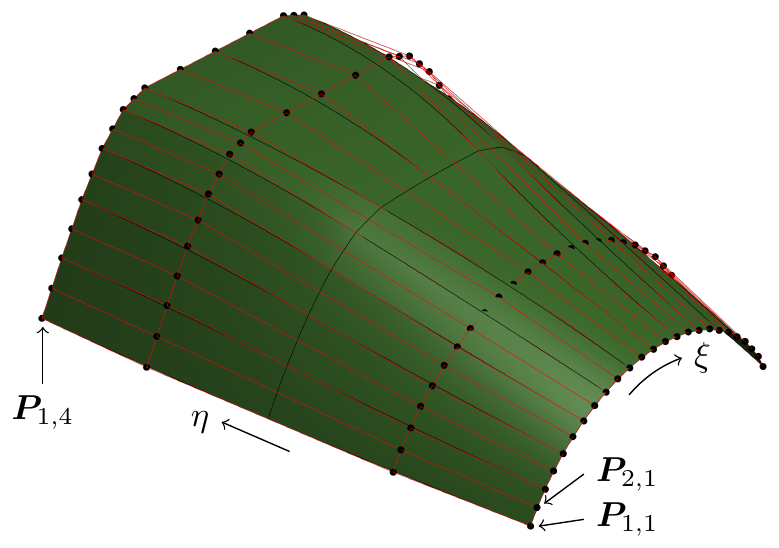}
		\caption{Illustration of the control polygon mesh.}
		\label{Fig3:BeTSSi_BC_tailSection_cp}
	\end{subfigure}
	\caption{Illustration of the upper transition part of the tail.}
\end{figure}

\begin{figure}
	\centering    
	\begin{subfigure}{0.49\textwidth}
		\centering
		\includegraphics{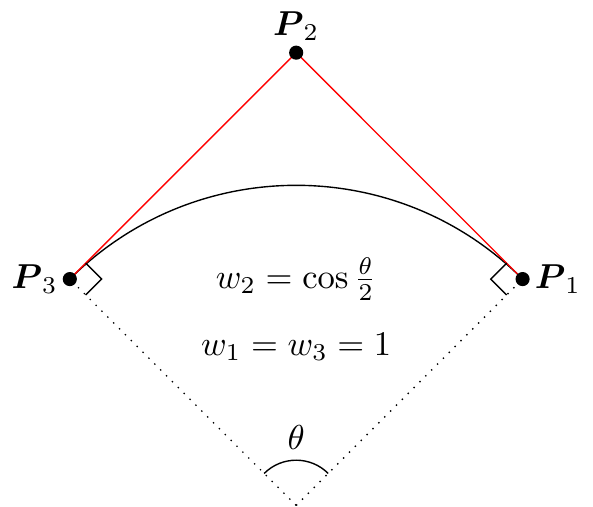}
		\caption{NURBS parametrization of arc of angle $\theta$ using three control points $\{\vec{P}_i\}_{i=1}^3$, the weights $\{w_i\}_{i=1}^3$ and the open knot vector $\Xi=\{0,0,0,1,1,1\}$.}
		\label{Fig3:arcParam1}
	\end{subfigure}
	~    
	\begin{subfigure}{0.49\textwidth}
		\centering
		\includegraphics{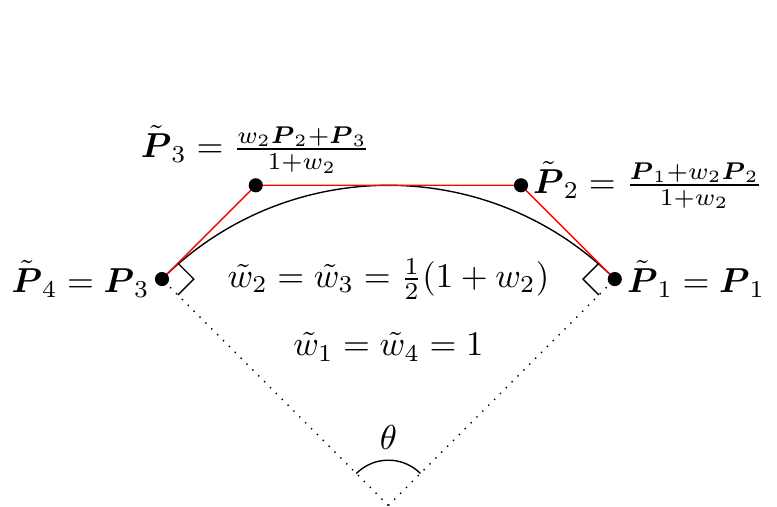}
		\caption{NURBS parametrization of arc of angle $\theta$ using four control points $\{\tilde{\vec{P}}_i\}_{i=1}^4$, the weights $\{\tilde{w}_i\}_{i=1}^4$ and the open knot vector $\tilde{\vec{t}}_\upxi=\{0,0,0,0.5,1,1,1\}$.}
		\label{Fig3:arcParam2}
	\end{subfigure}
	\caption{Two ways of parametrizing an arc using NURBS~\cite[p. 315]{Piegl1997tnb}.}
\end{figure}

The upper part of the tail section (highlighted in \Cref{Fig3:bettsi_upperPartOfTailSection}) is connected using a tensor NURBS surface of degree 2 such that it defines a smooth transition from the axisymmetric cone to the deck. More precisely, the upper part of the cone tail is divided into 12 arcs with angle $\frac{2\PI-\beta}{12}$, and the resulting points are connected to corresponding points on the transition to the deck from the axisymmetric hull. 
As illustrated in \Cref{Fig3:BeTSSi_BC_tailSection}, the NURBS patch is given by 24 elements. Thus, $4\cdot 25 = 100$ control points, $\vec{P}_{i,j}$, are needed as shown in \Cref{Fig3:BeTSSi_BC_tailSection_cp} (25 and 4 control points in the $\xi$ direction and $\eta$ direction, respectively).
The control points $\vec{P}_{1,j}$ and $\vec{P}_{25,j}$ for $j=1,2,3,4$ must be defined as in \Cref{Fig3:arcParam2}, while the control points $\vec{P}_{i,1}$ must be defined as in \Cref{Fig3:arcParam1}. The weights are defined by
\begin{equation*}
	w_{i,j} = \begin{cases}
		\tilde{w}_j & i\,\, \text{odd}\\
		\frac{\tilde{w}_j}{3}\left[(4-j)\cos\left(\frac{2\PI-\beta}{24}\right)+j-1\right] & i\,\, \text{even}
		\end{cases}		
\end{equation*}
where
\begin{equation*}
	\tilde{w}_j = \begin{cases}
		1 & j = 1,4\\
		\frac{1}{2}\left(1+\cos\frac{\alpha}{2}\right) & j = 2,3.
	\end{cases}
\end{equation*}
The locations of the control points $\vec{P}_{i,j}$, $j=2,3$ and $2\leq i\leq 24$, are determined by the requirement that the $x$ component is the same as $\vec{P}_{1,j}$ and the fact that the control polygon lines must be tangential to the surface both at the deck and the cone tail.

\subsection{NACA profiles}
The sail and the rudders are based on the NACA 00xx profiles \cite{Ladson1996cpt,Cummings2015gfa} (the first two digits indicate a symmetric airfoil, and the second two, the thickness-chord ratio). The NACA profiles are all based on the function
\begin{equation}\label{Eq3:NACA}
	f_t(x) = 5t\left(a_0\sqrt{x} + a_1x +a_2x^2+a_3x^3+a_4x^4\right).
\end{equation}
This function satisfies the condition $f_t(0)=0$ and should in addition satisfy
\begin{equation}\label{Eq3:NACAmainConditions}
	f_t(0.3) = \frac{t}{2}, \quad f_t'(0.3) = 0.
\end{equation}
In \cite{Ladson1996cpt,Cummings2015gfa} the coefficients are computed to be
\begin{align*}
	a_0 &= 0.2969\\
	a_1 &=-0.1260\\
	a_2 &=-0.3516\\
	a_3 &= 0.2843\\
	a_4 &=-0.1015.
\end{align*}
The conditions in \Cref{Eq3:NACAmainConditions} are approximated with a residual error of 0.0029\% and 0.013\%, respectively. Moreover, the additional condition $f_t(1) = 0.002$ is satisfied with a  residual error of 0.01\%. In order to have a zero-thickness trailing edge, i.e. $f_t(1)=0$, the original BeTSSi coefficients slightly modify the NACA coefficients to be
\begin{align*}
	a_0 &= 0.2969\\
	a_1 &=-0.1267\\
	a_2 &=-0.3523\\
	a_3 &= 0.2843\\
	a_4 &=-0.1022.
\end{align*}
The conditions in \Cref{Eq3:NACAmainConditions} are here approximated with a residual error of 0.025\% and 0.013\%, respectively. The fact that the conditions in \Cref{Eq3:NACAmainConditions} are approximated so poorly is problematic for an analysis suitable BeTSSi submarine as this results in tangential curves missing the NACA profiles with a significant error, resulting in elements with high aspect ratio or a redundant amount of elements in order to resolve these areas. This fact motivates a more precise definition of these coefficients.
%

Note that the leading-edge radius is given by
\begin{equation*}
	R_{\mathrm{le}} =\lim_{x\to 0^+} \left|\frac{\left[1+f_t'(x)^2\right]^{3/2}}{f_t''(x)}\right| =  \frac{25}{2}a_0^2t^2
\end{equation*}
and the included angle of the trailing edge by
\begin{equation*}
	\delta_{\mathrm{te}} = 2\tan^{-1}|f_t'(1)|.
\end{equation*}
Alternative conditions \cite{Cummings2015gfa}
\begin{equation}\label{Eq3:NACAconditions2}
	R_{\mathrm{le}} = \frac{25}{2}0.2969^2t^2,\quad \delta_{\mathrm{te}} = 2\tan^{-1}(5t\cdot 0.23385)
\end{equation}
yield the coefficients (for usage in double precision)
\begin{align*}
	a_0 &= 0.2969\\
	a_1 &\approx-0.128361732706295\\
	a_2 &\approx-0.335670924960620\\
	a_3 &\approx 0.251127048040123\\
	a_4 &\approx-0.083994390373209.
\end{align*}
Using
\begin{equation*}
	\delta_{\mathrm{te}} = 2\tan^{-1}(5t\cdot 0.243895)
\end{equation*}
yields coefficients slightly closer to the original BeTSSi coefficients.

\begin{figure}
	\centering
	\includegraphics{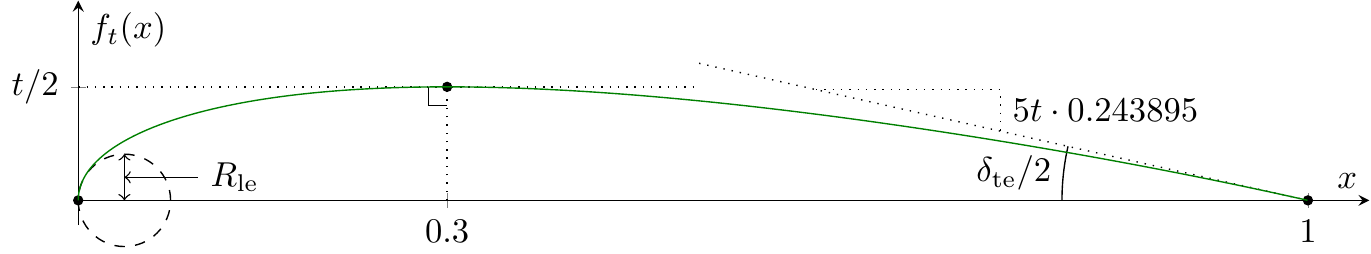}
	\caption{Illustration of the NACA profile used for the sail and the rudders. The five coefficients $a_i$ in \Cref{Eq3:NACA} are restricted by the conditions in \Cref{Eq3:NACAconditions3} as illustrated here.}
	\label{Fig3:NACA2}
\end{figure}

In summary, we shall use the conditions
\begin{equation}\label{Eq3:NACAconditions3}
	f_t(1) = 0,\quad  f_t(0.3) = \frac{t}{2}, \quad f_t'(0.3) = 0,\quad a_0 = 0.2969,\quad f_t'(1) = -5 t\cdot 0.243895
\end{equation}
which are illustrated in \Cref{Fig3:NACA2} and yields the coefficients (in double precision)
\begin{align*}
	a_0 &= 0.2969\\
	a_1 &\approx-0.12651673270629464\\
	a_2 &\approx-0.34981592496061949\\
	a_3 &\approx 0.28392704804012290\\
	a_4 &\approx-0.10449439037320877.
\end{align*}

Computing the relative error in the $L^2$-norm of the NACA profile based on these coefficients and the original NACA profile for the BeTSSi submarine yields an error of about $0.54\%$. Note that $f_t(\xi^2)$ is a polynomial of degree 8, such that the NACA profile can be exactly represented by a spline curve based on the parametrization $\vec{C}(\xi) = [\xi^2,f_t(\xi^2)]$. 

\subsection{Sail}
\begin{figure}
	\centering
	\includegraphics{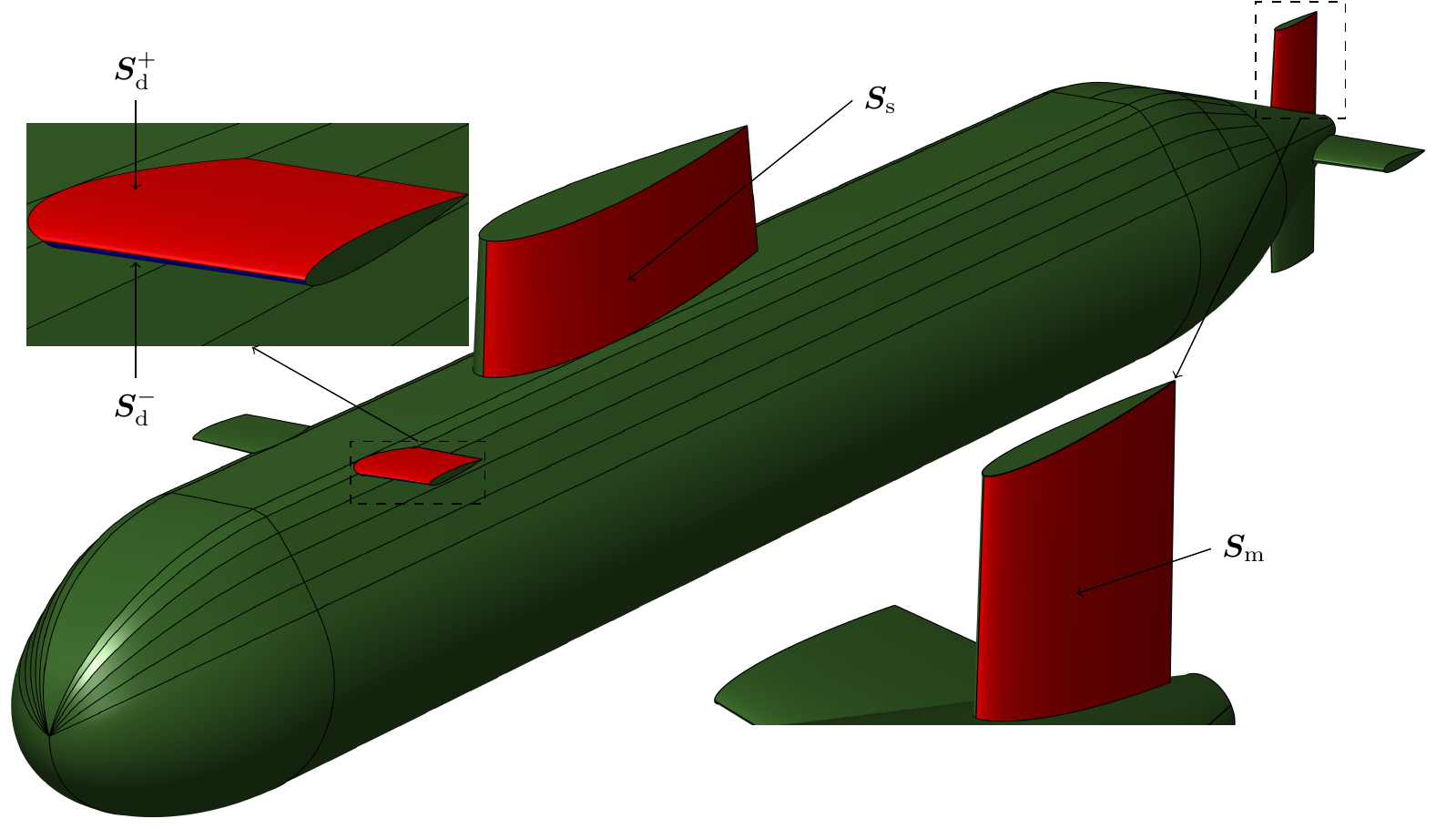}
	\caption{Illustration of the parametrizations $\vec{S}_{\mathrm{s}}$, $\vec{S}_{\mathrm{m}}$ and $\vec{S}_{\mathrm{d}}^\pm$ for the sail, the main rudders and the depth rudders, respectively.}
	\label{Fig3:rudders}
\end{figure}

Consider the port part ($y\geq 0$) of the sail. It can be parametrized by
\begin{equation}\label{Eq3:sail}
	\vec{S}_{\mathrm{s}}(\xi,\eta) = x_{\mathrm{s}}\vec{e}_{\mathrm{x}} + c\vec{e}_{\mathrm{z}} + \begin{bmatrix}
	-\left[l_{\mathrm{ls}} \xi^2+\eta\left(\delta_{\mathrm{s}}-(l_{\mathrm{ls}}-l_{\mathrm{us}})\xi^2\right)\right]\\
	l_{\mathrm{ls}}f_{t_{\mathrm{ls}}}(\xi^2) + \eta\left[l_{\mathrm{us}}f_{t_{\mathrm{us}}}(\xi^2)-l_{\mathrm{ls}}f_{t_{\mathrm{ls}}}(\xi^2)\right]\\
	\eta h_{\mathrm{s}}
	\end{bmatrix},\quad 0\leq \xi\leq 1,\quad 0\leq \eta \leq 1
\end{equation}
where
\begin{equation*}
	t_{\mathrm{us}}=\frac{b_{\mathrm{us}}}{l_{\mathrm{us}}},\quad t_{\mathrm{ls}}=\frac{b_{\mathrm{ls}}}{l_{\mathrm{ls}}},\quad\text{and}\quad b_{\mathrm{ls}}=2s.
\end{equation*}
This parametrization is illustrated in \Cref{Fig3:rudders}. The starboard part of the sail is obtained by mirroring the port side of the sail about the $xz$-plane. Finally, the roof is obtained by a linear loft between these two surfaces.

\subsection{Main rudders}
Consider the port part ($y\geq 0$) of the upper main rudder. It can be parametrized by
\begin{equation}\label{Eq3:mainRudders}
	\vec{S}_{\mathrm{m}}(\xi,\eta) = x_{\mathrm{m}}\vec{e}_{\mathrm{x}} + \begin{bmatrix}
	-l_{\mathrm{lm}} \xi^2-\delta_{\mathrm{m}}\eta\left(1-\xi^2\right)\\
	l_{\mathrm{lm}}f_{t_{\mathrm{lm}}}(\xi^2) + \eta\left[l_{\mathrm{um}}f_{t_{\mathrm{um}}}(\xi^2)-l_{\mathrm{lm}}f_{t_{\mathrm{lm}}}(\xi^2)\right]\\
	\eta h_{\mathrm{m}}
	\end{bmatrix},\quad 0\leq \xi\leq 1,\quad g(\xi)\leq \eta \leq 1
\end{equation}
where
\begin{equation*}
	\delta_{\mathrm{m}}=l_{\mathrm{lm}}-l_{\mathrm{um}},\quad t_{\mathrm{lm}}=\frac{b_{\mathrm{lm}}}{l_{\mathrm{lm}}},\quad\text{and}\quad t_{\mathrm{um}}=\frac{b_{\mathrm{um}}}{l_{\mathrm{um}}}
\end{equation*}
for a function $g$ (to be determined) representing the intersection between the rudder and the cone. The cone can be represented by
\begin{equation}\label{Eq3:cone}
	y^2+z^2=\left(x-x_{\mathrm{c}}\right)^2\tan^2\alpha,\quad x_{\mathrm{c}} = -(L+g_2+(b-h)\cot\alpha).
\end{equation}
Then, inserting the components of $\vec{S}_{\mathrm{m}}(\xi,\eta)$ in \Cref{Eq3:mainRudders} into \Cref{Eq3:cone} yields an equation in $\xi$ and $\eta$. This equation is quadratic in $\eta$ and has the solution $\eta=g(\xi)$ where
\begin{equation*}
	g(\xi) = \frac{-C_{\mathrm{b}}(\xi) + \sqrt{[C_{\mathrm{b}}(\xi)]^2-4C_{\mathrm{a}}(\xi)C_{\mathrm{c}}(\xi)}}{2C_{\mathrm{a}}(\xi)}
\end{equation*}
and
\begin{align*}
	C_{\mathrm{a}}(\xi) &= \left[l_{\mathrm{um}}f_{t_{\mathrm{um}}}(\xi^2)-l_{\mathrm{lm}}f_{t_{\mathrm{lm}}}(\xi^2)\right]^2+h_{\mathrm{m}}^2-\delta_{\mathrm{m}}^2(1-\xi^2)^2\tan^2\alpha\\
	C_{\mathrm{b}}(\xi) &= 2l_{\mathrm{lm}}f_{t_{\mathrm{lm}}}(\xi^2)\left[l_{\mathrm{um}}f_{t_{\mathrm{um}}}(\xi^2)-l_{\mathrm{lm}}f_{t_{\mathrm{lm}}}(\xi^2)\right] +2\tan^2\alpha\left(x_{\mathrm{m}}-l_{\mathrm{lm}}\xi^2-x_{\mathrm{c}}\right)\delta_{\mathrm{m}}\left(1-\xi^2\right)\\
	C_{\mathrm{c}}(\xi) &= [l_{\mathrm{lm}}f_{t_{\mathrm{lm}}}(\xi^2)]^2 - \tan^2\alpha\left(x_{\mathrm{m}}-l_{\mathrm{lm}}\xi^2-x_{\mathrm{c}}\right)^2.
\end{align*}
The trimming curve is then given by
\begin{equation*}
	\vec{r}_{\mathrm{m}}(\xi) = \vec{S}_{\mathrm{m}}(\xi,g(\xi)).
\end{equation*}
The parametrization $\vec{S}_{\mathrm{m}}$ is illustrated in \Cref{Fig3:rudders}. The starboard side of the upper main rudder is given by mirroring the port side of the main upper rudder about the $xz$-plane, and the top part of the rudder is connected by linear lofting. The other main rudders are obtained by rotations by angles of $\ang{90}$, $\ang{180}$ and $\ang{270}$ around the $x$-axis, respectively. Note that this trimming curve may not be represented exactly by NURBS basis functions, and hence, the BeTSSi submarine cannot be exactly represented by NURBS patches without trimming curves.

\subsection{Depth rudders}
Consider the port depth rudder ($y\geq 0$). The upper ($+$) part and lower ($-$) part can be parametrized by
\begin{equation}\label{Eq3:depthRudders}
	\vec{S}_{\mathrm{d}}^\pm(\xi,\eta) = \begin{bmatrix}
		x_{\mathrm{d}}\\
		s\\
		c-\frac{b_{\mathrm{ld}}}{2}
\end{bmatrix}	 + \begin{bmatrix}
	-l_{\mathrm{ld}} \xi^2-\delta_{\mathrm{d}}\eta\left(1-\xi^2\right)\\
	\eta h_{\mathrm{d}}\\
	\pm l_{\mathrm{ld}}f_{t_{\mathrm{ld}}}(\xi^2) \pm \eta\left[l_{\mathrm{ud}}f_{t_{\mathrm{ud}}}(\xi^2)-l_{\mathrm{ld}}f_{t_{\mathrm{ld}}}(\xi^2)\right]
	\end{bmatrix},\quad 0\leq \xi\leq 1,\quad g^\pm(\xi)\leq \eta \leq 1
\end{equation}
where
\begin{equation*}
	\delta_{\mathrm{d}}=l_{\mathrm{ld}}-l_{\mathrm{ud}},\quad t_{\mathrm{ld}}=\frac{b_{\mathrm{ld}}}{l_{\mathrm{ld}}},\quad t_{\mathrm{ud}}=\frac{b_{\mathrm{ud}}}{l_{\mathrm{ud}}},\quad h_{\mathrm{d}}=b-s\quad\text{and}\quad b_{\mathrm{ld}}=2\left[c-P_{\mathrm{p}}\left(s+\frac{C_3}{5}\right)\right].
\end{equation*}
The two panels to be trimmed by this surface are given by
\begin{equation}\label{Eq3:panels}
	D_1^\pm y + D_2^\pm z = D_3^\pm
\end{equation}
where
\begin{equation*}
	D_1^+ = \frac{b_{\mathrm{ld}}}{2},\quad D_2^+ = \frac{C_3}{5},\quad  D_3^+ =  D_1^+s+D_2^+c
\end{equation*}
and
\begin{equation*}
	D_1^- = c-P_{\mathrm{p}}\left(s+\frac{2C_3}{5}\right)-\frac{b_{\mathrm{ld}}}{2},\quad D_2^- = \frac{C_3}{5},\quad  D_3^- =  D_1^-\left(s+\frac{C_3}{5}\right)+D_2^-\left(c-\frac{b_{\mathrm{ld}}}{2}\right).
\end{equation*}
Then, inserting the components of $\vec{S}_{\mathrm{d}}^\pm(\xi,\eta)$ in \Cref{Eq3:depthRudders} into \Cref{Eq3:panels} yields an equation in $\xi$ and $\eta$. This equation is linear in $\eta$ and has the solution $\eta=g^\pm(\xi)$ where
\begin{equation*}
	g^\pm(\xi) = \frac{D_3^\pm - D_1^\pm s - D_2^\pm\left(c-\frac{b_{\mathrm{ld}}}{2}\pm l_{\mathrm{ld}}f_{t_{\mathrm{ld}}}(\xi^2)\right)}{D_1^\pm h_{\mathrm{d}}\pm D_2^\pm\left[l_{\mathrm{ud}}f_{t_{\mathrm{ud}}}(\xi^2) - l_{\mathrm{ld}}f_{t_{\mathrm{ld}}}(\xi^2)\right]}.
\end{equation*}
The trimming curves are then given by
\begin{equation*}
	\vec{r}_{\mathrm{d}}^\pm(\xi) = \vec{S}_{\mathrm{d}}^\pm(\xi,g^\pm(\xi)).
\end{equation*}
The parametrizations $\vec{S}_{\mathrm{d}}^\pm$ are illustrated in \Cref{Fig3:rudders}. The side part is again obtained by linear lofting. The starboard depth rudder is given by mirroring the port depth rudder about the $xz$-plane.

\section{An analysis suitable BeTSSi submarine}
\label{Sec3:BeTSSi_approximation}
Most of the BeTSSi submarine can be exactly represented by second order NURBS basis functions and will need no approximation for our analysis. The areas around the trimming curves, however, needs special care. Instead of incorporating the trimming curves in the analysis of the BeTSSi submarine, a reparametrization of the problematic areas is considered. This enables the possibility to represent the NACA profile with polynomial orders less than 8, which would otherwise be a rather significant restriction of the computational efficiency. A third reason for reparametrizing the submarine is to obtain an analysis suitable mesh around the non-Lipschitz areas (sides of the sail at the deck and the upper part of the depth rudders). The optimal way of parametrizing this area would be to have the same (we use linear) parametrization for the $x$-component as done in~\cite{Lipton2010roi}.

\begin{figure}
	\centering    
	\includegraphics{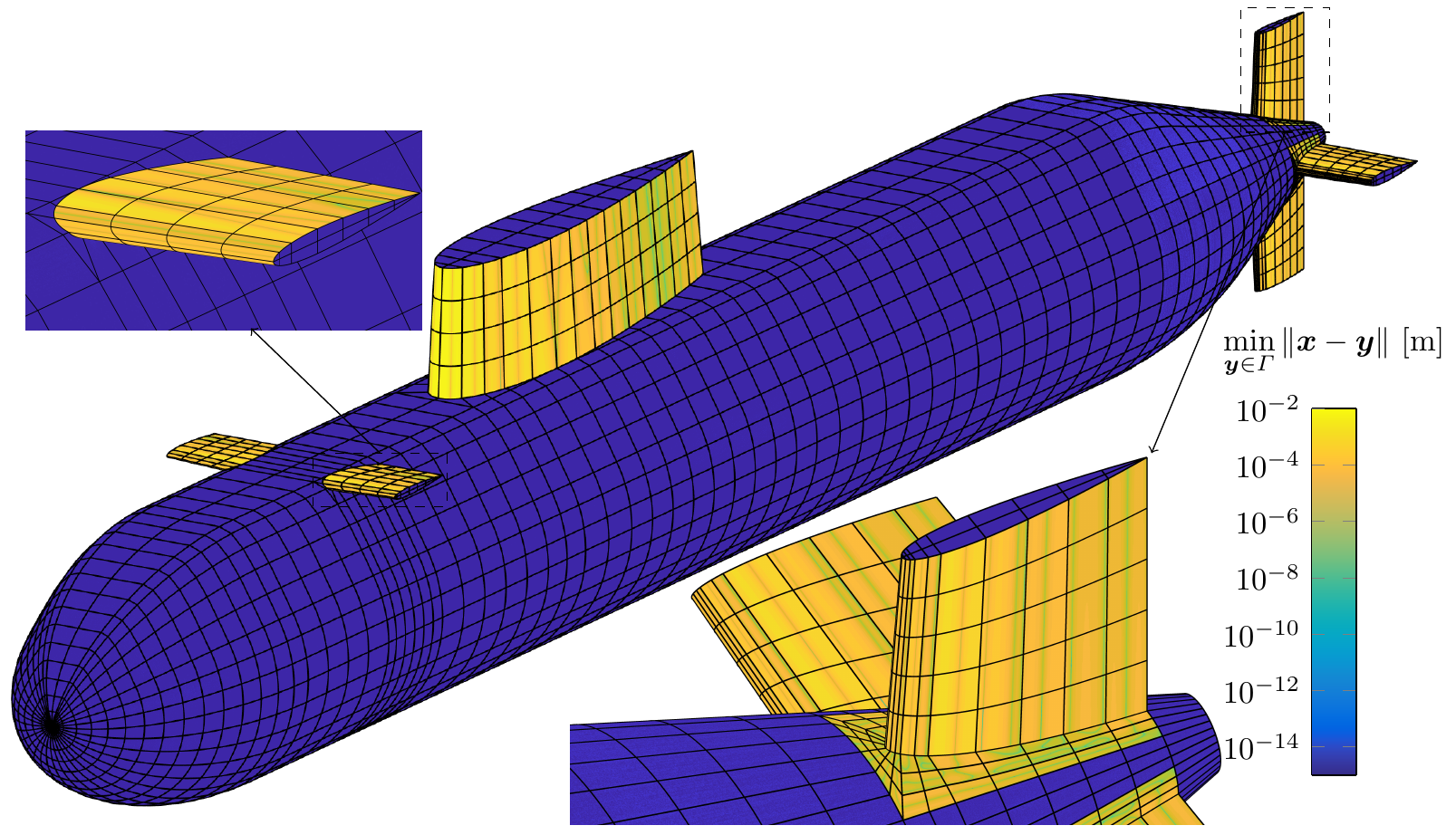}
	\caption{\textbf{An analysis suitable BeTSSi submarine}: The geometric surface approximation $\Gamma_{\check{p}}$ approximates the surface of the exact representation of the BeTSSi submarine $\Gamma$. Surface visualization of the mesh and geometric error for $\check{p}=2$. Most parts of the approximation are exact to machine epsilon precision.}
	\label{Fig3:geometricError}
\end{figure}

\begin{figure}
	\centering    
	\includegraphics{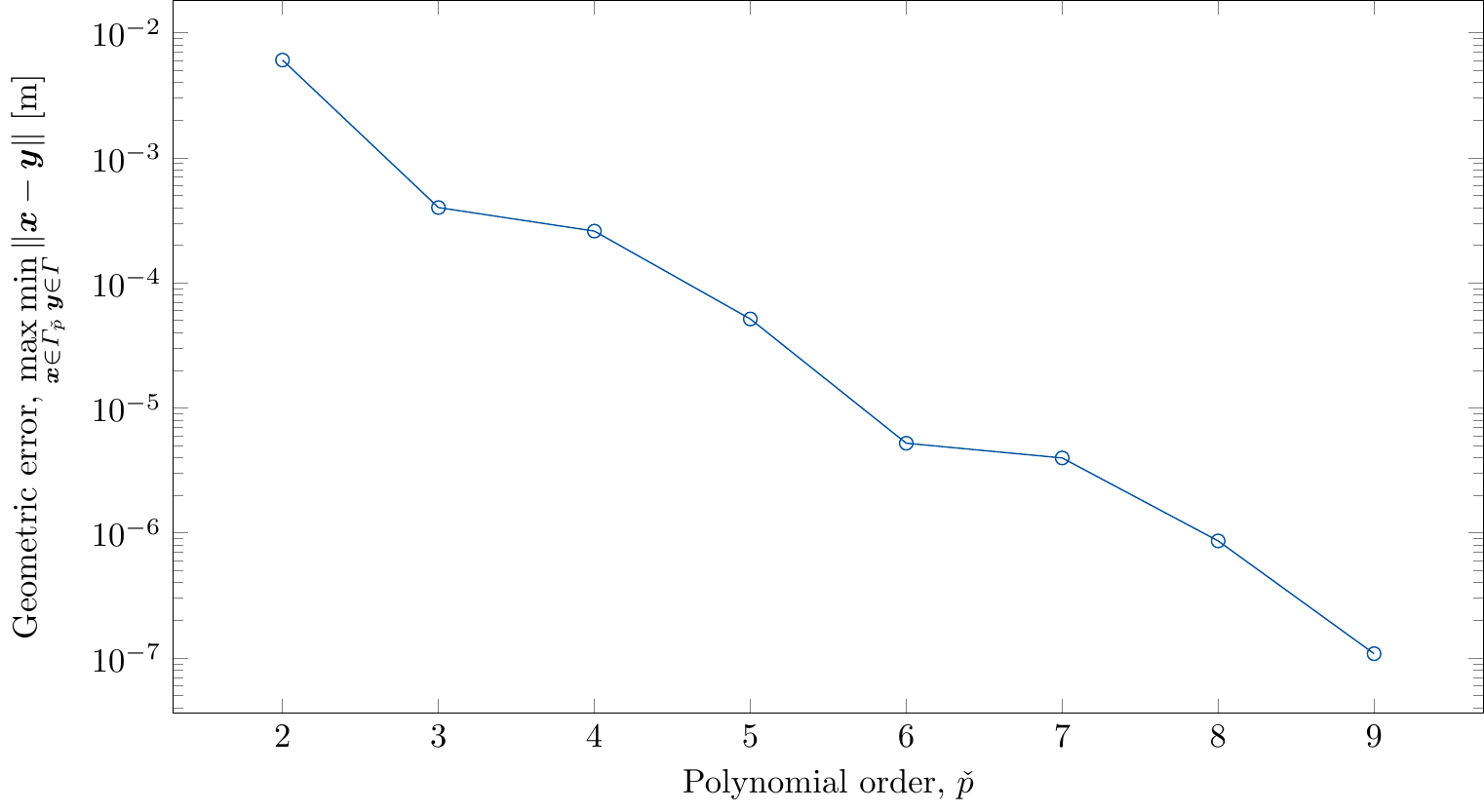}
	\caption{\textbf{An analysis suitable BeTSSi submarine}: The geometric surface approximation $\Gamma_{\check{p}}$ approximates the surface of the exact representation of the BeTSSi submarine $\Gamma$. Convergence plot showing exponential convergence to the exact geometry $\Gamma$.}
	\label{Fig3:geometricErrorsL2_1}
\end{figure}

\begin{figure}
	\centering    
	\includegraphics{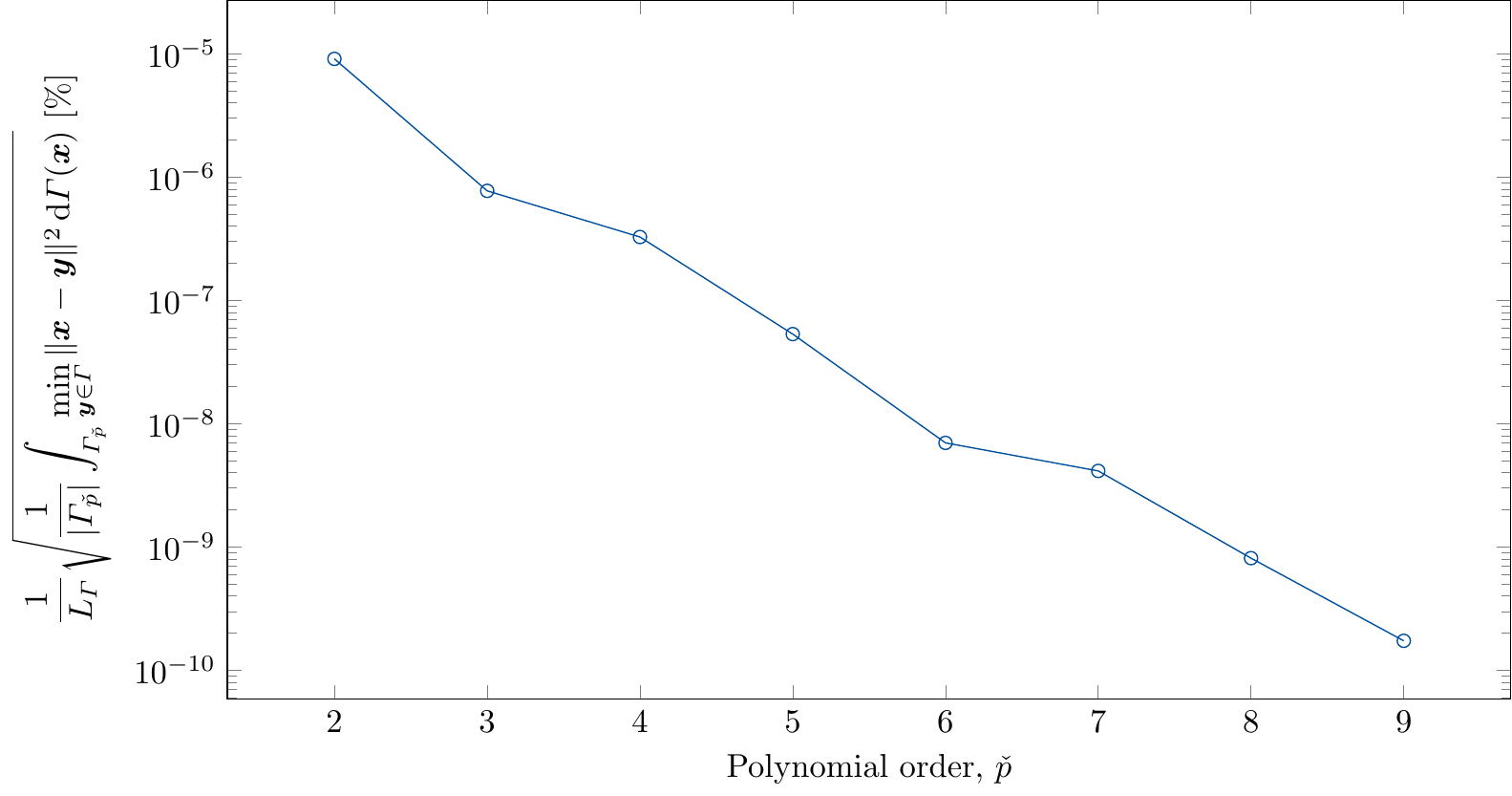}
	\caption{\textbf{An analysis suitable BeTSSi submarine}: Same as \Cref{Fig3:geometricErrorsL2_1} but in another norm. Here, the characteristic length of the geometry is given by $L_\Gamma = a+L+g_2+g_3$.}
	\label{Fig3:geometricErrorsL2_2}
\end{figure}

The approximations are done by performing a least squares of the trimmings curves. For the sail and the depth rudders, the surrounding areas are linear, and can be exactly represented based on the resulting NURBS-curve. For the main rudders, the surrounding areas are approximated by interpolation in such a way that the neighboring (exact) NURBS patches remain unaltered (illustrated in \Cref{Fig3:geometricError}). The interpolation was here preferred above the least squares as it resulted in more analysis suitable basis functions. The upper and lower curves of the sail/rudders are lofted linearly. \Cref{Fig3:geometricErrorsL2_1,Fig3:geometricErrorsL2_2} show the exponential convergence to the exact geometry.

All NURBS patches are conforming such that there is no need to handle master/slave faces by adding constraint equations as described in~\cite[p. 87-91]{Cottrell2006iao}. This results in redundant degrees of freedom, and the optimal mesh certainly requires a solution to this problem. Two very good alternatives include T-splines~\cite{Scott2011tsa} and LR B-splines~\cite{Johannessen2014iau}).

For the sake of brevity, the authors refer to~\cite{Venas2019bts} instead of giving an exact description of every minor detail in constructing this approximation. The exact BeTSSi submarine as well as the approximate submarines for $\check{p}=2,3,4$ are presented in the file formats \texttt{.step}, \texttt{.igs} and \texttt{.3dm} format. 

\begin{figure}
	\centering    
	\includegraphics{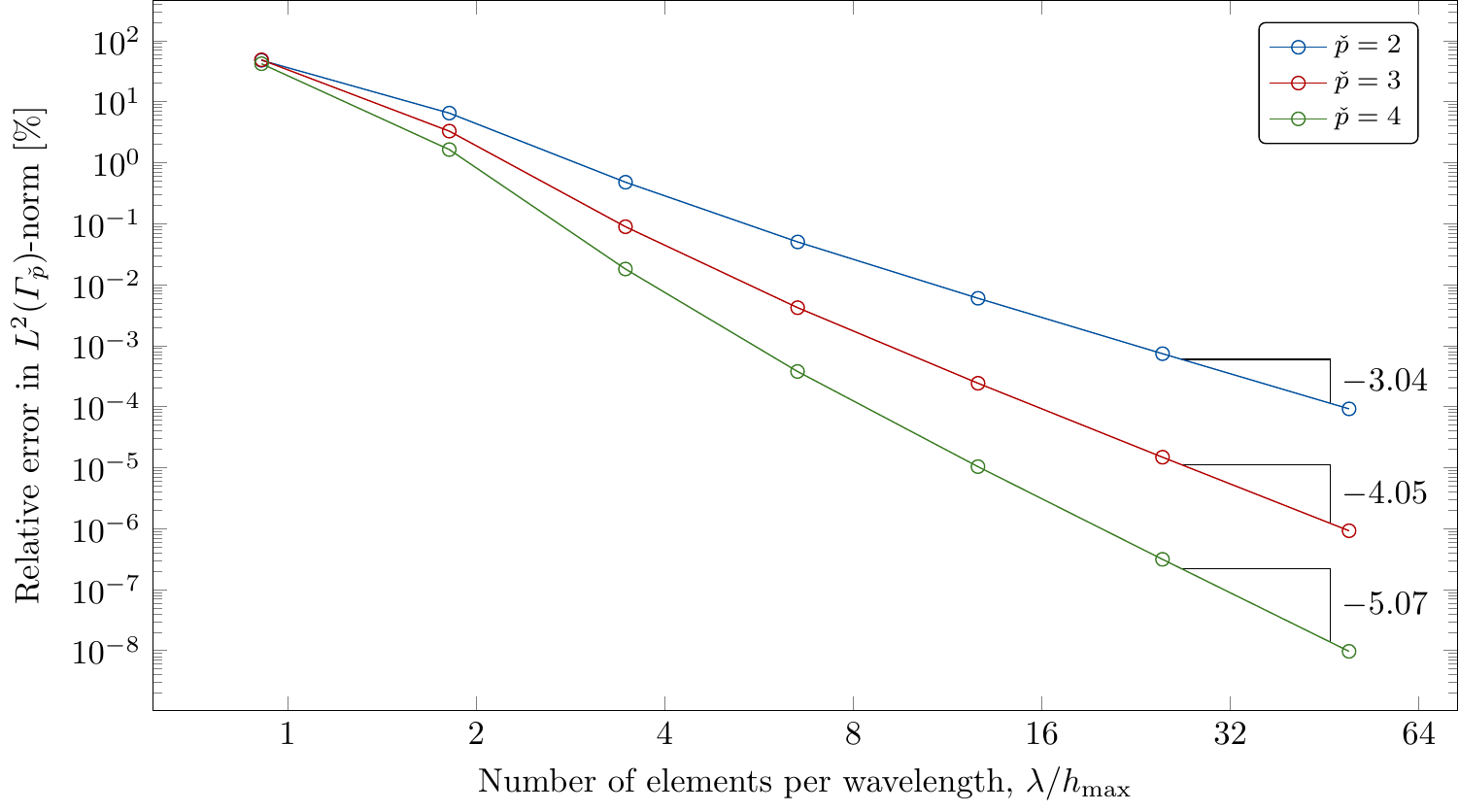}
	\caption{\textbf{An analysis suitable BeTSSi submarine}: Error of the best approximation for the manufactured solution presented in \Cref{Subsec3:manufactured} on $\Gamma_{\check{p}}$.}
	\label{Fig3:BCA_P_BA}
\end{figure}

By considering the manufactured solution in \Cref{Subsec3:manufactured} the numerical evidence observed from \Cref{Fig3:BCA_P_BA} indicates that the presence of non-Lipschitz domain does not affect the convergence rates (also observed in~\cite{Lipton2010roi}).

\section{Triangulation of the BeTSSi submarine}
\label{Sec3:BeTSSi_triangulation}
Triangularized versions of the exact BeTSSi submarine in \texttt{.stl} (both ASCII and binary) and \texttt{.bdf} format can be found in~\cite{Venas2019bts} where the triangulations is an optimization of meshes created in \COMSOL (surface mesh corresponding to the \COMSOL volume meshes considered in this work). An overview of the triangularization meshes can be found in \Cref{Tab3:BeTSSiTriangularization}. Since these meshes are used by WTD in the simulations they have provided for this work, they are denoted by ${\cal M}_m^{\textsc{wtd}}$.
\begin{table}
	\centering
	\caption{\textbf{Triangularization of the BeTSSi submarine}: Data for the meshes.}
	\label{Tab3:BeTSSiTriangularization}
	\begin{tabular}{c S[table-format = 8.0] S[table-format = 8.0] S[table-format = 1.3,round-mode=places,round-precision=3] S[table-format = 1.3,round-mode=places,round-precision=3] S[table-format = 1.4,round-mode=places,round-precision=4] S[table-format = 3.1,round-mode=places,round-precision=1] S[table-format = 3.1,round-mode=places,round-precision=1] S[table-format = 1.4,round-mode=places,round-precision=4]}
		\toprule
		Mesh & {\# triangles} & {\# vertices}  & {$h_{\mathrm{max}}^{(1)}$ [$\si{m}$]}  & {$h_{\mathrm{max}}^{(2)}$ [$\si{m}$]}  & {$\alpha_{\mathrm{min}}$ [$\ang{}$]}   & {$\alpha_{\mathrm{max}}$ [$\ang{}$]} & {$R_{\mathrm{max}}$} & {$S_{\mathrm{min}}$}\\
		\hline
		${\cal M}_{1}^{\textsc{wtd}}$ & 4140 & 2072 & 2.10916 & 1.89333 & 2.01783 & 124.215 & 28.4006 & 0.0336305\\
		${\cal M}_{2}^{\textsc{wtd}}$ & 10406 & 5205 & 1.03183 & 1.00496 & 1.00864 & 124.681 & 56.742 & 0.0168107\\
		${\cal M}_{3}^{\textsc{wtd}}$ & 31104 & 15554 & 0.542723 & 0.498592 & 0.544169 & 124.958 & 105.253 & 0.00906949\\
		${\cal M}_{4}^{\textsc{wtd}}$ & 106888 & 53446 & 0.280763 & 0.256928 & 0.282947 & 124.372 & 202.497 & 0.00471578\\
		${\cal M}_{5}^{\textsc{wtd}}$ & 400886 & 200445 & 0.138583 & 0.130041 & 0.142808 & 121.315 & 401.21 & 0.00238013\\
		${\cal M}_{6}^{\textsc{wtd}}$ & 1584014 & 792009 & 0.0722726 & 0.0691461 & 0.070946 & 124.969 & 807.597 & 0.00118243\\
		\bottomrule
	\end{tabular}
\end{table}
The \textit{resolution} (res) parameter $\lambda/h_{\mathrm{max}}^{(2)}$ (at $f = \SI{1}{kHz}$) is used in the file names. In~\Cref{Tab3:BeTSSiTriangularization}, $h_{\mathrm{max}}^{(1)}$ is defined as the maximum of the diameters of the smallest circle that inscribes the triangular element. For the $i^{\mathrm{th}}$ triangle with side lengths $l_{i,1}$, $l_{i,2}$ and $l_{i,3}$, it is given by
\begin{equation*}
	h_{\mathrm{max}}^{(1)} = \max_i\frac{2l_{i,1}l_{i,2}l_{i,3}}{\sqrt{(l_{i,1}+l_{i,2}+l_{i,3})(l_{i,1}+l_{i,2}-l_{i,3})(l_{i,1}+l_{i,3}-l_{i,2})(l_{i,2}+l_{i,3}-l_{i,1})}}.
\end{equation*}
\COMSOL uses another common definition of the element size, namely the largest side length of the triangle
\begin{equation*}
	h_{\mathrm{max}}^{(2)} = \max_{i,j} l_{i,j}.
\end{equation*}
The three angles of a triangle may be computed by
\begin{equation*}
	\alpha_{i,1} = \cos^{-1}\left(\frac{l_{i,2}^2+l_{i,3}^2-l_{i,1}^2}{2l_{i,2}l_{i,3}}\right),\quad
	\alpha_{i,2} = \cos^{-1}\left(\frac{l_{i,1}^2+l_{i,3}^2-l_{i,2}^2}{2l_{i,1}l_{i,3}}\right)\quad\text{and}\quad
	\alpha_{i,3} = \cos^{-1}\left(\frac{l_{i,1}^2+l_{i,2}^2-l_{i,3}^2}{2l_{i,1}l_{i,2}}\right),
\end{equation*}
such that the maximum and minimum angle are given by
\begin{equation*}
	\alpha_{\mathrm{max}} = \max_{i,j}\alpha_{i,j}\quad\text{and}\quad\alpha_{\mathrm{min}} = \min_{i,j}\alpha_{i,j},
\end{equation*}
respectively. The maximum aspect ratio is defined by
\begin{equation*}
	R_{\mathrm{max}} = \max_i \frac{\max_j l_{i,j}}{\min_j l_{i,j}}
\end{equation*}
and the minimum skewness is defined by
\begin{equation*}
	S_{\mathrm{min}} = \min_{i,j}\left[1-\max\left(\frac{\alpha_{i,j}-\alpha_{\mathrm{e}}}{\ang{180}-\alpha_{\mathrm{e}}}, \frac{\alpha_{\mathrm{e}}-\alpha_{i,j}}{\alpha_{\mathrm{e}}},\right)\right], \quad\alpha_{\mathrm{e}}= \ang{60}.
\end{equation*}
The main take-away here is the inevitability of the increase in the aspect ratio (and the reduction in skewness) during refinement. This is because of the presence of non-Lipschitz domains.

\section*{References}
\bibliographystyle{TK_CM} 
\bibliography{references}

\end{document}